\documentclass[12pt,a4paper,reqno]{amsart}
\usepackage[english]{babel}
\usepackage{amssymb,latexsym,amsfonts,amsthm,upref,amsmath}
\usepackage[foot]{amsaddr}
\usepackage[margin=1in]{geometry} 
\usepackage[dvipsnames,x11names]{xcolor}
 \definecolor{myblue}{HTML}{003399}
\usepackage{hyperref}
\hypersetup{colorlinks,citecolor=myblue,filecolor=black,linkcolor=myblue,urlcolor=myblue}
\usepackage{enumerate}  
\usepackage{tikz}
\usepackage{float}
\usepackage{cleveref}
\usepackage{mathtools}
\usepackage{cite}
\usepackage{filecontents}
\usepackage{enumitem}
\makeatletter
\newcommand{\leqnomode}{\tagsleft@true}
\newcommand{\reqnomode}{\tagsleft@false}
\newcommand{\cev}[1]{\reflectbox{\ensuremath{\vec{\reflectbox{\ensuremath{#1}}}}}}
\makeatother
\newtheorem*{thm*}{Theorem}
\newtheorem*{lem*}{Lemma}
\newtheoremstyle{prim}{}{}{\normalfont}{}{\bfseries}{.}{ }{}
\newtheoremstyle{stil}{}{}{\slshape}{}{\bfseries}{.}{ }{}
\theoremstyle{stil}
\newtheorem{thm}{Theorem}[section]
\newtheoremstyle{defi}{}{}{}{}{\bfseries}{.}{ }{}
\theoremstyle{defi}
\newtheorem{defn}[thm]{Definition}
\theoremstyle{defi}
\newtheorem{rem}[thm]{Remark}
\theoremstyle{stil}
\newtheorem*{mthm*}{Main Theorem}
\newtheorem*{kor*}{Corollary}
\newtheorem{pro}[thm]{Proposition}
\theoremstyle{stil}
\newtheorem{lem}[thm]{Lemma}
\theoremstyle{stil}
\newtheorem{kor}[thm]{Corollary}
\theoremstyle{prim}

\newenvironment{prf}{\noindent \textit{Proof.}}{\null\hfill$\qed$\hskip
2mm\vskip 2mm}

\newcommand{\K}{ {\rm K}}
\newcommand{\I}{ {\rm I}}

\newcommand{\J}{ {\rm J}}
\newcommand{\Jt}{ \widetilde{\rm J}}

\newcommand{\modd}{ \,{\rm mod\,\,}}
\newcommand{\F}{ {\rm F}}

\newcommand{\Ft}{ \widetilde{\rm F}}

\newcommand{\Uu}{ {\rm U} }

\newcommand{\Uh}{ {\rm U}_{\hspace{-1pt}h\hspace{-1pt}} (\widehat{\gl}_N)}
\newcommand{\Uhg}{ {\rm U}_{\hspace{-1pt}h\hspace{-1pt}} (\widehat{\g}_N)}
\newcommand{\Uhsl}{ {\rm U}_{\hspace{-1pt}h\hspace{-1pt}} (\widehat{\sll}_N)}
\newcommand{\Uhcri}{ {\rm U}_{\hspace{-1pt}h\hspace{-1pt}} (\widehat{\gl}_N)_{\text{\rm cri}}}
\newcommand{\Uhplus}{ \wvr{{\rm U}}_{\hspace{-1pt}h\hspace{-1pt}}^{+\hspace{-1pt}} (\widehat{\gl}_N)}
\newcommand{\Uhplusg}{ \wvr{{\rm U}}_{\hspace{-1pt}h\hspace{-1pt}}^{+\hspace{-1pt}} (\widehat{\g}_N)}
\newcommand{\Uhslplus}{ \wvr{{\rm U}}_{\hspace{-1pt}h\hspace{-1pt}}^{+\hspace{-1pt}} (\widehat{\sll}_N)}
\newcommand{\Vccgl}{\mathcal{V}_{\hspace{-1pt}c}(\gl_N)}
\newcommand{\Vccglg}{\mathcal{V}_{\hspace{-1pt}c}(\g_N)}
\newcommand{\Vccgll}{\Vc_{\hspace{-1pt}c}(\gl_N)}
\newcommand{\Vcccri}{\mathcal{V}_{\hspace{-1pt}\text{\rm cri}}(\gl_N)}

\newcommand{\R}{ {\overline{R}}}

\newcommand{\vac}{\mathop{\mathrm{\boldsymbol{1}}}}

\newcommand{\lk}{\lambda}

\newcommand{\gl}{\mathfrak{gl}}
\newcommand{\sll}{\mathfrak{sl}}
\newcommand{\g}{\mathfrak{g}}

\newcommand{\z}{\mathfrak{z}}

\newcommand{\CC}{\mathbb{C}}

\newcommand{\ZZ}{\mathbb{Z}}

\newcommand{\Lc}{\mathcal{L}}

\newcommand{\Sc}{\mathcal{S}}

\newcommand{\Tc}{\mathcal{T}}

\newcommand{\Vc}{\wvr{\mathcal{V}}}

\newcommand{\Vcccrib}{\overline{\mathcal{V}}_{\hspace{-1pt}\text{\rm cri}}(\gl_N)}
\newcommand{\VccNb}{\overline{\mathcal{V}}_{\hspace{-1pt}-N}(\gl_N)}

\newcommand{\wtld}{\widetilde}
\newcommand{\wht}{\widehat}
\newcommand{\wvr}{\overline}

\newcommand{\ot}{\otimes}
\newcommand{\ts}{\hspace{1pt}}
\newcommand{\qdet}{ {\rm qdet}\hspace{1pt}}
\newcommand{\tr}{ {\rm tr}}
\newcommand{\sgn}{ \mathop{\rm sgn}}

\newcommand{\ndo}{\mathop{\mathrm{End}}}
\newcommand{\om}{\mathop{\mathrm{Hom}}}

\newcommand{\rez}{\mathop{\mathrm{Res}}}

\newcommand{\diag}{\mathop{\mathrm{diag}}}
\newcommand{\cdotrl}{\mathop{\hspace{-2pt}\underset{\text{RL}}{\cdot}\hspace{-2pt}}}
\newcommand{\cdotlr}{\mathop{\hspace{-2pt}\underset{\text{LR}}{\cdot}\hspace{-2pt}}}

\newcommand{\iotaop}{\mathop{\iota}}
\newcommand{\iotaopz}{\mathop{\iota_{z}}}
\newcommand{\iotaopx}{\mathop{\iota_{x}}}
\newcommand{\iotaopjd}{\mathop{\iota_{z_1,z_2}}}
\newcommand{\iotaopdj}{\mathop{\iota_{z_2,z_1}}\hspace{-1pt}}
\newcommand{\iotan}{\mathop{\iota_{z_1,\ldots ,z_n}}}
\newcommand{\iotau}{\mathop{\iota_{u}}}
\newcommand{\iotauv}{\mathop{\iota_{u,v}}}
\newcommand{\iotasigma}{\mathop{\iota_{u_{\sigma_1} ,\ldots ,u_{\sigma_n}}}}
\newcommand{\iotazuv}{\mathop{\iota_{z,u_i,v_{j-n}}}}
\newcommand{\iotauiv}{\mathop{\iota_{u_i,v_{j-n}}}}
\newcommand{\iotazui}{\mathop{\iota_{z,u_i}}}
\newcommand{\iotaxuv}{\mathop{\iota_{x,u,v}}}

\newcommand{\iotazuvvh}{\mathop{\iota_{z,u,v,h}}}
\newcommand{\iotaall}{\mathop{\iota_{z_1,z_2,u_1,\ldots,u_n,v_1,\ldots ,v_m}}}
\newcommand{\iotazzo}{\mathop{\iota_{x,x_0}}}
\newcommand{\iotaztri}{\mathop{\iota_{z_0,u_1,u_2}}}
\newcommand{\jota}{\mathop{\iota_{z_2,z_0,u_1,\ldots ,u_n}}}

\newcommand{\fand}{\quad\text{and}\quad}
\newcommand{\Fand}{\qquad\text{and}\qquad}

\newcommand{\non}{\nonumber}
\newcommand{\beq}{\begin{equation}}
\newcommand{\eeq}{\end{equation}}
\newcommand{\ben}{\begin{equation*}}
\newcommand{\een}{\end{equation*}}

\makeatletter
\def\smalloverbrace#1{\mathop{\vbox{\m@th\ialign{##\crcr\noalign{\kern3\p@}%
  \tiny\downbracefill\crcr\noalign{\kern3\p@\nointerlineskip}%
  $\hfil\displaystyle{#1}\hfil$\crcr}}}\limits}
\makeatother

\makeatletter
\def\smallunderbrace#1{\mathop{\vtop{\m@th\ialign{##\crcr
   $\hfil\displaystyle{#1}\hfil$\crcr
   \noalign{\kern3\p@\nointerlineskip}% 
   \tiny\upbracefill\crcr\noalign{\kern3\p@}}}}\limits}
\makeatother

\begin{document}

\title[On the quantum affine vertex algebra associated with trigonometric \texorpdfstring{$R$}{R}-matrix 
]
{On the quantum affine vertex algebra associated with trigonometric \texorpdfstring{$R$}{R}-matrix 
}

\author{Slaven Ko\v{z}i\'{c}} 
\address{Department of Mathematics, Faculty of Science, University of Zagreb, Bijeni\v{c}ka cesta 30, 10\,000 Zagreb, Croatia}
\email{kslaven@math.hr}
\keywords{Quantum affine algebra, Quantum vertex algebra, $\phi$-Coordinated module, Quantum current}
\subjclass[2010]{17B37, 17B69, 81R50}

\begin{abstract}
We apply the theory of $\phi$-coordinated modules, developed by H.-S. Li, to the Etingof--Kazhdan quantum affine vertex algebra   associated with the trigonometric $R$-matrix of type $A$.   
We prove, for  a certain  associate $\phi$ of the one-dimensional additive formal group, 
that any    $\phi$-coordinated  module for the level $c\in\mathbb{C}$ quantum affine vertex algebra
 is naturally equipped with a structure of 
  restricted level $c$ module    for the quantum affine algebra in type $A$ and vice versa. Moreover, we show that any $\phi$-coordinated module is  irreducible with respect to the  action of the quantum affine vertex algebra if and only if it is irreducible with respect to the corresponding action of the quantum affine algebra.
	In the end, we discuss   relation between the centers  of the quantum affine algebra   and the quantum affine vertex algebra.
\end{abstract}

\maketitle

\allowdisplaybreaks
%\begin{samepage}...\end{samepage}

%\tableofcontents
%\newpage

%%%%%%%%%%%%%%%%%%%%%%%%%%%%%%%%%%%
%%%%%%%%%%%%%%%%%%%%%%%%%%%%%%%%%%%
\section*{Introduction}\label{intro}
\numberwithin{equation}{section}
%%%%%%%%%%%%%%%%%%%%%%%%%%%%%%%%%%%
%%%%%%%%%%%%%%%%%%%%%%%%%%%%%%%%%%%

The notion of  {\em vertex algebra}, originally introduced by Borcherds \cite{B}, presents a remarkable connection between  mathematics and theoretical physics.
The  vertex algebra theory led to important breakthroughs in multiple areas such as automorphic forms, finite simple groups and soliton equations; see, e.g., the books by E. Frenkel and Ben-Zvi \cite{FBZ}, I. Frenkel, Lepowsky and Meurman \cite{FLM} and Kac \cite{Kac2}.
Some of the most extensively studied examples of vertex algebras come from the theory of affine Kac--Moody Lie algebras; see the books by Kac \cite{Kac} and Lepowsky and Li \cite{LLi}.
Motivated by a parallel between the development of the   theories of affine Lie algebras and quantum affine algebras,
as well as by further applications   to two-dimensional statistical models and the quantum Yang--Baxter equation, 
I. Frenkel and Jing \cite{FJ} formulated a fundamental problem of 
generalizing the vertex algebra theory to the quantum case. 
 
The notion of {\em quantum vertex algebra}  was introduced by Etingof and Kazhdan \cite{EK}   based on the ideas of E. Frenkel and Reshetikhin \cite{FR2}.
The   examples of quantum vertex algebras were constructed  in \cite{EK} as quantizations of the quasiclassical structure on the  universal affine vertex algebra in type $A$  when the  classical $R$-matrix is of  rational, trigonometric and elliptic  type.
Recently, a structure theory of quantum vertex algebras was   developed by De Sole, Gardini and Kac \cite{DGK} and the Etingof--Kazhdan construction was generalized to the rational $R$-matrix in types $B$, $C$ and $D$ by Butorac, Jing and the author \cite{BJK}. 
On the other hand, several more general related notions, in particular, of   {\em $h $-adic nonlocal vertex algebra} and  of  its   {\em  module},     were introduced and extensively studied by Li \cite{Li}. 
They present  analogues of the corresponding notions, coming   from the Li      nonlocal vertex algebra theory \cite{LiG1} and  the Bakalov--Kac
      field algebra   theory  \cite{BK},  which are defined over the commutative ring $\CC[[h]]$, thus being compatible with  Etingof--Kazhdan's theory.  
Moreover, the  notion of $h$-adic nonlocal vertex algebra  module, which presents   a     generalization of 
vertex algebra module, appears to provide the right setting for the study of representations of       double Yangians and of  Etingof--Kazhdan's quantum  vertex algebras associated with the rational $R$-matrix; see \cite{Li} and \cite{c11} respectively. However, Li's subsequent results \cite{Li1} suggest that the solution of the original Frenkel--Jing  problem of associating    quantum vertex algebras to quantum affine algebras requires a new concept of {\em $\phi$-coordinated module}. Following such an approach,    Li,  Tan and Wang \cite{LTW}  recently established a correspondence between restricted  modules  for the  Ding--Iohara algebra of level $0$ associated with the affine Lie algebra $\wht{\mathfrak{sl}}_2$ \cite{DI}  and $\phi$-coordinated modules for certain quantum vertex algebra.

The definition of  a  $\phi$-coordinated module $W$ for a quantum vertex algebra $V$, as given in \cite{Li1}, is characterized by a certain  deformed version of the weak associativity property. Roughly speaking, it requires that the expressions
$$
\big((z_1-z_2)^p\ts Y_W(u,z_1)Y_W(v,z_2)\big)\big|_{z_1=\phi(z_2,z_0)}  \big. 
\fand (\phi(z_2,z_0) -z_2)^p\ts Y_W\left(Y(u,z_0)v,z_2\right)
$$
coincide for all $u,v\in V$, where $Y(z)$ is the vertex operator map on $V$, $Y_W(z)$ the $\phi$-coordinated module map, $\phi(z_2,z_0)\in\CC((z_2))[[z_0]]$ an  associate of the one-dimensional additive formal group and $p\geqslant 0$ an integer depending on $u,v$. While  setting $\phi(z_2,z_0)=z_2 +z_0$ leads to the usual weak associativity property, a different choice of the associate appears to be required in order  to adapt the theory to quantum affine algebras; see \cite{Li1}.

Let $\g_N=\gl_N ,\sll_N$.
In this paper, we  
consider
the    quantum   affine vertex algebra $\Vc_{\hspace{-1pt}c}(\g_N)$   associated with the trigonometric $R$-matrix, as defined by Etingof and Kazhdan \cite{EK}. 
We should mention that $\Vc_{\hspace{-1pt}c}(\g_N)$ can be also regarded as an associative algebra over  $\CC[[h]]$, which is topologically generated by the coefficients of certain Taylor series   organized into the  matrix   $T^+(u)\in\ndo\CC^N \ot \Vc_{\hspace{-1pt}c}(\g_N) [[u]]$,   subject to certain dual Yangian-type defining relations.
As a quantum vertex algebra, its vertex operator map $Y(z)$
is given in the form of {\em quantum currents} $\Tc (u)$, which go back to Reshetikhin and Semenov-Tian-Shansky \cite{RS}. Furthermore, the $\Sc$-locality property for $\Vc_{\hspace{-1pt}c}(\g_N)$, which is a quantum analogue of the locality in the corresponding affine vertex algebra, comes from the   {\em quantum current commutation relation} which, in this particular setting, can be expressed as
\begin{align}
&\Tc_1(u_1)R_{21} (e^{-u_1 +u_2 -hc} ) \Tc_2(u_2) R_{21}(e^{-u_1 +u_2})^{-1}\non\\
&\qquad\qquad=
R_{12}(e^{-u_2 +u_1})^{-1}\Tc_2(u_2)R_{12}(e^{-u_2+u_1-hc})\Tc_1(u_1),\label{intr1}
\end{align}
where $R(x)=R_{12}(x)$ is the trigonometric $R$-matrix of type $A$.\footnote{We explain the precise meaning of   relations  \eqref{intr1} and \eqref{intr2} in Subsection \ref{subsec021}.}
On the other hand, the original  quantum current commutation relation in \cite{RS} is given in the multiplicative form,
\begin{align}
&\Lc_1(x_1) R_{21}(x_2 e^{-hc}/ x_1) \Lc_2(x_2)  R_{21}(x_2 /x_1)^{-1}\non\\
&\qquad\qquad=
R_{12}(x_1 /x_2)^{-1}\Lc_2(x_2 ) R_{12}(x_1 e^{-hc}/ x_2)\Lc_1(x_1).\label{intr2}
\end{align}
Its significance comes from Ding's quantum current realization of the quantum affine algebra  $\Uhg$ \cite{D}, which relies on the famous Ding--Frenkel isomorphism \cite{DF}. The   algebra generators  are given as coefficients of matrix entries of the {\em quantum current} $\Lc(x)$, so   that
$
\Lc(x)$ belongs to  $ \ndo\CC^N\ot\Uhg [[x^{\pm 1}]]
$,
while the defining relations at the level $c\in\CC$ are given by   \eqref{intr2}, along with one more family of relations in the $\g_N= \sll_N$ case.
As in \cite[Sect. 5]{Li1}, in this paper we  consider   the  $\phi$-coordinated $\Vc_{\hspace{-1pt}c}(\g_N)$-modules for the associate $\phi(z_2,z_0)=z_2 e^{z_0}$, which connects commutation relations \eqref{intr1} and \eqref{intr2}. More specifically, by applying the substitutions $x_i=ze^{u_i}$ with $i=1,2$, multiplicative relation \eqref{intr2} takes the additive form  as in \eqref{intr1}. It is worth noting that both  additive and   multiplicative forms of the trigonometric $R$-matrix naturally occur in the theories of quantum groups and exactly solvable models; see \cite{PS,FRT,J}.

As with the rational $R$-matrix case \cite{c11}, the multiple copies of quantum currents $\Lc(x_i)$ with $i=1,\ldots ,n$ can be organized into the operators $\Lc_{[n]}(x_1,\ldots,x_n)$ in the variables $x_1,\ldots ,x_n$ which satisfy certain   generalized version of commutation relation \eqref{intr2}.  
Roughly speaking, such operators take place of the normal-ordered products of $n$ quantum currents.
 In particular, for any {\em restricted} $\Uhg$-module, i.e. for any module  $W$ such that $\Lc(x)w$ belongs to $\ndo\CC^N\ot W ((x))[[h]]$ for all $w\in W$, the series $\Lc_{[n]}(x_1,\ldots,x_n)w$ possesses only finitely many negative powers of the variables $x_1,\ldots ,x_n$ modulo $h^k$ for all $k\geqslant 0$ and $w\in W$.
By combining Ding's quantum current realization    \cite{D} with Li's theory of $\phi$-coordinated modules \cite{Li1} and Cherednik's fusion procedure for the trigonometric $R$-matrix \cite{C} in the $\g_N=\sll_N$ case, 
we
 establish the following correspondence between restricted modules for the quantum affine algebra and $\phi$-coordinated modules for the Etingof--Kazhdan quantum   vertex algebra, which is the main result of this paper.

\begin{mthm*}\label{mainthm1}
Let $\g_N=\gl_N ,\sll_N$.
Let $W$ be a restricted $\Uhg$-module of level $c\in\CC$. There exists a unique structure of $\phi$-coordinated $\Vc_{\hspace{-1pt}c}(\g_N)$-module on $W$, where $\phi(z_2,z_0) = z_2 e^{z_0}$, such that
\beq\label{formula}
Y_W(T_{[n]}^+(u_1,\ldots ,u_n)\vac,z) = \Lc_{[n]}(x_1,\ldots ,x_n)\big|_{x_1 = ze^{u_1},\ldots, x_n = ze^{u_n}}\big. \quad\text{for all }n\geqslant 1.
\eeq
Conversely, let $(W,Y_W)$ be a $\phi$-coordinated $\Vc_{\hspace{-1pt}c}(\g_N)$-module, where $\phi(z_2,z_0) = z_2 e^{z_0}$. There exists a unique structure of  restricted  $\Uhg$-module of level $c$ on $W$ such that
\beq\label{moduleformula}
\Lc (z) = Y_W (T^+(0)\hspace{-1pt}\vac, z).
\eeq
Moreover, a topologically free  $\CC[[h]]$-submodule $W_1$ of $W$ is a $\phi$-coordinated $\Vc_{\hspace{-1pt}c} (\g_N)$-submodule of $W$ if and only if $W_1$ is an $\Uhg$-submodule of $W$.
\end{mthm*}

In order to establish this correspondence,  some minor modifications had to be made to the definitions of quantum affine algebra and $\phi$-coordinated module. More specifically,   both notions were  redefined over  the  ring $\CC[[h]]$ and suitably completed, so that they are compatible 
with   Etingof--Kazhdan's definition of quantum vertex algebra.  

In the end, we recollect that the universal affine vertex algebra,  which governs the representation theory of the corresponding affine Lie algebra $\wht{\g}_N$, is constructed on the vacuum module over the universal enveloping algebra $\Uu(\wht{\g}_N)$; see \cite{FZ,Lian}. In contrast, $\Vc_{\hspace{-1pt}c} (\g_N)$ is not the vacuum module over  $\Uhg$, although its quantum vertex algebra structure turns into the corresponding affine vertex algebra in the classical limit. Furthermore, it is not clear whether the {\em vacuum module  $\Vccglg$   at the level $c$} over the quantum affine algebra  $\Uhg$ possesses any natural quantum vertex algebra-like structure that governs the representation theory of $\Uhg$.
However, we have the following simple consequence of the \hyperref[mainthm1]{Main Theorem}:

\begin{kor}\label{maincor}
Let $\g_N=\gl_N ,\sll_N$.
The vacuum module $\Vccglg$ over the quantum affine algebra $\Uhg$ is a $\phi$-coordinated $\Vc_{\hspace{-1pt}c}(\g_N)$-module. Moreover, $\Vccglg$ is an irreducible $\Uhg$-module if and only if it is an irreducible $\phi$-coordinated $\Vc_{\hspace{-1pt}c} (\g_N)$-module.
\end{kor}

The paper is organized as follows.  In Sections \ref{sec02} and \ref{sec0102}, we introduce the notation and provide  preliminary definitions and results on restricted $\Uhg$-modules   and on $\phi$-coordinated  $\Vc_{\hspace{-1pt}c} (\g_N)$-modules respectively.
In Section \ref{sec05}, we prove  the \hyperref[mainthm1]{Main Theorem}.
Finally, in Section \ref{newsec02}, we  discuss a connection between the families of central elements of the quantum affine algebra and the quantum affine vertex algebra established by $\phi$-coordinated module map \eqref{formula}.

%%%%%%%%%%%%%%%%%%%%%%%%%%%%%%%%%%%
%%%%%%%%%%%%%%%%%%%%%%%%%%%%%%%%%%%
\section{Restricted modules for the quantum affine algebra}\label{sec02}
%
%\setcounter{equation}{0}
%%%%%%%%%%%%%%%%%%%%%%%%%%%%%%%%%%%
%%%%%%%%%%%%%%%%%%%%%%%%%%%%%%%%%%%

In this section, we recall some basic properties  of the trigonometric $R$-matrix of type $A$.
Next, we recall Ding's quantum current realization of the quantum affine algebra in type $A$ and the corresponding notion of restricted module. Also, we derive certain   properties of the quantum currents which are required in the following sections. Finally, we demonstrate how the \hyperref[mainthm1]{Main Theorem} implies Corollary \ref{maincor}.

%%%%%%%%%%%%%%%%%%%%%%%%%%%%%%%%%%%
%%%%%%%%%%%%%%%%%%%%%%%%%%%%%%%%%%%
\subsection{Trigonometric \texorpdfstring{$R$}{R}-matrix}\label{sec01}
%%%%%%%%%%%%%%%%%%%%%%%%%%%%%%%%%%%
%%%%%%%%%%%%%%%%%%%%%%%%%%%%%%%%%%%

We  use the standard tensor notation, i.e. for any  
$$A=\sum_{i,j,k,l=1}^N a_{ijkl} \ts e_{ij}\ot e_{kl}\,\in\, \ndo\CC^N\ot\ndo\CC^N$$
and indices $r,s=1,\ldots , m$  such that $r\neq s$, where $m\geqslant 2$ and $e_{ij}\in\ndo \CC^N$ are the matrix units, we denote by $A_{rs}$ the element of the algebra $(\ndo\CC^N)^{\ot m}$, 
\beq\label{notation}
A_{rs}=\sum_{i,j,k,l=1}^N a_{ijkl}\ts (e_{ij})_r (e_{kl})_s,\quad\text{where}\quad
(e_{ij})_p = 1^{\ot (p-1)} \ot e_{ij} \ot 1^{\ot{(m-p)}}.
\eeq

Let $N\geqslant 2$ be an integer and $h$ a formal parameter. Introduce the  trigonometric $R$-matrix of type $A$ by 
\begin{align}
\R(x) =&\sum_{i=1}^N e_{ii}\ot e_{ii} 
+ e^{-h/2} \frac{1-x}{1-e^{-h}x}\sum_{\substack{i,j=1\\i\neq j}}^N e_{ii}\ot e_{jj}\non\\
&+\frac{\left(1-e^{-h}\right)x}{1-e^{-h}x}\sum_{\substack{i,j=1\\i> j}}^N e_{ij}\ot e_{ji}
+\frac{1-e^{-h}}{1-e^{-h}x}\sum_{\substack{i,j=1\\i< j}}^N e_{ij}\ot e_{ji}.\label{Rbar}
\end{align}
$R$-matrix \eqref{Rbar} can be regarded as a rational function in the variables $x$ and $e^{h/2}$, i.e. as an element of  $(\ndo\CC^N )^{\ot 2} (x,e^{h/2})$. 
 It satisfies the {\em Yang--Baxter equation}
\beq\label{YBE}
\R_{12}(x/y)\R_{13}(x)\R_{23}(y)=\R_{23}(y)\R_{13}(x)\R_{12}(x/y) 
\eeq 
and it possesses the {\em unitarity property}
\beq\label{unitrig}
\R_{12}(x)\ts \R_{21}(1/x) =1,
\eeq
where, in accordance with \eqref{notation}, the subscripts indicate the copies in the tensor product algebra $(\ndo\CC^N)^{\ot m}$ with $m=3$ in \eqref{YBE} and $m=2$ in \eqref{unitrig}.

Recall
the  formal Taylor Theorem, 
\beq\label{taylor}
b(z+z_0)=e^{z_0\frac{\partial}{\partial z} b(z)}=
\sum_{k=0}^\infty \frac{z_0^k}{k!} \frac{\partial^k}{\partial z^k} b(z)\quad\text{for}\quad b(z)\in V[[z^{\pm 1}]],
\eeq
where $V$ is a  vector space.
Due to   \eqref{taylor}, we can regard the $R$-matrix $\R(x)$ as an element of $(\ndo\CC^N )^{\ot 2} [[x,h]]$
via the expansion
\beq\label{exp}
\frac{1}{1-e^{ah}x}=
\frac{1}{1-(x+(e^{ah}-1)x)}=
\sum_{k=0}^{\infty} \frac{(e^{ah}-1)^k x^k}{k!}\frac{\partial^k}{\partial x^k}\left(\frac{1}{1-x}\right), \quad a\in\CC,
\eeq
where
\beq\label{expansions}
e^{ah}=\sum_{k\geqslant 0} (ah)^k /k!\in\CC[[h]]\fand
(1-x)^{-1}=\sum_{k\geqslant 0}x^k\in\CC[[x]].
\eeq

Due to \cite{FR}, there exists a unique series 
$f_q (x)$ in $\CC(q)[[x]]$
such that
\beq\label{fqh}
f_q (xq^{2N}) = f_q (x)\frac{\left(1-xq^2\right)\left(1-xq^{2N-2}\right)}{\left(1-x\right)\left(1-xq^{2N}\right)}.
\eeq
As demonstrated in \cite{KM}, the series $f_q (x)$ can be expressed as
\beq\label{f2}
f_q (x)=1+\sum_{k=1}^{\infty} f_{q,k} \left(\frac{x}{1-x}\right)^k,
\eeq
where all  $f_{q,k} (q-1)^{-k} \in\CC(q)$ are regular at $q=1$. Hence,  applying the substitution $q=e^{h/2}$ to \eqref{f2}  and 
using  the expansions in \eqref{expansions} 
  we obtain 
\beq\label{f3}
f(x)\coloneqq 1+\sum_{k=1}^{\infty} f_{k} \left(\frac{x}{1-x}\right)^k\in\CC[[x,h]],\quad\text{where}\quad f_k \coloneqq (f_{q,k})\left|_{q=e^{h/2}}\right.\in h^k\CC[[h]] .
\eeq
By \cite[Equation (2.11)]{KM} series \eqref{f3} satisfies 
\beq\label{fqhqf}
f (x)f  (xe^h)\ldots f  (xe^{(N-1)h})=\frac{1-x}{1-xe^{(N-1)h}}.
\eeq

The normalized $R$-matrix
\beq\label{R}
R(x)=f(x)\R(x)\,\in\, \ndo\CC^N \ot \ndo\CC^N[[x,h]]
\eeq
possesses the {\em crossing symmetry properties}
\beq\label{csym}
R(xe^{Nh})^{t_1}  D_1 ( R(x)^{-1})^{t_1}=D_1\fand (R(x)^{-1})^{t_2} D_2 R(xe^{Nh})^{t_2} = D_2,
\eeq
where $D $ denotes the diagonal matrix
\beq\label{matrix497}
D=\diag\left(e^{ (N-1)h/2 },e^{ (N-3)h/2},\ldots ,e^{- (N-1)h/2} \right)
\eeq
 and $t_i$ denotes the transposition applied on the tensor factor $i=1,2$;
see \cite{FR}.

Express the $R$-matrix $R(x)$ defined by \eqref{R} as
\beq\label{R2}
R(x)=g(x)R^+(x),\quad \text{where}\quad  g(x)=\frac{f(x) }{1-e^{-h}x},\quad R^+(x)=\left(1-e^{-h}x\right)\R(x).
\eeq
Clearly, $R^+ (x)$ is a polynomial with respect to  the variable $x$, i.e.   $R^+ (x)$ belongs to $(\ndo\CC^N)^{\ot 2}   [[h]][x]$.
On the other hand, 
as $(e^{ah}-1) x\in xh\CC[[h]]$,    we conclude by     \eqref{exp} and \eqref{f3} that  $g(x)$ admits the presentation
\beq\label{g}
g(x)=\sum_{k=0}^{\infty} g_k \frac{x^k}{\left(1-x\right)^{k+1}},\quad\text{where}\quad 
g_k\in h^k\CC[[h]]\fand g_0 =1.
\eeq

Denote by  $\CC_* (z_1,\ldots ,z_n)$
 the localization  of the ring of Taylor series
$\CC[[z_1,\ldots ,z_n]]$ at $\CC[z_1,\ldots ,z_n]^{\times}$.
Consider the unique embedding $ \CC_* (z_1,\ldots ,z_n)\to \CC((z_1))\ldots ((z_n))$. Extending  the embedding   to the $h$-adic completion of $\CC_* (z_1,\ldots ,z_n)$ we obtain the map
\beq\label{iotas}
\iotan\colon \CC_* (z_1,\ldots ,z_n)[[h]]\to \CC((z_1))\ldots ((z_n))[[h]].
\eeq
As in \cite{KM}, we now  apply the substitution $x=e^u$ to the normalized $R$-matrix $R(x)$ given by \eqref{R}. 
First, replacing the variable $x$ by $e^u$ in \eqref{g} we obtain
$$
 g(e^u)=\sum_{k=0}^{\infty} g_k \frac{e^{ku}}{\left(1-e^u\right)^{k+1}}=
\sum_{k=0}^{\infty} g_k \frac{e^{ku}\left(\frac{u}{1-e^{u}}\right)^{k+1}}{u^{k+1}}
\in\CC_*(u)[[h]] 
$$
since all numerators $e^{ku}u^{k+1}(1-e^{u})^{-k-1}$ belong to $\CC[[u]]$ and $g_k\in h^k\CC[[h]]$.
By applying the embedding $\iotau$   we get $\iotau  g(e^u) \in \CC((u))[[h]]$.
Next, as   $R^+ (x)$ is a polynomial with respect to the variable $x$, by applying the substitution $x=e^u$ we obtain $R^+(e^u)$, which belongs to $(\ndo\CC^N)^{\ot 2}\, [[h,u]]$. Finally, there exists a unique $\psi\in 1+h\CC[[h]]$ such that  the $R$-matrix
\beq\label{rplusg}
R(e^u)\coloneqq \psi \iotau g(e^u)   R^+ (e^u)  \,\in\, \ndo\CC^N\ot\ndo\CC^N((u))[[h]]
\eeq
possesses the {\em unitarity property}
\beq\label{uni}
R_{12}(e^u) R_{21}(e^{-u}) =1
\eeq
and 
the {\em crossing symmetry properties}
\beq\label{csym2}
R(e^{u+Nh})^{t_1}  D_1 ( R(e^u)^{-1})^{t_1}=D_1\fand (R(e^u)^{-1})^{t_2} D_2 R( e^{u+Nh})^{t_2} = D_2;
\eeq
see \cite[Prop. 1.2]{EK4} and \cite[Prop. 2.1]{KM}.
Of course,  $R$-matrix \eqref{rplusg} also  
satisfies the {\em Yang--Baxter equation}
\beq\label{yberat}
R_{12}(e^u) R_{13}(e^{u+v}) R_{23}(e^v)=R_{23}(e^v)R_{13}(e^{u+v}) R_{12}(e^u).
\eeq

In what follows, whenever it is clear from the context, we      omit   the  embedding symbol $\iotaop$ and write, e.g.,   $f(e^u)$ instead of $\iotau f(e^u)$. Furthermore, in the multiple variable case, we   employ the usual expansion convention where the choice of the embedding is determined by the order of the variables. For example, if $\sigma$ is a permutation in the symmetric group $\mathfrak{S}_n$, then $f(e^{u_{\sigma_1} +\ldots +u_{\sigma_n}})$   denotes $\iotasigma  f(e^{u_{\sigma_1} +\ldots +u_{\sigma_n}})\in\CC((u_{\sigma_1}))\ldots ((u_{\sigma_n} ))[[h]]$. In particular,  by $R_{13}(e^{u+v})$ in \eqref{yberat} is  denoted  $\iotauv g(e^{u+v})   R_{13}^+ (e^{u+v}) $.

%%%%%%%%%%%%%%%%%%%%%%%%%%%%%%%%%%%
\subsection{Quantum affine algebra}\label{subsec021}
%%%%%%%%%%%%%%%%%%%%%%%%%%%%%%%%%%%

Ding's quantum current realization of the quantum affine algebra of type $A$ was given in \cite[Prop. 3.1]{D}. We slightly modify the original  definition \cite[Def. 3.1]{D} in order to make the setting compatible with the quantum vertex algebra theory; see Remark \ref{comparison} for more details.
Our exposition starts in  parallel with \cite[Subsection 2.1]{c11}, where a certain quantum current algebra associated with the suitably normalized Yang $R$-matrix was introduced. We omit some simple proofs as they present a straightforward generalization of the arguments from the aforementioned paper to the trigonometric case.

For any integer $N\geqslant 2$ denote by $\F(N)$  the associative algebra over the ring $\CC[[h]]$ generated by the elements $1$, $C$ and $\lk_{ij}^{(r)}$, where $i,j=1,\ldots, N$ and $r\in\ZZ$, subject to the   defining relations
$$C\cdot a=a\cdot C\fand 1\cdot a=a\cdot 1=a\qquad\text{for all }a\in \F(N),$$
i.e. $1$ is the unit and $C$ is a central element in  $\F(N)$.
Introduce the  Laurent series
\beq\label{lambda}
\lk_{ij}(x)=\delta_{ij}-h\sum_{r\in\ZZ}\lk_{ij}^{(r)}x^{-r-1}\,\in\, \F(N)[[x^{\pm 1}]],\quad\text{where }  i,j=1,\ldots ,N,
\eeq
and   arrange them   into the matrix $\Lc(x)\in\ndo\CC^N\ot\F(N)[[x^{\pm 1}]]$, 
\beq\label{LAMBDA}
\Lc(x)=\sum_{i,j=1}^N e_{ij}\ot \lk_{ij}(x).
\eeq

We now introduce certain completion of the algebra $\F(N)$  which is suitable for expressing the defining relations  for the quantum affine algebra.
For an  integer $p\geqslant 1$ let $\I_p(N)$ be the left ideal in  $\F(N)$ generated by  all  $\lk_{ij}^{(r)}$, where  $i,j=1,\ldots, N$ and $r\geqslant p-1$. Define the completion of  $\F(N)$ as the inverse limit
$$\Ft (N) = \lim_{\longleftarrow} \ts\F(N)\ts /\ts \I_p(N).$$
The algebra  $\Ft (N)$ is naturally equipped with the $h$-adic topology and its $h$-adic completion is equal to $\Ft (N)[[h]]$. 
For any integer $p\geqslant 1$ let $\I_p^h (N)$ be the $h$-adically completed left ideal in $\Ft(N)[[h]]$  generated by $\I_p(N)$ and the element $h^p \cdot 1$.

We   generalize the tensor notation from \eqref{notation} to the matrix $\Lc(x)$ so that the subscript indicates the copy in the corresponding tensor product algebra,
\beq\label{notation2}
\Lc_r (x)=\sum_{i,j=1}^N 1^{\ot (r-1)}\ot e_{ij}\ot 1^{\ot (m-r)}  \ot \lk_{ij}(x)
\,\in\, (\ndo\CC^N)^{\ot m} \ot \F (N)[[x^{\pm 1}]].
\eeq
Employing such notation for $m=2$ and $r=1,2$ we introduce the  expressions
\begin{align*}
&\Lc_{[2]}^{(1)}(x,y)=\Lc_1(x) R_{21}(ye^{-hC}/x) \Lc_2(y)  R_{21}(y/x)^{-1},\\
&\Lc_{[2]}^{(2)}(x,y)= R_{12}(x/y)^{-1}\Lc_2(y) R_{12}(xe^{-hC}/y)\Lc_1(x).
\end{align*}
In accordance with   the discussion in Subsection \ref{sec01}, the $R$-matrices $R_{21}(ye^{ahC}/x)^{\pm 1}$ and $R_{12}(xe^{ahC}/y)^{\pm 1}$ with $a\in \CC$ are regarded as Taylor series with respect to  $y/x$ and $x/y$ respectively.
By arguing as in \cite[Lemma 2.1]{c11}, one can prove

\begin{lem}\label{LR}
The expressions $\Lc_{[2]}^{(1)}(x,y)$ and $\Lc_{[2]}^{(2)}(x,y)$ 
are well-defined elements of
$$ 
\ndo\CC^N  \ot\ndo\CC^N  \ot \Ft(N)[[ x^{\pm 1},y^{\pm 1},h]].
$$
Moreover, for any integer $p\geqslant 1$ both $\Lc_{[2]}^{(1)}(x,y) $ and $\Lc_{[2]}^{(2)}(y,x)$\footnote{Notice the swapped variables in this term.} modulo  $\I_p^h(N) $  belong to
$$
 \ndo\CC^N    \ot \ndo\CC^N    \ot \F(N) [[x^{\pm 1} ]]((y)).
$$
\end{lem}

By Lemma \ref{LR},  there exist elements 
$\lk_{i\ts j\ts k\ts l}^{(r,s;t)}$  in  $\Ft(N)[[h]]$ 
  such that
\begin{align*}
&\Lc_{[2]}^{(t)}(x,y)=\sum_{i,j,k,l=1}^N\sum_{r,s\in\mathbb{Z}} e_{ij}\ot e_{kl}\ot \lk_{i\ts j\ts k\ts l}^{(r,s;t)}\ts x^{-r-1}y^{-s-1} \quad\text{for}\quad t=1,2.
\end{align*}
Let $\J(N)$ be the  ideal in  the algebra $\Ft(N)[[h]]$ generated by all  elements
\beq\label{defrel} 
\lk_{i\ts j\ts k\ts l}^{(r,s;1)}- \lk_{i\ts j\ts k\ts l}^{(r,s;2)},\qquad \text{where} \qquad r,s\in\ZZ\fand i,j,k,l=1,\ldots ,N. 
\eeq
Introduce the completion of $\J(N)$ as the inverse limit
$$
\Jt (N)=\lim_{\longleftarrow} \ts \J(N)\ts /\ts \J(N)\cap \I_p (N).
$$
Note that the $h$-adic completion $[\Jt(N)][[h]]$ of
$$
[\Jt(N)]=\textstyle \left\{a\in\Ft (N)[[h]]\,:\, h^n a\in \Jt(N)\text{ for some integer }n\geqslant 0  \right\}
$$
  is also an ideal in  $\Ft(N)[[h]]$.
Following \cite[Def. 3.1]{D}, we define the  {\em (completed) quantum affine algebra}  $\Uh$ as the quotient of the algebra $\Ft(N)[[h]]$ by  the ideal  $[\Jt (N)][[h]]$,
\beq\label{quotient}
\Uh\ts =\ts \Ft(N)[[h]]\ts /\ts [\Jt (N)][[h]].
\eeq

Denote the images of the elements $1$, $C$ and $\lk_{ij}^{(r)}$ 
in quotient \eqref{quotient} again by $1$, $C$ and $\lk_{ij}^{(r)}$. 
Also,   denote by $\lk_{ij}(x)$ and $\Lc(x)$  the corresponding series 
in $\Uh[[x^{\pm 1}]]$ and   $\ndo\CC^N \ot \Uh [[x^{\pm 1}]]$ respectively. 
Defining relations \eqref{defrel} for the    algebra $\Uh$
 can be expressed  by  the
{\em quantum current 
commutation relation} 
\begin{align}\label{qc}
\Lc_1(x) R_{21}(ye^{-hC}/x) \Lc_2(y)  R_{21}(y/x)^{-1}=R_{12}(x/y)^{-1}\Lc_2(y) R_{12}(xe^{-hC}/y)\Lc_1(x),
\end{align}
as given by  
 Reshetikhin and Semenov-Tian-Shansky \cite{RS}. 
As the images of the elements $\lk_{i\ts j\ts k\ts l}^{(r,s;1)}$ and $\lk_{i\ts j\ts k\ts l}^{(r,s;2)}$ in quotient \eqref{quotient} coincide,   we   denote them by $\lk_{i\ts j\ts k\ts l}^{(r,s)}$. Also, we  write 
\beq\label{eldva}
\Lc_{[2]}(x,y)=\sum_{i,j,k,l=1}^N\sum_{r,s\in\mathbb{Z}} e_{ij}\ot e_{kl}\ot  \lk_{i\ts j\ts k\ts l}^{(r,s)}\ts x^{-r-1}y^{-s-1}
%\,\in\,  (\ndo\CC^N)^{\ot 2} \ot \Uh [[x^{\pm 1},y^{\pm 1}]]
.
\eeq
 and $\Lc_{[1]}(x)=\Lc (x)$.
Observe that  the both sides of    relation \eqref{qc} coincide with  $\Lc_{[2]}(x,y)$.
Motivated by \cite{RS}, we refer to the series $\Lc(x)$ as {\em quantum currents}.
Our next goal is to derive a certain generalized version of \eqref{qc}
consisting of $n+m$ quantum currents.

For  integers $n,m\geqslant 1$ introduce the functions depending on the variable $z$ and the
families of variables
$x=(x_1,\dots,x_n)$ and $y=(y_1,\dots,y_m)$ with values in
the space
$(\ndo\mathbb{C}^{N})^{\ot  n} \otimes
(\ndo\mathbb{C}^{N})^{\ot  m}
$
by
\begin{align}
R_{nm}^{12}(zxe^{ah}/y)= \prod_{i=1,\dots,n}^{\longrightarrow} 
\prod_{j=n+1,\ldots,n+m}^{\longleftarrow} R_{ij}(zx_i e^{ah}/y_{j-n}),\label{rnm12}\\
R_{nm}^{21}(ye^{ah}/zx)= \prod_{i=1,\dots,n}^{\longleftarrow} 
\prod_{j=n+1,\ldots,n+m}^{\longrightarrow} R_{ji}( y_{j-n}e^{ah}/zx_i),\label{rnm123}
\end{align}
where $a\in\CC$ and the arrows indicate the order of the factors. 
For example, we have
$$R_{22}^{12}(zx/y)=R_{14}R_{13}R_{24}R_{23} \fand R_{22}^{21}(y/xz)=R'_{32}R'_{42}R'_{31}R'_{41},$$
where $R_{ij}= R_{ij}(zx_i /y_{j-n})$ and $R'_{ji}=R_{ji}( y_{j-n}/zx_i)$.
The corresponding functions  associated with the $R$-matrix $R^+(x)$ given by \eqref{R2},  $R_{nm}^{+12}(zxe^{ah}/y)$ and $R_{nm}^{+21}(ye^{ah}/zx)$, can be defined analogously. Note that the evaluations of \eqref{rnm12} and \eqref{rnm123} at $z=1$ are well-defined. We denote them by $R_{nm}^{12}(xe^{ah}/y)$ and $R_{nm}^{21}(ye^{ah}/x)$ respectively.
Next, for any  integer $n\geqslant 1$ and the family of variables $x=(x_1,\ldots ,x_n)$   define the functions with values in $(\ndo\CC^N)^{\ot n}$   by
\begin{align}
&R_{[n,a]}(x)=\prod_{i=1,\ldots ,n-1}^{\longrightarrow}\prod_{j=i+1,\ldots ,n}^{\longrightarrow}R_{ji}(x_j e^{-ah} /x_i)^{-1},\label{r1}\\
&\cev{R}_{[n,a]}(x)=\prod_{i=1,\ldots ,n-1}^{\longleftarrow}\prod_{j=i+1,\ldots ,n}^{\longleftarrow}R_{ji}(x_j e^{-ah} /x_i)^{-1},\label{r2}
\end{align}
where $a\in\CC $  and the arrows again indicate the order of the factors.
For example, we have
$$
R_{[n,a]}(x)=R_{21}R_{31}R_{41}R_{32}R_{42}R_{43}
\fand
\cev{R}_{[n,a]}(x)=R_{43}R_{42}R_{32}R_{41}R_{31}R_{21},
$$
where
$R_{ji}=R_{ji}(x_j e^{-ah} /x_i)^{-1}$.
If $a=0$,  we omit the second subscript and  write 
$$
R_{[n]}(x)=R_{[n,0]}(x)\fand\cev{R}_{[n]}(x)=\cev{R}_{[n,0]}(x).
$$ 
Finally, for any integer $n\geqslant 2$ we  generalize $\Lc_{[2]}(x,y)$, as given by \eqref{eldva}, by setting
\beq\label{Ln}
\Lc_{[n]}(x)=\hspace{-4pt}\prod_{i=1,\ldots ,n}^{\longrightarrow}\hspace{-4pt} \left(\Lc_{i}(x_i)R_{i+1\ts i }(x_{i+1}e^{-hC}/x_i) \ldots R_{n\ts i}(x_{n}e^{-hC}/x_i)  \right)\,\cdot\, \cev{R}_{[n]}(x).
\eeq

Denote by  $\I_p^h (\widehat{\gl}_N), \I_p (\widehat{\gl}_N)$   the images of   the left ideals $\I_p^h (N),\I_p  (N) \subset \Ft(N)[[h]]$ in the   algebra $\Uh$ with respect to the canonical map $\Ft(N)[[h]]\to \Uh$. In the next proposition, we use the superscripts $1,2,3$   to indicate the following tensor factors:
$$
\smalloverbrace{(\ndo\CC^N)^{\ot n}}^{1} \ot \smalloverbrace{(\ndo\CC^N)^{\ot m}}^{2}\ot \smalloverbrace{\Uh}^{3}.
$$
The   proposition can be proved by using Lemma \ref{LR}, Yang--Baxter equation \eqref{YBE}, quantum current  commutation relation \eqref{qc} and arguing as in \cite[Prop. 2.4 and 2.5]{c11}. 

\begin{pro}\label{qcgenpro}
For any integers $n,m\geqslant 1$  and the families of variables $x=(x_1,\ldots ,x_n)$ and $y=(y_1,\ldots ,y_m)$ we have:
\begin{enumerate} 
\item The expression $\Lc_{[n]}(x)$ is a well-defined element of
$$(\ndo\CC^N)^{\ot n} \ot \Uh[[x_1^{\pm 1},\ldots ,x_n^{\pm 1}]].$$
\item For any $p\geqslant 1$ the element $\Lc_{[n]}(x)$ modulo $\I_p^h (\widehat{\gl}_N)$ belongs to
$$(\ndo\CC^N)^{\ot n} \ot \Uh ((x_1,\ldots ,x_n)).$$
\item The following quantum current commutation relation holds:
\begin{align}
&\Lc_{[n]}^{13}(x) R_{nm}^{21} (ye^{-hC}/x ) \Lc_{[m]}^{23}(y)R_{nm}^{21} ( y/x)^{-1}\non\\
 &\qquad\qquad=R_{nm}^{12}(x  / y)^{-1}\Lc_{[m]}^{23}(y)R_{nm}^{12}(xe^{-hC} / y) \Lc_{[n]}^{13}(x).\label{qcgen}
\end{align}
Moreover, both sides of \eqref{qcgen} coincide with $\Lc_{[n+m]}(x,y)$.
\end{enumerate}
\end{pro}
Generalizing \eqref{eldva}  
we denote the coefficients of the matrix entries in \eqref{Ln} as follows:
$$
\Lc_{[n]}(x)=\sum_{i_1,j_1,\ldots,i_n,j_n=1}^N\sum_{r_1,\ldots ,r_n\in\mathbb{Z}} e_{i_1 j_1}\ot\ldots \ot e_{i_n j_n}\ot  \lk_{i_1\ts j_1\ldots i_n\ts j_n}^{(r_1,...,r_n)}\ts x^{-r_1-1}_1\ldots x_n^{-r_n-1}.
$$

Our next goal is to introduce the quantum affine algebra associated with the affine Lie algebra $\wht{\mathfrak{sl}}_N$. 
Let $P^h $ be the $h$-permutation operator,
$$
P^h = \sum_{i=1}^N e_{ii}\ot e_{ii} + e^{h/2}\sum_{\substack{i,j=1\\i> j}}^N e_{ij}\ot e_{ji} +e^{-h/2}\sum_{\substack{i,j=1\\i< j}}^N e_{ij}\ot e_{ji}\in\ndo\CC^N \ot \ndo\CC^N [[h]].
$$
Consider the action of the symmetric group $\mathfrak{S}_n$ on the space $(\CC^N)^{\ot n}$ which is given by $\sigma_i\mapsto P_{\sigma_i}^h=P_{i\ts i+1}^h$ for $i=1,\ldots ,n-1$, where
$\sigma_i$ is the transposition $(i,i+1)$. For a reduced decomposition of a permutation $\sigma=\sigma_{i_1}\ldots \sigma_{i_k} \in \mathfrak{S}_n$ set $P_\sigma^h= P^h_{\sigma_{i_1}}\ldots P^h_{\sigma_{i_k}}$. Let $A^{(n)}$ be the image of the normalized anti-symmetrizer with respect to this action, so that
\beq\label{anti}
A^{(n)}=\frac{1}{n!}\sum_{\sigma\in\mathfrak{S}_n} \sgn\sigma \cdot P_\sigma^h .
\eeq

Define the {\em quantum determinant} of the matrix $\Lc(x)$ by
\beq\label{qdet497}
\qdet \Lc (x)=\tr_{1,\ldots ,N} \,A^{(N)}\ts\Lc_{[N]}(x_1,\ldots ,x_N)\big|_{x_1 = x,\ldots, x_N = xe^{-(N-1)h}}\big. \ts D_1\ldots D_N ,
\eeq
where the trace is taken over all $N$ copies of $\ndo\CC^N$ and the matrix $D$ is given by \eqref{matrix497}. The quantum determinant is a formal power series in the variable $x$ with coefficients in the quantum affine algebra, i.e. $\qdet \Lc (x) $ belongs to $  \Uh[[x^{\pm 1}]]$. Indeed, the substitution $x_1 = x,\ldots, x_N = xe^{-(N-1)h}$ in \eqref{qdet497} is well-defined due to 
the second assertion of 
Proposition \ref{qcgenpro}. Furthermore, all coefficients $d_r$ of the quantum determinant
\beq\label{qdetc}
\qdet \Lc (x)= 1 - h\sum_{r\in\ZZ} d_r x^r
\eeq
belong to the center of the quantum affine algebra  at the level $c\in\CC$; see Proposition \ref{qdetpro}.

 Let $\I_\qdet$ be the ideal in the algebra $\Uh$ generated by the elements $d_r$, where $r\in \ZZ$.
Introduce its completion  as the inverse limit
$$
\wtld{\I}_\qdet=\lim_{\longleftarrow} \ts \I_\qdet\ts /\ts \I_\qdet\cap \I_p (\gl_N).
$$
The $h$-adic completion $[\wtld{\I}_\qdet][[h]]$ of
$$
[\wtld{\I}_\qdet]=\textstyle \left\{a\in\Uh\,:\, h^n a\in \wtld{\I}_\qdet\text{ for some integer }n\geqslant 0  \right\}
$$
  is also an ideal in   $\Uh$.
Define the {\em (completed) quantum affine algebra}  $\Uhsl$ as the quotient of the algebra $\Uh$ by  the relation $\qdet \Lc (x)= 1$, i.e. 
$$
\Uhsl\ts =\ts \Uh\ts/\ts  [\wtld{\I}_\qdet][[h]].
$$

\begin{rem}\label{comparison}
In Ding's definition \cite[Def. 3.1]{D}, the quantum affine algebra is introduced as an associative algebra over the field $\CC(q)$. However, as our goal is to study quantum vertex algebras associated to quantum affine algebras, we used the identification $q=e^{h/2}$ and introduced the quantum affine algebra as a  suitably completed associative algebra over the commutative ring $\CC[[h]]$. Thus we established the setting compatible with  Etingof--Kazhdan's notion of quantum vertex algebra \cite[Sect. 1.4]{EK}, which, in particular, is required to be a topologically free $\CC[[h]]$-module; see also Li's notion of $h$-adic quantum vertex algebra \cite[Def. 2.20]{Li}.
Furthermore, in contrast with Ding's realization, we use  normalized  $R$-matrix \eqref{R} instead of \eqref{Rbar}. Such choice of the $R$-matrix enables the constructions of certain large families of central elements of the quantum affine algebra at the critical level and of the topological generators of the quantum Feigin--Frenkel center, as demonstrated in \cite{FJMR} and \cite{KM} respectively; see also Section \ref{newsec02}. 
\end{rem}

%%%%%%%%%%%%%%%%%%%%%%%%%%%%%%%%%%%
\subsection{Restricted modules}\label{subsec022}
%%%%%%%%%%%%%%%%%%%%%%%%%%%%%%%%%%%

Recall that a $\CC[[h]]$-module $W$ is said to be {\em torsion-free} if $h w= 0$ for all nonzero $w\in W$ and that $W$ is said to be {\em separable} if $\cap_{n\geqslant 1} h^n W=0$. Moreover, $W$ is said to be {\em topologically free} if it is separable, torsion-free and complete with respect to $h$-adic topology; see  \cite[Chapter XVI]{Kas}. 

Let $\g_N=\gl_N ,\sll_N$. By arguing as in \cite[Prop. 2.2]{c11} one can show that the algebra $\Uhg$ is topologically free.
Define 
a {\em restricted} $\Uhg$-module $W$ as a topologically free $\CC[[h]]$-module such that
\beq\label{restricted}
\Lc(x)w\in \ndo\CC^N\ot W((x))[[h]]\quad \text{for all } w\in W.
\eeq
\begin{pro}\label{restricted496}
Let $W$ be a restricted $\Uhg$-module. Then for any $n\geqslant 1$ and the variables $x=(x_1,\ldots ,x_n) $ we have
\beq\label{restrictedf}
\Lc_{[n]}(x)w\in (\ndo\CC^N)^{\ot n}\ot W((x_1,\ldots ,x_n))[[h]]\quad \text{for all } w\in W.
\eeq
\end{pro}

\begin{prf}
Apply  quantum current commutation relation \eqref{qc} on an arbitrary element  of some restricted module. For every integer $k\geqslant 0$ the left hand side contains finitely many negative powers of the variable $y$ modulo $h^k$ while the right hand side contains finitely many negative powers of the variable $x$ modulo $h^k$. Hence the statement of the proposition holds for $n=2$. The case $n>2$ is proved by induction on $n$ which relies on   \eqref{qcgen}.
\end{prf}

\begin{rem}\label{linremark}
Note that  \eqref{restrictedf} implies  $\Lc_{[n]}(x)\in\ndo\CC^N\ot\om(W,W((x_1,\ldots,x_n))[[h]])$ for all $n\geqslant 1$. Hence we can apply  the substitutions $x_1 = ze^{u_1}, \ldots, x_n=ze^{u_n}$, thus getting
\begin{align}
\Lc_{[n]}(x_1,\ldots ,x_n)\big|_{x_1 = ze^{u_1},\ldots , x_n=ze^{u_n}}\big. \in
\ndo\CC^N\ot\om(W,W((z))[[h,u_1,\ldots ,u_n]]).\label{lin}
\end{align}
We will often denote the expression in \eqref{lin} more briefly by $\Lc_{[n]}(x)\left|_{x_i = ze^{u_i}}\right.$.
\end{rem}
 
As usual, an
$\Uhg$-module $W$ is said to be of {\em level $c$} if the central element $C\in \Uhg$ acts on $W$ as a  scalar multiplication by some $c\in\CC$.
Denote by $\Uhg_c$ the {\em quantum affine algebra at the level $c$}, i.e. the quotient of $\Uhg$ by the ideal generated by the element $C-c$.
 Let $\K_c $ be the left ideal in the  algebra $\Uhg_c$ generated by all elements 
$$ \lk_{i_1\ts j_1\ts\ldots\ts i_n\ts j_n }^{(r_1,\ldots ,r_n)}\quad\text{such that}\quad r_k\geqslant 0\text{  for at least one integer  }k=1,\ldots,n,$$
 where $n\geqslant 1$, $i_1,\ldots, i_n,j_1,\ldots ,j_n=1,\ldots , N$ and $r_1,\ldots,r_n\in \ZZ$.
Introduce the completion of $\K_c  $ as the inverse limit
$$\wtld{\K}_c  =\lim_{\longleftarrow} \ts \K_c  \ts /\ts \K_c  \cap \I_p (\wht{\g}_N).$$
Then the $h$-adic completion $[\wtld{\K}_c][[h]]$ of
$$[\wtld{\K}_c]=\left\{a\in \Uhg_c\,:\, h^n a\in \wtld{\K}_c\text{ for some }n\geqslant 0\right\}$$
is also a left ideal in $\Uhg_c$.
Define the {\em vacuum module $\Vccglg$ at the level $c$} over the quantum affine algebra $\Uhg$  as the   quotient of   $\Uhg_c $ by its left ideal  $[\wtld{\K}_c][[h]]$,   
\beq\label{quotient2}
 \Vccglg\,=\,\Uhg_c\,/\, [\wtld{\K}_c][[h]].
\eeq

Observe that the canonical map $\Uhg_c\to \Vccglg$ maps  the left ideal $\I_p^h(\wht{\g}_N)$   to $h^p  \Vccglg$. Denote by $\vac$ the image of the unit $1 \in\Uhg$   with respect to this map.

\begin{pro}\label{free2}
The vacuum module $\Vccglg$ is a topologically free $\CC[[h]]$-module. Moreover, it is a restricted $\Uhg$-module.
\end{pro}
 
\begin{prf}
The first assertion is verified by arguing as in \cite[Prop. 2.2]{c11}. As for the second assertion, we first observe that all elements
\beq\label{elsofform}
\lk_{i_1\ts j_1\ts\ldots\ts i_n\ts j_n }^{(r_1,\ldots ,r_n)}\vac\qquad \text{such that}\qquad n\geqslant 0\fand r_k< 0 \text{ for all }k=1,\ldots,n
\eeq 
span an $h$-adically dense $\CC[[h]]$-submodule of $\Vccglg$.
Indeed, this follows from the fact that each monomial
$\lk_{k_1\ts l_1}^{(s_1)}\ldots \lk_{k_m\ts l_m}^{(s_m)}\vac\in \Vccglg$ can be expressed
using   elements \eqref{elsofform}. 
This is  done by employing  crossing symmetry properties \eqref{csym} and  invertibility of the trigonometric $R$-matrix to move all  $R$-matrices which appear on the right hand side of
\beq\label{qcvac}
\Lc_{[a+b]}(x,y)\vac=\Lc_{[a]}^{13}(x) R_{ab}^{21} (ye^{-hC}/x ) \Lc_{[b]}^{23}(y)R_{ab}^{21} ( y/x)^{-1}\vac,
\eeq
where $a+b=m$,  $x=(x_1,\ldots ,x_a)$ and $y=(y_1,\ldots ,y_b)$,
to the left hand side  (for more details see Remark \ref{csrem}), and then taking the  coefficient of 
$
x_1^{-s_1-1}\hspace{-2pt}\ldots x_{a}^{-s_a -1} y_1^{-s_{a+1}-1}\hspace{-2pt}\ldots y_b^{-s_m-1}
$
at the matrix entry $e_{k_1\ts l_1}\ot \ldots \ot e_{k_m\ts l_m}$. 
Note that    \eqref{qcvac} follows from Proposition \ref{qcgenpro}.

Therefore, it is sufficient to check that
$\Lc(z)w$ belongs to $\ndo\CC^N\ot \Vccglg((z))[[h]]$ for all  $w\in \Vccglg$ of the form as in \eqref{elsofform}. However, as   \eqref{qcvac} contains only nonnegative powers of the variables $x_1,\ldots ,x_a,y_1,\ldots ,y_b$, this follows    by setting $a=1$ and $b=n$ in \eqref{qcvac}, then moving   $R_{1n}^{21} (ye^{-hC} /x )$ and $R_{1n}^{21} ( y/x)^{-1}$ to the left hand side and, finally, by  taking the  coefficient of 
$ y_1^{-r_{1}-1}\ldots y_n^{-r_n-1}
$
at the matrix entries $e_{ij}\ot e_{i_1\ts j_1}\ot \ldots \ot e_{i_n\ts j_n}$ for $i,j=1,\ldots ,N$.
\end{prf}
Observe that Proposition \ref{free2} and the \hyperref[mainthm1]{Main Theorem} imply Corollary \ref{maincor}.

%%%%%%%%%%%%%%%%%%%%%%%%%%%%%%%%%%%
%%%%%%%%%%%%%%%%%%%%%%%%%%%%%%%%%%%
\section{\texorpdfstring{$\phi$}{phi}-Coordinated modules  for the quantum affine vertex algebra}\label{sec0102}
%%%%%%%%%%%%%%%%%%%%%%%%%%%%%%%%%%%
%%%%%%%%%%%%%%%%%%%%%%%%%%%%%%%%%%%
In this section, we recall Etingof--Kazhdan's construction of the quantum affine vertex algebra associated with trigonometric $R$-matrix in type $A$. Next, we suitably modify Li's definition of $\phi$-coordinated module, thus establishing the   setting for the    \hyperref[mainthm1]{Main Theorem}.

%%%%%%%%%%%%%%%%%%%%%%%%%%%%%%%%%%%
\subsection{Quantum affine vertex algebra}\label{subsec012}
%%%%%%%%%%%%%%%%%%%%%%%%%%%%%%%%%%%

We follow \cite{EK3,EK4} to introduce the $R$-matrix algebras $\Uhplusg$; see also \cite{FRT,RS}.
Let  $\Uhplus$ be the associative algebra over the ring $\CC[[h]]$ generated by elements $t_{ij}^{(-r)}$, where $i,j=1,\ldots ,N$ and $r=1,2,\ldots$, subject to the  defining relations
\beq\label{rtt}
R(e^{u-v})\ts T_{1}^+(u)\ts T_2^+ (v)=  T_2^+ (v)\ts T_{1}^+(u)\ts R(e^{u-v}),
\eeq
where $T^+(u) \in \ndo\CC^N \ot \Uhplus[[u]]$ is given by
$$
T^+(u) =\sum_{i,j=1}^N e_{ij}\ot t_{ij}^+ (u)\qquad\text{for}\qquad t_{ij}^+(u)=\delta_{ij}-h\sum_{r=1}^{\infty}t_{ij}^{(-r)}u^{r-1}\in\Uhplus[[u]].
$$
As in \eqref{notation2}, 
we use subscripts  in \eqref{rtt} to indicate copies in the tensor product algebra
$\ndo\CC^N \ot \ndo\CC^N \ot\Uhplus$. 
Note that the $R$-matrix
$R(e^{u-v})$ in defining relation \eqref{rtt} can be replaced by $R^+ (e^{u-v})$.

Define the {\em quantum determinant} of the matrix $T^+(u)$ by
\beq\label{qdetvoa}
\qdet T^+ (u)=\tr_{1,\ldots ,N} \,A^{(N)}\ts T_1^+(u)\ldots T_N^+(u-(N-1)h)\ts D_1\ldots D_N ,
\eeq
where the trace is taken over all $N$ copies of $\ndo\CC^N$ and the matrix $D$ is given by \eqref{matrix497}. The quantum determinant   $ \qdet T^+ (u)$ belongs to $\Uhplus [[u]]$. Moreover, its coefficients $\delta_r$, which are given by
\beq\label{detkoef}
\qdet T^+ (u)=1-h\sum_{r\geqslant 0}\delta_r u^r,
\eeq
belong to the center of the algebra $\Uhplus$; see proof of \cite[Prop. 3.10]{KM}.
Define the algebra $\Uhslplus$ as 
the quotient of $\Uhplus$ over the $h$-adically completed ideal generated by the elements $\delta_0,\delta_1,\ldots$ 
Hence we have the following relation in $\Uhslplus$:
\beq\label{qdetvoa2}
\qdet T^+ (u)=1 .
\eeq

Let $\g_N =\gl_N, \sll_N$.
For positive integers $n$ and $m$ we extend the notation in \eqref{rnm12} and \eqref{rnm123} by introducing the functions depending on the variable $z$ and
the families of variables $u = (u_1 ,...,u_n )$ and $v = (v_1 ,...,v_m )$ with values in the space
$(\ndo\CC^N )^{\ot n} \ot (\ndo\CC^N )^{\ot m}$
 by
\begin{align}
&R_{nm}^{12}(e^{z+u-v+ah})= \prod_{i=1,\dots,n}^{\longrightarrow} 
\prod_{j=n+1,\ldots,n+m}^{\longleftarrow}   R_{ij} (e^{z+u_i-v_{j-n}+ah}),\label{rnm12exp}\\
&R_{nm}^{21}(e^{z+u-v+ah})= \prod_{i=1,\dots,n}^{\longleftarrow} 
\prod_{j=n+1,\ldots,n+m}^{\longrightarrow} R_{ji}( e^{z+u_i-v_{j-n}+ah}),\label{rnm123exp} 
\end{align}
where $a\in\CC$. Note that the expansion convention, as introduced at the end of  Subsection \ref{sec01}, is   applied  on every factor on the right hand side, i.e. 
$$R_{ij} (e^{z+u_i-v_{j-n}+ah}) = \psi \iotazuv g (e^{z+u_i-v_{j-n}+ah})\ts R_{ij}^+ (e^{z+u_i-v_{j-n}+ah}).$$
If the variable $z$ is omitted in \eqref{rnm12exp} or \eqref{rnm123exp}, the embeddings $\iotauiv$ are applied on the corresponding normalizing functions $g(e^{u_i-v_{j-n}+ah})$ instead.
The functions $R_{nm}^{+12}(e^{z+u-v+ah})$ and $R_{nm}^{+21}(e^{z+u-v+ah})$ corresponding to the $R$-matrix $R^+(x)$ given by \eqref{R2} can be defined analogously.
Denote by $\vac$ the unit in the algebra $\Uhplusg $. We recall \cite[Lemma 2.1]{EK}:
\begin{lem}\label{lemma21}
For any $c\in\CC$ there exists a unique operator series
$$T^*(u)\in\ndo\CC^N \ot \om( \Uhplusg [[h]],\Uhplusg ((u))[[h]] )$$
such that   for all $n\geqslant 0$ we have
\begin{align}
&T^{*}_{1} (u)\ts T_{2}^{+}(v_1)\ldots T_{n+1}^+(v_n)\vac\non\\
&\qquad\qquad = R_{1n}^{12}(e^{u-v+hc/2})^{-1} T_{2}^{+}(v_1)\ldots T_{n+1}^+(v_n) R_{1n}^{12}(e^{u-v-hc/2})\vac.\label{tstar}
\end{align}
\end{lem}

In order to indicate action \eqref{tstar}, which is uniquely determined by the scalar $c\in\CC$, we denote the topologically free $\CC[[h]]$-module $\Uhplusg [[h]]$ by $\Vc_{\hspace{-1pt}c}(\g_N)$.  
Following \cite{EK}, we introduce the operators on $(\ndo\CC^N)^{\ot n} \ot \Vc_{c}(\g_N)$ by
\begin{align*}
T_{[n]}^{+}(u|z)=T_{1}^{+}(z+u_1)\ldots T_{n}^{+}(z+u_n)\fand T_{[n]}^*(u|z)=T_{1}^*(z+u_1)\ldots T_{n}^*(z+u_n).
\end{align*}
By the expansion convention from Subsection \ref{sec01}, the operator $T_{[n]}^* (u|z)$ contains only nonnegative powers of the variables $u_1,\ldots ,u_n$ as the embeddings $\iotazui$ are applied on its corresponding factors.
If the variable $z$ is omitted, we  write
\beq\label{rnm1234t}
T_{[n]}^{+}(u )=T_{1}^{+}(u_1)\ldots T_{n}^{+}(u_n)\fand T_{[n]}^*(u)=T_{1}^*(u_1)\ldots T_{n}^*(u_n).
\eeq
The next proposition, as given in  \cite[Prop. 2.2]{EK}, is  verified using \eqref{rtt} and \eqref{tstar}.
In  relations \eqref{rtt2}--\eqref{rtt3}, the superscripts $1,2,3$ indicate the  tensor factors as follows:
$$
\smalloverbrace{(\ndo\CC^N)^{\ot n}}^{1} \ot \smalloverbrace{(\ndo\CC^N)^{\ot m}}^{2}\ot \smalloverbrace{\Vc_c(\g_N)}^{3}.
$$
For example, the superscripts $1,3$ in $T_{[n]}^{*13}(u|z_1)$ indicate that  the operator
$T_{[n]}^{*}(u|z_1)$
 is applied on the tensor factors $1,\ldots ,n$ and $n+m+1$. 

\begin{pro} For any integers $n,m\geqslant 1$ and the families of variables $u=(u_1,\ldots ,u_n)$
 and $v=(v_1,\ldots ,v_m)$ the following equalities hold on $\Vc_c(\g_N)$:
\begin{align}
&R_{nm}^{12}(e^{z_1-z_2+u-v})T_{[n]}^{*13}(u|z_1)T_{[m]}^{*23}(v|z_2) 
=\,T_{[m]}^{*23}(v|z_2)T_{[n]}^{*13}(u|z_1)R_{nm}^{12}(e^{z_1-z_2+u-v}),\label{rtt2}\\
&R_{nm}^{12}(e^{z_1-z_2+u-v})T_{[n]}^{+13}(u|z_1)T_{[m]}^{+23}(v|z_2) 
=\,T_{[m]}^{+23}(v|z_2)T_{[n]}^{+13}(u|z_1)R_{nm}^{12}(e^{z_1-z_2+u-v}),\label{rtt1}\\ 
&R_{nm}^{\ts 12}(e^{z_1-z_2+u-v+hc/2})T_{[n]}^{*13}(u|z_1)T_{[m]}^{+23}(v|z_2)
=\,T_{[m]}^{+23}(v|z_2)T_{[n]}^{*13}(u|z_1)R_{nm}^{\ts 12}(e^{z_1-z_2+u-v-hc/2}).
\label{rtt3}
\end{align}
\end{pro}

From now on, the tensor products are understood as $h$-adically completed.
The notion of {\em quantum vertex algebra} was introduced by Etingof and Kazhdan  \cite{EK}. It is defined as a quadruple $(V,Y,\vac,\Sc)$ such that
\begin{enumerate}[wide, labelwidth=!, labelindent=0pt]
\item[1. ]  $V$ is a topologically free $\mathbb{C}[[h]]$-module.

\item[2. ]
$Y=Y(z)$ is the {\em vertex operator map}, i.e. a $\mathbb{C}[[h]]$-module map 
\begin{align*}
Y \colon V\ot V&\to V((z))[[h]]\\
u\ot v&\mapsto Y(z)(u\ot v)=Y(u,z)v=\sum_{r\in\mathbb{Z}} u_r v \ts z^{-r-1}
\end{align*}
which satisfies the {\em weak associativity}:
for any $u,v,w\in V$ and $n\in\mathbb{Z}_{\geqslant 0}$
there exists $p\in\mathbb{Z}_{\geqslant 0}$
such that
\begin{align}
&(z_0 +z_2)^p\ts Y(u,z_0 +z_2)Y(v,z_2)\ts w  - (z_0 +z_2)^p\ts Y\big(Y(u,z_0)v,z_2\big)\ts w
\in h^n V[[z_0^{\pm 1},z_2^{\pm 1}]].\label{associativity}
\end{align}

\item[3. ] $\vac$ is the {\em vacuum vector}, i.e. a distinct element of $V$ satisfying
\beq\label{v1}
Y(\vac ,z)v=v,\quad Y(v,z)\vac\in V[[z]]\fand \lim_{z\to 0} Y(v,z)\ts\vac =v\quad\text{for all }v\in V,
\eeq

\item[4. ] $\Sc=\Sc(z)$ is the {\em braiding map}, i.e. a $\mathbb{C}[[h]]$-module map
$ V\otimes V\to V\otimes V\otimes\mathbb{C}((z))[[h]]$ which satisfies
  the $\mathcal{S}$-{\em locality}:
for any $u,v\in V$ and $n\in\mathbb{Z}_{\geqslant 0}$ there exists
$p\in\mathbb{Z}_{\geqslant 0}$ such that   for all $w\in V$
\begin{align}
(z_1-z_2)^{p}\ts Y(z_1)\big(1\otimes Y(z_2)\big)\big(\mathcal{S}(z_1 -z_2)(u\otimes v)\otimes w\big)&
\nonumber\\
 -(z_1-z_2)^{p}\ts Y(z_2)\big(1\otimes Y(z_1)\big)(v\otimes u\otimes w)
&\in h^n V[[z_1^{\pm 1},z_2^{\pm 1}]].\label{locality}
\end{align}
\end{enumerate}
The given data  should  posses several other properties which  we omit as they  are not used in this paper;  for a complete definition see  \cite[Sect. 1.4]{EK}.
Finally, we recall  Etingof--Kazhdan's construction \cite[Thm. 2.3]{EK} in the trigonometric $R$-matrix case:
\begin{thm}\label{EK:qva}
For any $c\in \CC$
there exists a unique  quantum vertex algebra  structure on $\Vc_c(\g_N)$
  such that the vertex operator map $Y$ is given by
\beq\label{qva1}
Y\big(T_{[n]}^+ (u)\vac,z\big)=T_{[n]}^+ (u|z)\ts T_{[n]}^* (u|z+hc/2)^{-1}, 
\eeq
the vacuum vector is $\vac\in \Vc_c(\g_N)$
  and the braiding  map $\mathcal{S}(z)$ is defined by the relation  
\begin{align}
\mathcal{S}(z)\big(R_{nm}^{  12}(e^{z+u-v})^{-1}  T_{[m]}^{+24}(v) 
R_{nm}^{  12}(e^{z+u-v-h  c})  T_{[n]}^{+13}(u)(\vac\otimes \vac) \big)&\nonumber\\
 =T_{[n]}^{+13}(u)  R_{nm}^{  12}(e^{z+u-v+h  c})^{-1} 
 T_{[m]}^{+24}(v)  R_{nm}^{  12}(e^{z+u-v})(\vac\otimes \vac)\label{qva2}&
\end{align}
for operators on
$
(\ndo\mathbb{C}^{N})^{\otimes n} \otimes
(\ndo\mathbb{C}^{N})^{\otimes m}\otimes \Vc_c(\g_N) \ot \Vc_c(\g_N) $.
\end{thm}

\begin{rem}\label{csrem}
Crossing symmetry properties \eqref{csym2} of $R$-matrix \eqref{rplusg} can be expressed using the ordered product notation as
\beq\label{csym_equiv}
(D_1 R(e^{u+hN}) D_1^{-1})\cdotrl R(e^{u})^{-1}=1
\fand
(D_2 R(e^u)^{-1} D_2^{-1})\cdotlr R(e^{u+hN})=1,
\eeq
where the subscript RL (LR) indicates that the first tensor factor of $D_i R(e^u)^{-1} D_i^{-1}$, $i=1,2$,  is applied from the right (left) while the second tensor factor  is applied from the left (right). 
Such  notation naturally extends to the products of multiple  $R$-matrices such as \eqref{rnm12exp} and \eqref{rnm123exp}. For example, by \eqref{csym_equiv}, we have
\beq\label{crossingexample}
\left(D_{[m]}^2 R_{nm}^{  12}(e^{z+u-v-h (N+ c)})^{-1}(D_{[m]}^2)^{-1}\right)
\cdotlr
R_{nm}^{  12}(e^{z+u-v-h  c})  =1,
\eeq
where $D_{[m]}^2 =1^{\ot n}\ot D^{\ot m}$ and the subscript LR now indicates that the tensor factors $1,\ldots ,n$ ($n+1,\ldots ,n+m$) are applied from the left (right). 
As with \eqref{csym_equiv},
one can write crossing symmetry properties \eqref{csym} of $R$-matrix \eqref{R} using the ordered product notation. As before, the notation naturally extends to the    multiple  $R$-matrix products such as
 \eqref{rnm12}--\eqref{r2}. 
\end{rem}

Combining  
\eqref{qva2} and \eqref{crossingexample} we find the explicit formula for the action of the braiding,
\begin{align}
&\mathcal{S}(z)\big(   T_{[n]}^{+13}(u)T_{[m]}^{+24}(v) (\vac\otimes \vac) \big) =\left(D_{[m]}^2 R_{nm}^{  12}(e^{z+u-v-h (N+ c)})^{-1}(D_{[m]}^2)^{-1}\right)\non
 \\
 & \qquad\cdotlr \big( R_{nm}^{  12}(e^{z+u-v}) T_{[n]}^{+13}(u)  R_{nm}^{  12}(e^{z+u-v+h  c})^{-1} 
 T_{[m]}^{+24}(v)  R_{nm}^{  12}(e^{z+u-v})(\vac\otimes \vac)\big). \label{sop1}
\end{align}

\begin{rem}\label{remark17}
As with \eqref{exp}, by formal Taylor Theorem \eqref{taylor} we have
$$
\frac{1}{1-xe^{u-v+ah}}=\sum_{k=0}^{\infty} \frac{(e^{u-v+ah}-1)^k x^k}{k!}\frac{\partial^k}{\partial x^k}\left(\frac{1}{1-x}\right).
$$
Therefore, due to    \eqref{R2},
we can regard the $R$-matrix $R(xe^{u-v+ah})$ as an element of $(\ndo\CC^N)^{\ot 2}(x)[[u,v,h]]$ for any $a\in\CC$, i.e. as a rational function in the variable $x$. Clearly, applying the embedding $\iotaxuv$  we obtain an element of $(\ndo\CC^N)^{\ot 2}((x))[[u,v,h]]$. 
\end{rem}

We now extend the notation   \eqref{rnm12exp} and \eqref{rnm123exp} by introducing the functions depending on the variable $x$ and
the families of variables $u = (u_1 ,...,u_n )$ and $v = (v_1 ,...,v_m )$ with values in the space
$(\ndo\CC^N )^{\ot n} \ot (\ndo\CC^N )^{\ot m}$
 by
\begin{align}
&R_{nm}^{12}(xe^{u-v+ah})= \prod_{i=1,\dots,n}^{\longrightarrow} 
\prod_{j=n+1,\ldots,n+m}^{\longleftarrow}   R_{ij} (xe^{u_i-v_{j-n}+ah}),\label{rnm12expx}\\
&R_{nm}^{21}(xe^{u-v+ah})= \prod_{i=1,\dots,n}^{\longleftarrow} 
\prod_{j=n+1,\ldots,n+m}^{\longrightarrow} R_{ji}( xe^{u_i-v_{j-n}+ah}),\label{rnm123expx} 
\end{align}
where $a\in\CC$. In accordance with Remark \ref{remark17}, the $R$-matrices in  \eqref{rnm12expx} and \eqref{rnm123expx} are regarded as rational functions in  the variable $x$.  We   use the map given by the following lemma in Definition \ref{phimod} below,  to introduce the notion of $\phi$-coordinated $\Vc_{\hspace{-1pt}c}(\g_N)$-module.
\begin{lem}\label{eshet}
There exists a unique $\CC[[h]]$-module map 
$$\wht{\Sc}(x)\colon \Vc_{\hspace{-1pt}c}(\g_N)\ot \Vc_{\hspace{-1pt}c}(\g_N) \to  \Vc_{\hspace{-1pt}c}(\g_N)\ot \Vc_{\hspace{-1pt}c}(\g_N) (x)[[h]]$$ 
such that 
\begin{align}
&\wht{\Sc}(x)\big(   T_{[n]}^{+13}(u)T_{[m]}^{+24}(v) (\vac\otimes \vac) \big) =\left(D_{[m]}^2 R_{nm}^{  12}(xe^{u-v-h (N+ c)})^{-1}(D_{[m]}^2)^{-1}\right)\non
 \\
  &\qquad\cdotlr \big( R_{nm}^{  12}(xe^{u-v}) T_{[n]}^{+13}(u)  R_{nm}^{  12}(xe^{u-v+h  c})^{-1} 
 T_{[m]}^{+24}(v)  R_{nm}^{  12}(xe^{u-v})(\vac\otimes \vac)\big). \label{sop1hat}
\end{align}
Moreover, the map $\wht{\Sc}(x)$ satisfies
\begin{align}
\wht{\mathcal{S}}(x)\big(R_{nm}^{  12}(xe^{u-v})^{-1}  T_{[m]}^{+24}(v) 
R_{nm}^{  12}(xe^{u-v-h  c})  T_{[n]}^{+13}(u)(\vac\otimes \vac) \big)&\nonumber\\
 =T_{[n]}^{+13}(u)  R_{nm}^{  12}(xe^{u-v+h  c})^{-1} 
 T_{[m]}^{+24}(v)  R_{nm}^{  12}(xe^{u-v})(\vac\otimes \vac)&.\label{qva2hat}
\end{align}
\end{lem}

\begin{prf}
The map $\wht{\mathcal{S}}(x)$ is well-defined by \eqref{sop1hat}, i.e.   it maps the ideal of relations \eqref{rtt}, and \eqref{qdetvoa2} in the $\g_N=\sll_N$ case, to itself.
Indeed, 
this follows by a straightforward calculation which relies on 
the identity 
$$
R_{12}(e^{u-v})D_1 D_2 =D_2 D_1 R_{12}(e^{u-v})
$$
and
the following version of Yang--Baxter equation   \eqref{YBE}:
$$
R_{12}(e^{u-v}) R_{13} (xe^{u +\alpha h}) R_{23} (xe^{ v+\alpha h})
=
R_{23} (xe^{ v+\alpha h})R_{13} (xe^{u +\alpha h})R_{12}(e^{u-v}),\quad \alpha\in\CC.
$$
Moreover, the proof in the $\g_N=\sll_N$ case employs identity \eqref{fqhqf} and some  properties of the anti-symmetrizer $A^{(N)}$, which are given by \eqref{497a}, \eqref{497d} and
$$
A^{(N)}T_1^+(u) T_2^+(u-h)\ldots T_N^+(u-(N-1)h)=T_N^+(u-(N-1)h)\ldots  T_2^+(u-h)T_1^+(u)A^{(N)},
$$
see  \cite[Equality (3.12)]{KM}.
As for relation \eqref{qva2hat}, it follows from \eqref{sop1hat} and the equality
$$
\left(D_{[m]}^2 R_{nm}^{  12}(xe^{u-v-h (N+ c)})^{-1}(D_{[m]}^2)^{-1}\right)\cdotlr
R_{nm}^{  12}(xe^{u-v-h  c}) =1,
$$
which is verified by using crossing symmetry properties \eqref{csym}; recall Remark \ref{csrem}.
\end{prf}

%%%%%%%%%%%%%%%%%%%%%%%%%%%%%%%%%%%
\subsection{\texorpdfstring{$\phi$}{phi}-Coordinated modules}
%%%%%%%%%%%%%%%%%%%%%%%%%%%%%%%%%%%

The notion of $\phi$-coordinated   module, 
%for a nonlocal vertex algebra, 
where $\phi$ is an associate of the one-dimensional additive formal group, was introduced by Li  \cite{Li1}.
As in \cite[Sect. 5]{Li1}, throughout this paper we    consider the associate
\beq\label{phi}
\phi(z_2,z_0) = z_2 e^{z_0}.
\eeq
Before we proceed to the definition of $\phi$-coordinated   module, we  introduce some notation.
Let $V$ be a topologically free $\CC[[h]]$-module and $a_1,\ldots ,a_n,k>0$ arbitrary integers. Suppose that some element $A$ of $\om(V,V[[z_1^{\pm 1},z_2^{\pm 1},u_1,\ldots ,u_n]])$ can be expressed as
\begin{gather}
A=B+ u_1^{a_1} C_1 +\ldots u_n^{a_n} C_n +h^k C_{n+1} \qquad\text{for some}\label{decompmod}\\
B\in\om(V,V((z_1,z_2))[[u_1,\ldots ,u_n,h]]),\,\, C_1,\ldots , C_{n+1}\in\om(V,V[[z_1^{\pm 1},z_2^{\pm 1},u_1,\ldots ,u_n]]). \non
\end{gather}
To indicate the fact that $A $ possesses a decomposition  as in \eqref{decompmod}, we   write
\beq\label{asups}
A\in\om(V,V((z_1,z_2))[[u_1,\ldots ,u_n]] )\mod u_1^{a_1},\ldots, u_n^{a_n},h^k.
\eeq
Note that the substitution
\beq\label{bsupst}
B\big|_{z_1=\phi(z_2,z_0)}\big. =\jota\left( B(z_1,z_2,u_1,\ldots , u_n)\big|_{z_1=\phi(z_2,z_0)}\big.\right)
\eeq
is well-defined even though the  substitution  $\textstyle A\big|_{z_1=\phi(z_2,z_0)}\big.$ does not exist in general. 
In what follows, the substitution $z_1=\phi(z_2,z_0)$ is always understood as in \eqref{bsupst}, i.e. the given expression is expanded in nonnegative powers of the variable $z_0$.
In order to   simplify our notation, we  denote \eqref{bsupst} as
\beq\label{asupst}
A \big|_{z_1=\phi(z_2,z_0)}^{\modd u_1^{a_1},\ldots, u_n^{a_n},h^k}\big. =A(z_1,z_2,u_1,\ldots , u_n)\big|_{z_1=\phi(z_2,z_0)}^{\modd u_1^{a_1},\ldots, u_n^{a_n},h^k}\big. .
\eeq
The element $B$ as in \eqref{decompmod}  is clearly  unique modulo
$$
\sum_{i=1}^n u_i^{a_i}\om(V,V[[z_1^{\pm 1},z_2^{\pm 1},u_1,\ldots ,u_n]]) + h^k \om(V,V[[z_1^{\pm 1},z_2^{\pm 1},u_1,\ldots ,u_n]]).
$$

Let $\g_N =\gl_N, \sll_N$.
The following definition 
of $\phi$-coordinated $\Vc_{\hspace{-1pt}c}(\g_N)$-module 
is based on \cite[Def. 3.4]{Li1}, which we slightly modify in order to make it compatible with Etingof--Kazhdan's quantum vertex algebra theory;  see Remark  \ref{napomena} for more details.
\begin{defn}\label{phimod}
A {\em $\phi$-coordinated $\Vc_{\hspace{-1pt}c}(\g_N)$-module} is a pair $(W,Y_W)$ such that $W$ is a topologically free $\CC[[h]]$-module and $Y_W=Y_W(z)$ is  a $\mathbb{C}[[h]]$-module map
\begin{align*}
Y_W \colon \Vc_{\hspace{-1pt}c}(\g_N)\ot W&\to W((z))[[h]]\\
u\ot w&\mapsto Y_W(z)(u\ot w)=Y_W(u,z)w=\sum_{r\in\mathbb{Z}} u_r w \ts z^{-r-1} 
\end{align*}
which satisfies
 $Y_W(\vac,z)w=w $ for all $w\in W$; 
the {\em weak associativity}: for any $u,v\in \Vc_{\hspace{-1pt}c}(\g_N)$ and $k\in\mathbb{Z}_{\geqslant 0}$ there exists $p\in\mathbb{Z}_{\geqslant 0}$ such that
\begin{align}
&(z_1-z_2)^p\ts Y_W(u,z_1)Y_W(v,z_2)\in\om (W,W((z_1,z_2)) )\mod h^k\Fand\label{associativitymod0}\\
&\big((z_1-z_2)^p\ts Y_W(u,z_1)Y_W(v,z_2)\big)\big|_{z_1=\phi(z_2,z_0)}^{\modd h^k}  \big. \non\\
&\qquad- (\phi(z_2,z_0) -z_2)^p\ts Y_W\left(Y(u,z_0)v,z_2\right)\ts
\in\ts  h^k \om(W,W[[z_0^{\pm 1},z_2^{\pm 1}]]);\label{associativitymod}
\end{align}
and
the {\em $\wht{\Sc}$-locality}:
for any $u,v\in \Vc_{\hspace{-1pt}c}(\g_N)$ and $k\in\mathbb{Z}_{\geqslant 0}$ there exists $p\in\mathbb{Z}_{\geqslant 0}$ such that    
\begin{align}
&(z_1-z_2)^{p}\ts Y_W(z_1)\big(1\otimes Y_W(z_2)\big)\iotaopjd\big(\wht{\mathcal{S}}(z_1/z_2)(u\otimes v)\otimes w\big)\label{localitymod}\\
&\qquad-(z_1-z_2)^{p}\ts Y_W(z_2)\big(1\otimes Y_W(z_1)\big)(v\otimes u\otimes w)\ts
\in\ts h^k W[[z_1^{\pm 1},z_2^{\pm 1}]] \quad\text{for all }w\in W.\non
\end{align}

Let $W_1$ be a topologically free $\CC[[h]]$-submodule of $W$. A pair $(W_1,Y_{W_1})$ is said to be a {\em $\phi$-coordinated $\Vc_{\hspace{-1pt}c}(\g_N)$-submodule} of $W$ if $Y_W(v,z)w$ belongs to $W_1$ for all $v\in \Vc_{\hspace{-1pt}c}(\g_N)$ and $w\in W_1$, where $Y_{W_1}$ denotes the restriction and corestriction of $Y_W$,
$$Y_{W_1} (z)=Y_{W} (z)\Big|_{\Vc_{\hspace{-1pt}c}(\g_N)\ot W_1}^{ W_1((z))[[h]]}\Big.\colon \Vc_{\hspace{-1pt}c}(\g_N)\ot W_1 \to W_1 ((z))[[h]].$$
\end{defn}
   
\begin{rem}\label{assocremark}
Regarding the weak associativity,  note that  \eqref{associativitymod0} and \eqref{associativitymod} employ  the notation introduced in \eqref{asups} and \eqref{asupst} for $n=0$, i.e.   there are no variables $u_1,\ldots ,u_n$.
Next, observe that the $\wht{\Sc}$-locality already implies that there exists $p\in\mathbb{Z}_{\geqslant 0}$ such that \eqref{associativitymod0} holds. However, we still include this requirement in the definition as it ensures that the  integer $p$ is   large enough so that the substitution $z_1=\phi(z_2,z_0)$ in \eqref{associativitymod} is well-defined. Finally, the motivation for expressing the weak associativity in the   form as in \eqref{associativitymod0} and \eqref{associativitymod} is given in \cite[Rem. 3.2]{Li1}.
\end{rem}

\begin{rem}\label{napomena}
As with the quantum affine algebra in the previous section, we introduce the notion of $\phi$-coordinated module over the ring $\CC[[h]]$ instead of a field in order to make it compatible with the Etingof--Kazhdan quantum vertex algebra theory; cf. original definition \cite[Def. 3.4]{Li1}.
Furthermore, unlike the original definition,   we require that the $\phi$-coordinated module map $Y_W(z)$  possesses    $\wht{\Sc}$-locality property \eqref{localitymod}. The general theory developed by Li  suggests that  \eqref{localitymod} might be omitted from the definition, due to the fact that the vertex operator map $Y(z)$ already possesses   $\Sc$-locality property \eqref{locality}; see \cite[Prop. 5.6]{Li1}. However, we  include the $\wht{\Sc}$-locality  in the definition in order to emphasize the importance of quantum current commutation relation \eqref{qc}. More specifically, in the proof of the \hyperref[mainthm1]{Main Theorem}, we   derive the $\wht{\Sc}$-locality property  directly from the quantum current commutation relation; see Lemma \ref{lemlem4}.
\end{rem}

Introduce the series
$$\delta(z)=\sum_{k\in\ZZ} z^k\in\CC[[z^{\pm 1}]]\fand 
\log(1+z)=-\sum_{k=1}^{\infty} \frac{(-z)^{k}}{k}\in z\CC[[z]].
$$
The following Jacobi-type identity was  established in \cite[Prop. 5.9]{Li1}. Although, in contrast with \cite{Li1}, we consider quantum vertex algebras and $\phi$-coordinated modules defined over the ring $\CC[[h]]$, the next proposition can be proved by  arguing as in the proofs of \cite[Lemma 5.8]{Li1} and \cite[Prop. 5.9]{Li1}.

\begin{pro}\label{Jacobi_prop}
Let $W$ be a $\phi$-coordinated $\Vc_{\hspace{-1pt}c}(\g_N)$-module, where $\phi(z_2,z_0)=z_2 e^{z_0}$. 
For any $u,v\in\Vc_{\hspace{-1pt}c}(\g_N)$  we have
\begin{align}
&(z_2 z)^{-1} \delta\left(\frac{z_1 -z_2}{z_2 z}\right) Y_W(z_1) (1 \ot Y_W(z_2))(u\ot v) \label{Jacobi1}\\
&\qquad -(z_2 z)^{-1} \delta\left(\frac{z_2-z_1}{-z_2 z}\right) Y_W(z_2)(1\ot Y_W(z_1))  \iotaopdj \wht{\Sc}(z_2 /z_1)(v\ot u)\label{Jacobi2}\\
=&\, z_1^{-1}\delta\left(\frac{z_2 (1+z)}{z_1}\right) Y_W\left(Y(u,\log(1+z))v,z_2\right).\label{Jacobi3}
\end{align}
\end{pro}

%%%%%%%%%%%%%%%%%%%%%%%%%%%%%%%%%%%
%%%%%%%%%%%%%%%%%%%%%%%%%%%%%%%%%%%
\section{Proof of the \texorpdfstring{\hyperref[mainthm1]{Main Theorem}}{Main Theorem}}\label{sec05}
%%%%%%%%%%%%%%%%%%%%%%%%%%%%%%%%%%%
%%%%%%%%%%%%%%%%%%%%%%%%%%%%%%%%%%%

In this section we prove the \hyperref[mainthm1]{Main Theorem}.
The proof   is   divided into four parts, Subsections \ref{sec051}--\ref{sec054}. 
In Subsection \ref{sec051}, we obtain some properties of the normalizing functions for the trigonometric $R$-matrix   which are required in the later stages of the proof; see Lemmas \ref{arofexlemma}--\ref{lemma4}.
In Subsection \ref{sec052}, we demonstrate how to establish the $\phi$-coordinated $\Vc_{\hspace{-1pt}c}(\gl_N)$-module structure on a restricted module of level $c$ for the quantum affine algebra $\Uh$; see Lemmas \ref{lemlem1}--\ref{lemlem4}. The key ingredient in this part of the proof is Ding's quantum current realization and, in particular, the fact that quantum current commutation relation \eqref{qc} resembles   $\wht{\Sc}$-locality property \eqref{localitymod}. In Subsection \ref{sec053}, we use Li's Jacobi-type identity, as given in Proposition \ref{Jacobi_prop}, to   establish the structure of restricted module of level $c$ for the quantum affine algebra $\Uh$ on a 
$\phi$-coordinated $\Vc_{\hspace{-1pt}c}(\gl_N)$-module; see Lemma \ref{lemlem5}.
Finally, we finish the proof in the $\g_N=\gl_N$ case by showing that the $\CC[[h]]$-submodules invariant with respect to the action of the quantum affine algebra and with respect to  the corresponding action of the quantum vertex algebra coincide; see Lemma \ref{lemlem6}.
In Subsection \ref{sec054}, we  use the fusion procedure for the two-parameter trigonometric $R$-matrix to extend the  results to the $\g_N=\sll_N$ case, thus completing the proof of the \hyperref[mainthm1]{Main Theorem}; see Lemmas \ref{profuzion}--\ref{lemlemE}.

%%%%%%%%%%%%%%%%%%%%%%%%%%%%%%%%%%%
\subsection{Normalizing functions  }\label{sec051}
%%%%%%%%%%%%%%%%%%%%%%%%%%%%%%%%%%%

Introduce the function $r(x)$ by
\beq\label{arofexrr31}
r(x)=- xe^h (1-e^{h}x)^{-1} f(x)^{-1},
\eeq
where   $f(x)$ is given by \eqref{f3}. 
\begin{lem}\label{arofexlemma}
The   function $r(x)\in\CC[[x,h]]$ satisfies
\beq\label{arofexrr}
R_{21}(x)^{-1}=r(x) R_{12}^+(1/x).
\eeq
Moreover, it admits the presentation
\beq\label{arofex}
r(x)=\sum_{k=0}^{\infty} r_k \frac{x^{k+1}}{(1-x)^{k+1}} \quad\text{such that}\quad r_{k}\in h^k \CC[[h]]\fand r_0=-e^{h}.
\eeq
\end{lem}

\begin{prf}
By combining unitarity property \eqref{unitrig} and \eqref{R2} we obtain
\begin{align*}
&R_{21}(x)^{-1}=\left(f(x)\R_{21}(x)\right)^{-1}=f(x)^{-1}\R_{21}(x)^{-1}
=f(x)^{-1} \R_{12}(1/x)\\
= &\,f(x)^{-1} (1-e^{-h}x^{-1})^{-1} R_{12}^+ (1/x)
=-xe^h (1-e^{h}x)^{-1} f(x)^{-1}  R_{12}^+ (1/x) =r(x) R_{12}^+(1/x),
\end{align*}
as required.
Next, by using \eqref{f3} we find   
\beq\label{lemaf}
f(x)^{-1}=\sum_{l=0}^{\infty} \left(-\sum_{k=1}^{\infty} f_k\left( \frac{x}{1-x}\right)^k \right)^l=1+\sum_{k=1}^{\infty} \beta_k\left( \frac{x}{1-x}\right)^k
\eeq
for some $\beta_k\in h^k\CC[[h]]$. 
It is clear that the product of \eqref{exp} for $a=1$, \eqref{lemaf} and $-xe^h $  is  equal to   $r(x)$ and, furthermore, that it admits presentation \eqref{arofex}.
\end{prf}

We use the following   lemma  in the proofs of weak associativity and  $\wht{\Sc}$-locality of the $\phi$-coordinated module map, as well as to establish the   restricted module structure on a $\phi$-coordinated $\Vc_{\hspace{-1pt}c}(\gl_N)$-module; see Lemmas \ref{lemlem3}, \ref{lemlem4} and \ref{lemlem5} respectively.

\begin{lem}\label{poleslemma}
Let $F=g^{\pm 1}$ or $F=r^{\pm 1}$. For any  integers $a_1,a_2,k> 0$ and $\alpha\in\CC$ there exists an integer $p\geqslant 0$ such that the coefficients of all monomials 
\beq\label{monomi}
u_1^{a_1'} u_2^{a_2'} h^{k'}, \qquad\text{where}\qquad 0\leqslant a_1'< a_1,\quad 0\leqslant a_2'< a_2\fand 0\leqslant k'< k,
\eeq
 in
$ (z_1 - z_2)^p F(z_1 e^{u_1 -u_2 +\alpha h} / z_2)$ 
belong to $\CC[z_1 , z_2^{\pm 1}]$ and such that the coefficients of all monomials \eqref{monomi} in
\beq\label{compare}
\left((z_1 - z_2)^p F(z_1 e^{u_1 -u_2 +\alpha h} / z_2)\right)\big|_{z_1=z_2 e^{z_0}}^{\modd u_1^{a_1}, u_2^{a_2}, h^{k}}\big. 
\fand   
z_2^p (e^{z_0}-1)^p  F(e^{z_0+u_1-u_2+\alpha h})
\eeq
coincide.
\end{lem}

\begin{prf}
Set  $\delta=0$ for $F=g$ and $\delta=1$ for $F=r$, i.e. $\delta=\delta_{F,r}$, so that we can consider both cases simultaneously. 
 Let $U=\CC[[x,x_0,h]]$.
Recall \eqref{g} and \eqref{arofex}.
As the map $\iotaopx$ commutes with partial differential operator $\partial / \partial x$,  by using    Taylor Theorem  \eqref{taylor} we find
$$
\iotazzo F(x+x_0) 
=  \sum_{s=0}^{\infty}   \frac{x_0^l}
{l!} \frac{\partial^l}{\partial x^l} \iotaopx F(x)
= \sum_{l,s=0}^{\infty} \frac{F_s\ts  x_0^l}{l!} \iotaopx
\frac{\partial^l}{\partial x^l}  \left(\frac{x^{s+\delta}}
{\left(1-x\right)^{s+1}}\right)\in U,
$$
 where $F_s =g_s$ for $F=g$ and $F_s =r_s$ for $F=r$.
By \eqref{g} and \eqref{arofex}, every $F_s$ belongs to $h^s \CC[[h]]$, so all summands with $s\geqslant k$ are trivial modulo $h^k U$.
Hence the given expression modulo $U_0\coloneqq x_0^{a_1+a_2+k-2} U + h^k U$ contains only finitely many nonzero summands and, consequently, only  finitely many   terms   $(1-x)^{s+1}$ in the denominator. Therefore, there exists an integer $p\geqslant 0$ such that
$$
 \iotazzo (1-x)^p F(x+x_0) =  \sum_{l,s=0}^{\infty} \frac{F_s\ts  x_0^l}{l!} \iotaopx\ts
  (1-x)^p
\frac{\partial^l}{\partial x^l}  \left(\frac{x^{s+\delta}}
{\left(1-x\right)^{s+1}}\right) \in  \CC[x,x_0,h]\mod U_0 ,
$$
where the equality holds modulo $U_0$ and
the map $\iotaopx$ can be omitted on the right hand side as $p$ can be chosen so that  
$(1-x)^p$ cancels all negative powers of $(1-x)$ modulo $U_0$.
By applying the substitution
$  (x,x_0)=(z_1 /z_2 ,  z_1 (e^{u_1-u_2+\alpha h}-1) /z_2 ) $
to 
\beq\label{modulodp}
\iotazzo (1-x)^p F(x+x_0)\mod U_0
\eeq
 and then
multiplying the resulting expression by $(-z_2)^p$ we get
\beq\label{drugiput2}
 (z_1 - z_2)^p F(z_1 e^{u_1 -u_2 +\alpha h} / z_2)\in\CC[z_1,z_2^{\pm 1},u_1,u_2,h]\mod V_0 
\eeq
for $V_0=u_1^{a_1}V + u_2^{a_2}V + h^{k}V$ and $V=\CC[[z_1,z_2^{\pm 1},u_1,u_2,h]]$,
thus proving the first assertion of the lemma.

Set $W_0 =u_1^{a_1}W + u_2^{a_2}W + h^{k}W$ for $W=\CC [[z_0, z_2 ,u_1,u_2,h]]$. As \eqref{modulodp}  is a polynomial in the variables $x$ and $x_0$, by applying the substitution 
$
(x,x_0)=(e^{z_0}, e^{z_0}(e^{u_1 -u_2 +\alpha h} -1))
$
to \eqref{modulodp} 
and then multiplying the resulting expression by  $(-z_2)^p$ we get
\beq\label{drugiput3}
z_2^p (e^{z_0}-1)^p  F(e^{z_0+u_1-u_2+\alpha h})\mod W_0,
\eeq
where, by the expansion convention from Subsection \ref{sec01}, $F(e^{z_0+u_1-u_2+\alpha h})$  stands for $\iotaztri F(e^{z_0+u_1-u_2+\alpha h})$.\footnote{Note that the expression $\iotazzo (1-x)^p F(x+x_0)$ is considered modulo $U_0$ because, otherwise, the aforementioned substitution would not be well-defined (although the same substitution is well-defined when applied to $(1-x)^p F(x+x_0)$ with $F(x+x_0)$ being regarded as a rational function with respect to the variables $x$ and $x_0$).}
Finally, as \eqref{drugiput2} modulo $V_0$ is a polynomial with respect to the variables $z_1 /z_2$ and $z_2$, by applying the substitution $z_1=z_2 e^{z_0}$  we again  obtain \eqref{drugiput3}, thus proving the second assertion of the lemma.

If $F=g^{-1}$ or $F=r^{-1}$, one  easily checks that 
\beq\label{onemoreref}
F(x)= \sum_{s=0}^\infty F_s \frac{x^{s-\delta_{ F,r^{-1}}}}{(1-x)^{s-1}}\quad\text{for some}\quad F_s\in h^s\CC[[h]],
\eeq
so the lemma is verified by arguing as above.
\end{prf}

We now recall a certain useful  consequence of \cite[Lemma 2.7]{Li1},  as given in \cite[Rem. 2.8]{Li1}:  
For any $A(z_1,z_2),B(z_1,z_2)\in\CC((z_1,z_2))$, the equality
\beq\label{rem28}
A(z_1,z_2)\big|_{z_1=z_2 e^{z_0}}\big. =B(z_1,z_2)\big|_{z_1=z_2 e^{z_0}}\big.
\quad\text{implies}\quad
A(z_1,z_2)  =B(z_1,z_2).
\eeq
Since the $\CC[[h]]$-module $\CC((z_1,z_2))[[h]]$ is separable, implication \eqref{rem28} clearly extends to  any
 $A(z_1,z_2),B(z_1,z_2)\in\CC((z_1,z_2))[[h]]$.

\begin{lem}\label{lemma3}
In $\CC((u))[[h]]$ we have
\beq\label{arofexg}
 r(e^{-u}) =\psi^2 g(e^u). 
\eeq
Moreover, for any  integers $a_1,a_2,k> 0$ and $\alpha\in\CC$ there exists an integer $p\geqslant 0$ such that the coefficients of all monomials 
\eqref{monomi}
 in
\beq\label{arofexg2}
(z_1 - z_2)^p \ts r(z_2 e^{-u_1+u_2-\alpha h}/z_1)\Fand (z_1 - z_2)^p\ts \psi^2\ts g(z_1 e^{u_1-u_2+\alpha h}/z_2)
\eeq
coincide.
\end{lem}

\begin{prf}
By \cite[Prop. 2.1]{KM} we have $ \psi^2 f(e^u)= f(e^{-u})^{-1}$. 
Therefore, using \eqref{arofexrr31}  we get
\begin{align*}
 r(e^{-u})&= -  e^{-u+h}(1-e^{-u+h})^{-1}f(e^{-u})^{-1}=- \psi^2  e^{-u+h}(1-e^{-u+h})^{-1}f(e^{u})\\
&= \psi^2 (1-e^{u-h})^{-1}f(e^{u})=\psi^2 g(e^u),
\end{align*}
as required, where the last equality follows from  \eqref{R2}.
Next,
by Lemma \ref{poleslemma} and  \eqref{arofexg}, there exists   $p\geqslant 0$    such that the coefficients of all monomials \eqref{monomi} in both expressions in  \eqref{arofexg2} belong to $\CC[z_1^{\pm 1},z_2^{\pm 1}] $ 
and such that the coefficients of all monomials \eqref{monomi} in
\begin{align*}
&\left( (z_1 - z_2)^p r(z_2 e^{-u_1+u_2 -\alpha h}/z_1)\right)\big|_{z_1=z_2 e^{z_0}}^{\modd u_1^{a_1}, u_2^{a_2}, h^{k}}\big.
\Fand\\
& \left( (z_1 - z_2)^p \psi^2 g(z_1 e^{u_1 -u_2+\alpha h}/z_2)\right)\big|_{z_1=z_2 e^{z_0}}^{\modd u_1^{a_1}, u_2^{a_2}, h^{k}}\big. 
\end{align*}
coincide. The second assertion  of the lemma now follows by   implication   \eqref{rem28}.
\end{prf}

The next lemma, which relies on Lemma \ref{lemma3}, will be used  in the proof of $\wht{\Sc}$-locality of the $\phi$-coordinated module map; see Lemma \ref{lemlem4}.

\begin{lem}\label{lemma4}
{\rm (1)} Let $F=g^{\pm 1}$ or $F=r^{\pm 1}$.
There exists  $\wht{F}(x,u,v)$ in $\CC(x)[[u,v,h]]$ such that for all $\alpha\in\CC$ the following equality in $\CC((z))[[u,v,h]]$ holds:  
\beq\label{fzuv}  
 \wht{F}(e^z,u,v-\alpha h)=  F(e^{z+u-v+\alpha h}).
\eeq
{\rm (2)}  For any integers $n,m>0$, the families of variables $u=(u_1,\ldots ,u_n)$, $v=(v_1,\ldots ,v_m)$ and $c\in\CC$ there exist  functions
$ \wht{G}(x,u,v), \wht{H}(x,u,v) \in\CC(x)[[u_1,\ldots ,u_n,v_1,\ldots ,v_m,h]]$ 
such that the following equalities  hold in $\CC((z))[[u_1,\ldots ,u_n,v_1,\ldots ,v_m,h]]$: 
\begin{align}
 &\qquad\qquad \wht{G}(e^z,u,v)=
 G(z,u,v)\fand  
\wht{H}(e^z,u,v)=
 H(z,u,v),\qquad\text{where}\non\\
 &G(z,u,v)=\prod_{i=1}^n \prod_{j=1}^m  
g(e^{z+u_i-v_j-h(N+c)})^{-1}
g(e^{z+u_i-v_j+hc})^{-1}
g(e^{z+u_i-v_j})^2 ,\label{gzuv}\\
 &H(z,u,v)=\prod_{i=1}^n \prod_{j=1}^m  
g(e^{-z+u_i-v_j-h(N+c)})^{-1} 
r(e^{z-u_i+v_j-hc})^{-1} 
g(e^{-z+u_i-v_j}) 
r(e^{z-u_i+v_j}).\label{fzuv2}
\end{align}
{\rm (3)} Let $a_1,\ldots ,a_n,b_1,\ldots ,b_m ,k>0$ be arbitrary integers and $\iotaop=\iotaall$ the embedding. There exists an integer $p\geqslant 0$ such that the   coefficients of all monomials 
$$
u_1^{a'_1}\ldots u_n^{a'_n} v_1^{b'_1}\ldots v_{m}^{b'_m} h^{k'},\qquad \text{where}\qquad 0\leqslant a'_{i} < a_i,\, \quad 0\leqslant b'_j < b_j \fand 0\leqslant k'< k
$$
in
$
(z_1 -z_2)^p  \iotaop \wht{G}(z_1 /z_2,u,v)$ and
$(z_1 -z_2)^p \iotaop \wht{H}(z_2 /z_1,u,v)$
coincide.
\end{lem}

\begin{prf}
Due to \eqref{g}, \eqref{arofex} and \eqref{onemoreref}, we can regard $g(x)^{\pm 1}$ and $r(x)^{\pm 1}$ as elements of $\CC(x)[[h]]$. 
Let $F=g^{\pm 1}$ or $F=r^{\pm 1}$ and write
$F(x)=\sum_{s=0}^{\infty} F_s(x) h^s$
for some $F_s(x)\in \CC(x)$.
Applying formal Taylor Theorem \eqref{taylor} to $z\mapsto  \iotaopz F_s(e^{z})$ we get for any $\alpha\in\CC$ 
$$
\iotazuvvh F_s(e^{z+u-v+\alpha h}) =\sum_{l= 0}^\infty \frac{(u-v+\alpha h)^l}{l!}\frac{\partial^l}{\partial z^l} \iotaopz F_s (e^z)\quad\text{in}\quad \CC((z))[[u,v,h]].
$$
The partial differential operator $\partial/\partial z$ commutes with the map $\iotaopz$ and   all $\frac{\partial^l}{\partial z^l}  F_s (e^z)$   can be naturally regarded as elements of $\CC(e^z)$. Hence we can introduce  functions $\wht{F}_{ l,s}(x)\in\CC(x)$ by the requirement $\wht{F}_{ l,s}(e^z)=\frac{\partial^l}{\partial z^l}F_s (e^z) $. The first statement of the lemma now clearly follows as the  function $\wht{F}(x,u,v)\in\CC(x)[[u,v,h]]$ satisfying \eqref{fzuv} can be defined by
$$\wht{F}(x,u,v)=\sum_{ s=0}^{\infty}\left(\sum_{ l=0}^{\infty}\frac{(u-v)^l  }{l!}   \wht{F}_{l,s}(x)\right)h^s .$$
The second statement  is proved by applying the first statement   on each factor of \eqref{gzuv} and \eqref{fzuv2}. Finally, by \eqref{arofexg} we have $G(z,u,v)=H(-z,u,v)$, so the third statement  follows by   Lemma \ref{lemma3}.
\end{prf}

%%%%%%%%%%%%%%%%%%%%%%%%%%%%%%%%%%%
\subsection{Establishing the   \texorpdfstring{$\phi$}{phi}-coordinated \texorpdfstring{$\Vc_{\hspace{-1pt}c}(\gl_N)$}{Vc(glN)}-module structure}\label{sec052}
%%%%%%%%%%%%%%%%%%%%%%%%%%%%%%%%%%%

Let $W$ be a  restricted $\Uh$-module of level $c\in \CC$. In this subsection, we prove the first assertion of the \hyperref[mainthm1]{Main Theorem}, i.e. we show  that \eqref{formula} defines a unique
 structure of $\phi$-coordinated $\Vc_{\hspace{-1pt}c}(\gl_N)$-module on $W$, where $\phi(z_2,z_0) = z_2 e^{z_0}$.
The proof is divided into four lemmas which  verify all requirements imposed  by Definition \ref{phimod}.

\begin{lem}\label{lemlem1}
Formula \eqref{formula}, together with $Y_W(\vac,z)=1_W$, defines a unique $\CC[[h]]$-module map $\Vc_{\hspace{-1pt}c}(\gl_N)\ot W \to W((z))[[h]]$.
\end{lem}

\begin{prf}
First, we note that the right hand side of \eqref{formula} is well-defined, as was discussed in Remark \ref{linremark}. Next, we recall that the algebra $\Uhplus$ is spanned by all coefficients of all matrix entries of $T_{[n]}^+(u)$, $n\geqslant 1$, and $\vac$; see \cite[Sect. 3.4]{EK3} or \cite[Prop. 2.4]{KM}\footnote{
We should mention that the notation in this paper slightly differs from \cite{KM}. In particular, the algebra $\textrm{U}(R)$, as defined in \cite[Sect. 2]{KM}, coincides with the algebra $\Uhplus$ defined in  Subsection \ref{subsec012}.
}.  
In order to prove the lemma, we  have to show that $v\mapsto Y_W(v,z)$ preserves the ideal of relations \eqref{rtt}. More specifically, it is sufficient to check that for any integers $n\geqslant 2$ and $i=1,2,\ldots ,n-1$ and the family of variables $u=(u_1,\ldots ,u_n)$   the expression
\begin{align}
R_{i\ts i+1}(e^{u_i -u_{i+1}} )T_{[n]}^+ (u)\vac - P_{i\ts i+1}T_{[n]}^+ (u_{i+1,i})\vac P_{i\ts i+1}R_{i\ts i+1}(e^{u_i -u_{i+1}} ),\label{expr01}
\end{align}
where $u_{ i+1,i}=(u_1,\ldots ,u_{i-1},u_{i+1},u_i,u_{i+2},\ldots ,u_n)$, belongs to the kernel of $v\mapsto Y_W(v,z)$. 

Let $x=(x_1,\ldots ,x_n)$ and $x_{i+1,i}=(x_1,\ldots ,x_{i-1},x_{i+1},x_i,x_{i+2},\ldots ,x_n)$. Using Yang--Baxter equation \eqref{YBE} and   commutation relation \eqref{qc} one can prove
the identity
\beq\label{expr03}
R_{i\ts i+1} (x_i /x_{i+1})\Lc_{[n]}(x)=P_{i\ts i+1}\Lc_{[n]}(x_{ i+1,i})P_{i\ts i+1}R_{i\ts i+1} (x_i /x_{i+1}).
\eeq
 By Proposition \ref{restricted496}, all matrix entries of  $\Lc_{[n]}(x)$ belong to $\om(W,W((x_1,\ldots ,x_n))[[h]])$, so all matrix entries in \eqref{expr03} belong to 
$$\om(W,W((x_{i+1}))((x_1,\ldots ,x_i,x_{i+2},\ldots ,x_n))[[h]]).$$
Recall  $R$-matrix decomposition \eqref{R2}. By \eqref{onemoreref} the function   $g(x_i /x_{i+1})^{-1}$ belongs to $\CC[x_{i+1}^{-1}][[h,x_{i}]]$,
so 
 we can multiply \eqref{expr03} by   $g(x_i /x_{i+1})^{-1}$, 
 thus getting
\beq\label{expr04}
R_{i\ts i+1}^+ (x_i /x_{i+1})\Lc_{[n]}(x)=P_{i\ts i+1}\Lc_{[n]}(x_{ i+1,i})P_{i\ts i+1}R_{i\ts i+1}^+ (x_i /x_{i+1}).
\eeq
Since the $R$-matrix $R_{i\ts i+1}^+ (x_i /x_{i+1})$ is a polynomial in $x_i/x_{i+1}$, all matrix entries of both sides in \eqref{expr04} belong to 
$\om(W,W ((x_1,\ldots   ,x_n))[[h]])$.
Therefore, we can apply the substitutions $x_i =ze^{u_i}$ with $i=1,\ldots ,n$ to \eqref{expr04}, thus getting  the following equality in $(\ndo\CC^N)^{\ot n} \ot \om(W,W((z))[[h,u_1,\ldots ,u_n]])$:
$$
R_{i\ts i+1}^+ (e^{u_i -u_{i+1}}) \cdot \left(\Lc_{[n]}(x)\right)\big|_{x_i = ze^{u_i}}\big.
=P_{i\ts i+1}\left(\Lc_{[n]}(x_{i+1,i})\right)\big|_{x_i = ze^{u_i}}\big.\ts P_{i\ts i+1} \ts R_{i\ts i+1}^+  (e^{u_i -u_{i+1}}). 
$$
Multiplying the equality by $\psi g (e^{u_i -u_{i+1}})\in\CC((u_{i+1}))[[h,u_{i}]]$ and using \eqref{rplusg} we find 
$$
R_{i\ts i+1} (e^{u_i -u_{i+1}}) 
 \left(\Lc_{[n]}(x)\right)\big|_{x_i = ze^{u_i}}\big.
-P_{i\ts i+1}\left(\Lc_{[n]}(x_{i+1,i})\right)\big|_{x_i = ze^{u_i}}\big. \ts P_{i\ts i+1}\ts R_{i\ts i+1}  (e^{u_i -u_{i+1}})=0. 
$$
 As the left   hand side  coincides  with the image of   \eqref{expr01}, with respect to $Y_W(z)$,   we conclude that  \eqref{formula} defines a  $\CC[[h]]$-module map  $\Vc_{\hspace{-1pt}c}(\gl_N)\ot W \to W[[z ^{\pm 1}]]$, as required.
Moreover, by Remark \ref{linremark}  
 its image  belongs to $ W((z))[[h]]$.
Finally, it is clear that the $\CC[[h]]$-module map $Y_W(z)$ is uniquely determined by  \eqref{formula}.
\end{prf}

The next lemma  follows from $\wht{\Sc}$-locality property \eqref{localitymod} which is verified in Lemma \ref{lemlem4} below; recall Remark \ref{assocremark}. 
Nonetheless, we   provide the direct proof    as  the underlying calculations are required in the proof of Lemma \ref{lemlem3}.

\begin{lem}\label{lemlem2}
The map $Y_W(z)$ satisfies \eqref{associativitymod0}, i.e. for any  $u,v\in \Vc_{\hspace{-1pt}c}(\gl_N)$ and $k\in\mathbb{Z}_{\geqslant 0}$ there exists $p\in\mathbb{Z}_{\geqslant 0}$ such that
\beq\label{1associativitymod0}
(z_1-z_2)^p\ts Y_W(u,z_1)Y_W(v,z_2)\in\om (W,W((z_1,z_2)) )\mod h^k .
\eeq
\end{lem}

\begin{prf}
For any integers $n,m\geqslant 1$ and   families of variables $u=(u_1,\ldots ,u_n)$ and $v=(v_1,\ldots ,v_m)$ we have
\beq\label{exprr01}
Y_W(T_{[n]}^{+13}(u)\vac,z_1 )Y_W(T_{[m]}^{+23}(v)\vac,z_2 )
=\left(\Lc_{[n]}^{13} (x)\right)\big|_{x_i =z_1e^{u_i}}\big.  
\left(\Lc_{[m]}^{23} (y)\right)\big|_{y_j =z_2e^{v_j}}\big. ,
\eeq
where $x=(x_1,\ldots ,x_n)$ and $y=(y_1,\ldots ,y_m)$. The coefficients in \eqref{exprr01} are  operators on the multiple tensor product with superscripts   $1,2,3$ indicating the tensor factors:
$$
\smalloverbrace{(\ndo\CC^N)^{\ot n}}^{1}
\ot
\smalloverbrace{(\ndo\CC^N)^{\ot m}}^{2}
\ot 
\smalloverbrace{W}^{3}.
$$
Let us rewrite the right hand side in \eqref{exprr01}.
The third assertion of Proposition \ref{qcgenpro} implies
\beq\label{exprr03}
\Lc_{[n]}^{13}(x) R_{nm}^{21} (ye^{-hc} / x) \Lc_{[m]}^{23}(y)=\Lc_{[n+m]}(x,y)R_{nm}^{21} (y/ x).
\eeq
By expressing the second crossing symmetry   relation in \eqref{csym}  in  the variable $x=y_j e^{-h(N+c)}/x_i$, then  applying  the transposition $t_2$ and finally conjugating the resulting equality by the permutation operator $P$
we find
$$
\left(D_1 R_{21}(y_j e^{-h(N+c)}/x_i)^{-1} D_1^{-1} \right) \cdotrl R_{21} (y_j e^{-hc}/x_i) =1.
$$
Furthermore, due to  Lemma \ref{arofexlemma}, we can write this equality as
$$
  r(y_j e^{-h(N+c)}/x_i) \left( D_1 R_{12}^+ (x_i e^{h(N+c)}/y_j)  D_1^{-1} \right) \cdotrl R_{21} (y_j e^{-hc}/x_i) =1.
$$
Hence we have
\begin{gather}
r(x,y) 
\left( D_{[n]}^1  
R_{nm}^{+12}(x e^{h(N+c)}/y) 
\left(D_{[n]}^1\right)^{-1}\right)\cdotrl  R_{nm}^{21} (ye^{-hc} / x) =1,\qquad \text{where}\label{temp9}\\
 D_{[n]}^1  =D^{\ot n}\ot 1^{\ot m}
\fand
r(x,y)=\prod_{i=1}^n \prod_{j=1}^m  r(y_j e^{-h(N+c)}/x_i) .\label{rfunction} 
\end{gather}
Using \eqref{temp9} we can move $ R_{nm}^{21} (ye^{-hc} / x)$ in \eqref{exprr03} to the right hand side, which gives us
\begin{align}
&\Lc_{[n]}^{13}(x)  \Lc_{[m]}^{23}(y)=
r(x,y) 
\left(  D_{[n]}^1  
R_{nm}^{+12}(x e^{h(N+c)}/y)
\left(D_{[n]}^1\right)^{-1}\right)\cdotrl\left(\Lc_{[n+m]}(x,y) R_{nm}^{21} (y/ x)\right) 
\non\\
=&\, r(x,y) \ts g(x,y)
\left(  D_{[n]}^1  
R_{nm}^{+12}(x e^{h(N+c)}/y)
\left(D_{[n]}^1\right)^{-1}\right)\cdotrl\left(\Lc_{[n+m]}(x,y) R_{nm}^{+21} (y/ x)\right),
\label{exprr04}
\end{align}
where the second equality comes from \eqref{R2} and the function $g(x,y)$ is given by
\beq\label{gfunction} 
g(x,y)=\prod_{i=1}^n \prod_{j=1}^m   g(y_j /x_i) .
\eeq

Let $a_1,\ldots ,a_n,b_1,\ldots ,b_m ,k>0$ be arbitrary integers. We now apply the substitutions 
\beq\label{subs}
x_i=z_1 e^{u_i}, \quad y_j=z_2 e^{v_j}\qquad\text{for}\qquad i=1,\ldots ,n,\quad  j=1,\ldots ,m
\eeq
 to \eqref{exprr04}, thus getting \eqref{exprr01}, and then   consider the coefficients of all monomials 
\beq\label{monomials}
u_1^{a'_1}\ldots u_n^{a'_n} v_1^{b'_1}\ldots v_{m}^{b'_m} h^{k'},\qquad \text{where} \qquad 0\leqslant a'_{i} < a_i,\quad 0\leqslant  b'_j < b_j \fand k'< k.
\eeq
First, as the $R$-matrix $R^+(w)$ is a polynomial with respect to the variable $w$, we conclude 
by Proposition \ref{restricted496} and Remark \ref{linremark} that
\begin{align}
&\left(
\left( 
D_{[n]}^1  
R_{nm}^{+12}(x e^{h(N+c)}/y)
\left(D_{[n]}^1\right)^{-1}
\right)
\cdotrl
\left(
\Lc_{[n+m]}(x,y)R_{nm}^{+21} (y/ x)
\right)
\right) 
\bigg|_{\substack{x_i =z_1e^{u_i},\, y_j =z_2e^{v_j}}}\bigg. \label{prviclan}\\
&\in
(\ndo\CC^N)^{\ot n}\ot(\ndo\CC^N)^{\ot m}\ot\om(W,W((z_1,z_2))[[u_1,\ldots,u_n,v_1,\ldots ,v_n, h]]). \non
\end{align}
Next, by Lemma \ref{poleslemma}  there exists an integer $p\geqslant 0$, which depends on the choice of integers $a_1,\ldots ,a_n,b_1,\ldots ,b_m,k$, such that the coefficients of all monomials \eqref{monomials} in
\beq\label{drugiclan}
(z_1 -z_2)^p   \big( r(x,y) g(x,y)\big)\big|_{x_i =z_1e^{u_i},\,y_j =z_2e^{v_j}} \big.
\eeq
belong to $\CC[z_1^{\pm 1},z_2^{\pm 1}]$.
Finally, we observe that the coefficients of all monomials \eqref{monomials} in the product of \eqref{prviclan} and \eqref{drugiclan} coincide with the corresponding coefficients in
$$
(z_1 -z_2)^p \ts Y_W(T_{[n]}^{+13}(u)\vac,z_1 )Y_W(T_{[m]}^{+23}(v)\vac,z_2 ).
$$
Therefore, by the preceding discussion,   these coefficients belong to 
$$
(\ndo\CC^N)^{\ot n}\ot(\ndo\CC^N)^{\ot m}\ot\om(W,W((z_1,z_2)))
,$$ which implies the statement of the lemma.
\end{prf}

\begin{lem}\label{lemlem3}
The map $Y_W(z)$ satisfies weak associativity \eqref{associativitymod}, i.e.
 for any $u,v\in \Vc_{\hspace{-1pt}c}(\gl_N)$ and $k\in\mathbb{Z}_{\geqslant 0}$ there exists $p\in\mathbb{Z}_{\geqslant 0}$ such that \eqref{1associativitymod0} holds and such that
\begin{align}
&\big((z_1-z_2)^p\ts Y_W(u,z_1)Y_W(v,z_2)\big)\big|_{z_1= z_2 e^{z_0}}^{\modd h^k}  \big. \non\\
&\qquad- z_2^p ( e^{z_0} -1)^p\ts Y_W\left(Y(u,z_0)v,z_2\right)\ts
\in\ts  h^k \om(W,W[[z_0^{\pm 1},z_2^{\pm 1}]]).\label{1associativitymod}
\end{align}
\end{lem}

\begin{prf}
Let $n,m,a_1,\ldots ,a_n,b_1,\ldots ,b_m ,k>0$ be arbitrary integers, $u=(u_1,\ldots ,u_n)$ and $v=(v_1,\ldots ,v_m)$ the families of variables.
Consider the coefficients of all monomials \eqref{monomials} in the expression
\beq\label{left}
\big((z_1-z_2)^p\ts Y_W(T_{[n]}^{+13} (u)\vac,z_1)Y_W(T_{[m]}^{+23} (v) \vac,z_2)\big)\big|_{z_1= z_2 e^{z_0}}^{\modd u_1^{a_1},\ldots , u_n^{a_n},v_1^{b_1},\ldots ,v_m^{b_m},h^k}\big.,
\eeq
which corresponds to the first summand in  \eqref{1associativitymod}. As   demonstrated in the proof of Lemma \ref{lemlem2}, they coincide with the coefficients of all monomials \eqref{monomials} in the product
\begin{align}
&\left( \hspace{-2pt}\left( \hspace{-1pt}
\left(  
D_{[n]}^1  
R_{nm}^{+12}(x e^{h(N+c)}/y)
\left(D_{[n]}^1\right)^{-1}
\right)
\cdotrl
\left(
\Lc_{[n+m]}(x,y) R_{nm}^{+21} (y/ x)
\right)\hspace{-1pt}
\right) 
\Big|_{\substack{x_i =z_1e^{u_i},\,y_j =z_2e^{v_j}}}\Big.\right) 
\bigg|_{z_1= z_2 e^{z_0}}\bigg. \label{tren1}\\
& \times
\left( (z_1 -z_2)^p \cdot \big( r(x,y) g(x,y)\big)\big|_{\substack{x_i =z_1e^{u_i},\,y_j =z_2e^{v_j}}}\big. \right) \Big|_{z_1= z_2 e^{z_0}}^{\modd u_1^{a_1},\ldots , u_n^{a_n},v_1^{b_1},\ldots ,v_m^{b_m},h^k}\Big. \label{tren2}
\end{align}
for a suitably chosen integer $p\geqslant 0$ (which depends on $a_1,\ldots ,a_n,b_1,\ldots ,b_m ,k$).  Recall that the functions $r$ and $g$ are given by \eqref{rfunction} and \eqref{gfunction}.
First, we observe that the coefficients of all monomials \eqref{monomials}  in  factor \eqref{tren1} coincide with the corresponding coefficients in
\beq\label{tren3}
\left(  
D_{[n]}^1  
R_{nm}^{+12}(e^{z_0 +u-v+h(N+c)})
\left(D_{[n]}^1\right)^{-1}
\right)
\cdotrl
\left(\hspace{-1pt}
\bigg(
\Lc_{[n+m]}(x,y)
\Big|_{\substack{x_i =z_1e^{u_i}\\y_j =z_2e^{v_j}}}\Big. 
\bigg)
\bigg|_{z_1= z_2 e^{z_0}}\bigg.
R_{nm}^{+21} (e^{-z_0-u+v})
\hspace{-1pt}\right)\hspace{-1pt}\hspace{-1pt}.
\eeq
Next, we turn to   factor \eqref{tren2}.  Due to  Lemma \ref{poleslemma}, we can assume that the integer $p$ was chosen so that the coefficients of all monomials \eqref{monomials} in factor \eqref{tren2} coincide with the coefficients of all monomials \eqref{monomials} in
$$
z_2^p (e^{z_0} -1)^p \prod_{i=1}^n \prod_{j=1}^m r(e^{-z_0-u_i+v_j-h(N+c)}) g(e^{-z_0-u_i+v_j}).
$$
Moreover, by \eqref{arofexg}, this is equal to
\beq\label{left3}
z_2^p (e^{z_0} -1)^p \psi^{2nm} \prod_{i=1}^n \prod_{j=1}^m g(e^{z_0+u_i-v_j+h(N+c)}) g(e^{-z_0-u_i+v_j}).
\eeq
Finally, we conclude that the coefficients of all monomials \eqref{monomials} in \eqref{left} coincide with the coefficients of the corresponding monomials in the product of 
\eqref{tren3} and \eqref{left3}.

Consider the expression
\beq\label{right1}
z_2^p (e^{z_0} -1)^p \ts Y_W(Y(T_{[n]}^{+13} (u)\vac,z_0)T_{[m]}^{+23} (v)\vac,z_2),
\eeq
which corresponds to the second summand in \eqref{1associativitymod}. By \eqref{qva1} it is equal to
\beq\label{exprt01}
z_2^p (e^{z_0} -1)^p \ts Y_W( T_{[n]}^{+13} (u|z_0)T_{[n]}^{*13} (u|z_0+hc/2)^{-1} T_{[m]}^{+23} (v)\vac,z_2).
\eeq
Since $T^*(u)\vac =\vac$, by combining relation \eqref{rtt3} and the first  crossing symmetry relation in \eqref{csym_equiv} we obtain
\begin{align}
&T_{[n]}^{*13} (u|z_0+hc/2)^{-1} T_{[m]}^{+23} (v)\vac\non \\
&\qquad=
\left(D_{[n]}^{1} R_{nm}^{12}(e^{z_0+u-v +h(N+c)}) (D_{[n]}^{1})^{-1}\right)\cdotrl\left(T_{[m]}^{+23} (v)\vac R_{nm}^{12}(e^{z_0+u-v})^{-1} \right).\label{ge3}
\end{align}
Introduce the functions
$$
g_1(u,v,z_0)=\psi^{nm}\prod_{i=1}^n \prod_{j=1}^m g(e^{z_0+u_i-v_j+h(N+c)}),
\quad
g_2(u,v,z_0)=\psi^{nm} \prod_{i=1}^n \prod_{j=1}^m g(e^{-z_0-u_i+v_j}).
$$
By \eqref{rplusg} we have
\beq\label{ge1}
R_{nm}^{12}(e^{z_0+u-v +h(N+c)})=g_1(u,v,z_0)R_{nm}^{+12}(e^{z_0+u-v +h(N+c)}).
\eeq
Furthermore, by combining  \eqref{rplusg} and unitarity property \eqref{uni}    we find
\beq\label{ge2}
R_{nm}^{12}(e^{z_0+u-v  })^{-1}=R_{nm}^{21}(e^{-z_0-u+v })=g_2(u,v,z_0)R_{nm}^{+21}(e^{-z_0-u+v }).
\eeq
Using \eqref{ge1} and \eqref{ge2} we rewrite the right hand side of \eqref{ge3} as
\begin{align}
&g_1(u,v,z_0)g_2(u,v,z_0) \left(D_{[n]}^{1} R_{nm}^{+12}(e^{z_0 +u-v+h(N+c)}) (D_{[n]}^{1})^{-1}\right)\cdotrl\left(T_{[m]}^{+23} (v)\vac R_{nm}^{+21}(e^{-z_0-u+v }) \right).\label{ge4}
\end{align}
Next, we employ  \eqref{ge4} and then \eqref{formula}  to express \eqref{exprt01} as
\begin{align}
&\left(D_{[n]}^{1} R_{nm}^{+12}(e^{z_0 +u-v +h(N+c)}) (D_{[n]}^{1})^{-1}\right)\cdotrl\Big(Y_W(T_{[n]}^{+13} (u|z_0)T_{[m]}^{+23} (v)\vac,z_2)  R_{nm}^{+21}(e^{-z_0-u+v })\Big)\non\\
&\qquad\times z_2^p (e^{z_0} -1)^p   \ts  g_1(u,v,z_0)\ts g_2(u,v,z_0)\non\\
=&\left(D_{[n]}^{1} R_{nm}^{+12}(e^{z_0 +u-v +h(N+c)}) (D_{[n]}^{1})^{-1}\right)\cdotrl
\Big(\Lc_{[n+m]} (x,y)
\Big|_{\substack{x_i =z_2e^{z_0+u_i}\\\hspace{-13pt}y_j =z_2e^{v_j}}}\Big. \ts
  R_{nm}^{+21}(e^{-z_0-u+v })\Big)\non\\
 &\qquad\times z_2^p (e^{z_0} -1)^p   \ts g_1(u,v,z_0)\ts g_2(u,v,z_0).\non
\end{align}
Note that $z_2^p (e^{z_0} -1)^p g_1(u,v,z_0)g_2(u,v,z_0)$ is equal to \eqref{left3} and that
$$
 \Lc_{[n+m]} (x,y)
\Big|_{\substack{x_i =z_2e^{z_0+u_i}\\\hspace{-13pt}y_j =z_2e^{v_j}}}\Big.  
=
\left(
\Lc_{[n+m]}(x,y)
\Big|_{\substack{x_i =z_1e^{u_i}\\y_j =z_2e^{v_j}}}\Big. 
\right)
\bigg|_{z_1= z_2 e^{z_0}}\bigg.
 .
$$
Therefore,   the product of 
\eqref{tren3} and \eqref{left3}  is equal to \eqref{right1}, so we conclude  that
the coefficients of all monomials \eqref{monomials} in \eqref{left} and in \eqref{right1} coincide, as required.
\end{prf}

\begin{lem}\label{lemlem4}
The map $Y_W(z)$ satisfies
   $\wht{\Sc}$-locality \eqref{localitymod}, i.e.
for any $u,v\in \Vc_{\hspace{-1pt}c}(\gl_N)$ and $k\in\mathbb{Z}_{\geqslant 0}$ there exists $p\in\mathbb{Z}_{\geqslant 0}$ such that   for all $w\in W$
\begin{align}
&(z_1-z_2)^{p}\ts Y_W(z_1)\big(1\otimes Y_W(z_2)\big)\iotaopjd\big(\wht{\mathcal{S}}(z_1/z_2)(u\otimes v)\otimes w\big)\non\\
&\qquad-(z_1-z_2)^{p}\ts Y_W(z_2)\big(1\otimes Y_W(z_1)\big)(v\otimes u\otimes w)
\in h^k W[[z_1^{\pm 1},z_2^{\pm 1}]].\label{1localitymod}
\end{align}
\end{lem}

\begin{prf}
Let $n,m,a_1,\ldots ,a_n,b_1,\ldots ,b_m ,k>0$ be arbitrary integers, $u=(u_1,\ldots ,u_n)$ and $v=(v_1,\ldots ,v_m)$ the families of variables. 
We will apply
$Y_W(z_1)\big(1\otimes Y_W(z_2)\big)\iotaopjd\wht{\mathcal{S}}(z_1/z_2)$, which corresponds to the first summand in \eqref{1localitymod}, to 
 \beq\label{elment}
T_{[n]}^{+13}(u)\ts T_{[m]}^{+24}(v)(\vac\ot\vac)
\eeq 
and then consider the coefficients of all monomials \eqref{monomials} in the resulting expression.
Note that the superscripts $1,2,3, 4$ in \eqref{elment} indicate the tensor factors as follows:
$$
  \smalloverbrace{(\ndo\CC^N)^{\ot n}}^{1} \ot 
	\smalloverbrace{(\ndo\CC^N)^{\ot m} }^{2} \ot 
	\smalloverbrace{\Vc_{\hspace{-1pt}c}(\gl_N)}^{3}\ot 
	\smalloverbrace{\Vc_{\hspace{-1pt}c}(\gl_N)}^{4} .
$$

Applying the map $\wht{\Sc}(z_1/z_2)$, given by  \eqref{sop1hat},
to \eqref{elment} and then
using \eqref{R2} 
we get 
\begin{align}
&
\wht{G}(z_1/z_2,u,v)\left(D_{[m]}^2 R_{nm}^{ + 12}(z_1 e^{u-v-h (N+ c)}/z_2)^{-1} (D_{[m]}^2)^{-1}\right)\label{sop5}
 \\
 & \cdotlr \left(  R_{nm}^{ + 12}(z_1 e^{u-v}/z_2) T_{[n]}^{+13}(u)  
R_{nm}^{  +12}(z_1   e^{u-v+h  c}/z_2) ^{-1}
 T_{[m]}^{+24}(v)   R_{nm}^{ + 12}(z_1 e^{u-v}/z_2)(\vac\otimes \vac)\right),\non
\end{align}
where the function $\wht{G}(x,u,v)$ is given by Lemma \ref{lemma4}. As we only consider the coefficients of  monomials \eqref{monomials}, it is sufficient to carry out the calculations modulo $U$, where
\begin{align*}
&U= \sum_{i=1}^n u_i^{a_i} V + \sum_{j=1}^m v_j^{b_j} V +h^k V\Fand \\
& V=(\ndo\CC^N)^{\ot (n+m)} \ot\Vc_{\hspace{-1pt}c}(\gl_N)^{\ot 2} ((z_1))((z_2))
[[u_1,\ldots ,u_n,v_1,\ldots ,v_m,h]] .
\end{align*}
By Lemma \ref{lemma4}, there exists an integer  $p\geqslant 0$ such that the image of the product of $(z_1-z_2)^p$ and \eqref{sop5} with respect to the map 
$\iotaop=\iotaall$
coincides with
\begin{align}
&
(z_1-z_2)^p\iotaop \wht{H}(z_2/z_1,u,v)\left(D_{[m]}^2 R_{nm}^{ + 12}(z_1 e^{u-v-h (N+ c)}/z_2)^{-1} (D_{[m]}^2)^{-1}\right)\label{ssp5}
 \\
 & \cdotlr \left(  R_{nm}^{+  12}(z_1 e^{u-v}/z_2) T_{[n]}^{+13}(u)  
R_{nm}^{  +12}(z_1   e^{u-v+h  c}/z_2) ^{-1}
 T_{[m]}^{+24}(v)   R_{nm}^{+  12}(z_1 e^{u-v}/z_2)(\vac\otimes \vac)\right)\non
\end{align}
modulo $U$, where the function $\wht{H}(x,u,v)$ is given by Lemma \ref{lemma4}. 
Note that  there are only finitely many monomials \eqref{monomials}. Therefore, due to Lemma \ref{poleslemma}, we can assume that $p=4p_0$ for some integer $p_0\geqslant 0$ such that all   coefficients of monomials \eqref{monomials} in 
\begin{alignat*}{2}
&\iotaop \left( (z_1-z_2)^{p_0}\ts  R_{nm}^{  12}(z_1 e^{u-v-h (N+ c)}/z_2)^{-1}\right) ,\quad
&&\iotaop\left( (z_1-z_2)^{p_0}\ts R_{nm}^{ 12}(z_1 e^{u-v}/z_2)\right) , \\
&\iotaop \left( (z_1-z_2)^{p_0}\ts R_{nm}^{  21}(z_2   e^{-u+v-h  c}/z_1)\right)  ,\quad
&&\iotaop\left(  (z_1-z_2)^{p_0}\ts R_{nm}^{ 21}(z_2 e^{-u+v}/z_1)^{-1} \right) 
\end{alignat*}
belong to
$
(\ndo\CC^N)^{\ot (n+m)} ((z_1,z_2))
$.
Due to the definition of 
the function 
$\wht{H}(x,u,v)$, see in particular \eqref{fzuv2}, we conclude by  \eqref{R2} and \eqref{arofexrr}  that the expression in \eqref{ssp5} equals
\begin{align}
&
 (z_1-z_2)^p \iotaop
\left(D_{[m]}^2 R_{nm}^{  12}(z_1 e^{u-v-h (N+ c)}/z_2)^{-1} (D_{[m]}^2)^{-1}\right)
\cdotlr \Big(  R_{nm}^{ 12}(z_1 e^{u-v}/z_2)\Big.\label{ssp6}
 \\
 &  \Big. \times T_{[n]}^{+13}(u)  
R_{nm}^{  21}(z_2   e^{-u+v-h  c}/z_1) 
 T_{[m]}^{+24}(v)   R_{nm}^{ 21}(z_2 e^{-u+v}/z_1)^{-1}(\vac\otimes \vac)\Big)\mod U.\non
\end{align}

Next,  we apply  $Y_W(z_1)(1\ot Y_W(z_2)) $  to \eqref{ssp6}, thus getting
\begin{align}
&
(z_1-z_2)^p 
\left(D_{[m]}^2 \iotaop\left( R_{nm}^{  12}(z_1 e^{u-v-h (N+ c)}/z_2)^{-1}\right) (D_{[m]}^2)^{-1}\right)
\cdotlr \Big(  \iotaop\left( R_{nm}^{ 12}(z_1 e^{u-v}/z_2)\right)\Big. \label{sssp7}
 \\
 &  \Big. \times 
\Lc_{[n]}^{13}(x)\big|_{x_i = z_1 e^{u_i}}  \big.
\iotaop\left( R_{nm}^{  21}(z_2   e^{-u+v-h  c}/z_1)\right) 
 \Lc_{[m]}^{23}(y)   \big|_{y_j = z_2 e^{v_j}}  \big.
\iotaop\left( R_{nm}^{ 21}(z_2 e^{-u+v}/z_1)^{-1}\right)\Big)  \non
\end{align}
modulo $U_0$, where $x=(x_1,\ldots ,x_n)$ and $y=(y_1,\ldots ,y_m)$ are the families of variables and
\begin{align*}
&U_0= \sum_{i=1}^n u_i^{a_i} V_0 + \sum_{j=1}^m v_j^{b_j} V_0 +h^k V_0\qquad\text{for} \\
& V_0=(\ndo\CC^N)^{\ot (n+m)}\ot  \om (W,W[[z_1^{\pm 1},z_2^{\pm 2}]])[[u_1,\ldots ,u_n,v_1,\ldots ,v_m]].
\end{align*}
By employing quantum current commutation relation \eqref{qcgen} we  rewrite \eqref{sssp7}
as
\begin{align}
&
\left(D_{[m]}^2 \iotaop\left((z_1-z_2)^{p_0}  R_{nm}^{  12}(z_1 e^{u-v-h (N+ c)}/z_2)^{-1}\right) (D_{[m]}^2)^{-1}\right)\non\\
&\cdotlr \Big(  \iotaop\left((z_1-z_2)^{p_0}  R_{nm}^{ 12}(z_1 e^{u-v}/z_2)\right) 
\iotaop\left( (z_1-z_2)^{p_0} R_{nm}^{ 12}(z_1 e^{u-v}/z_2)^{-1}\right)
\Big.  \label{ssp9}
 \\
 &   \Big. \times 
 \Lc_{[m]}^{23}(y)   \big|_{y_j = z_2 e^{v_j}}  \big.
\iotaop\left((z_1-z_2)^{p_0} R_{nm}^{ 12}(z_1 e^{u-v-hc}/z_2)\right)
\Lc_{[n]}^{13}(x)\big|_{x_i = z_1 e^{u_i}}  \big.
 \Big)\mod U_0. \non
\end{align}
Observe that all products in \eqref{ssp9} are well-defined modulo $U_0$ due to our choice of the integer $p=4p_0$ and Remark \ref{linremark}.
Canceling the $R$-matrices $   R_{nm}^{ 12}(z_1 e^{u-v}/z_2)^{\pm 1} $ and then using the following consequence of  the second crossing symmetry   relation in  \eqref{csym} which is verified by arguing as in Remark \ref{csrem},
$$
\left(D_{[m]}^2  R_{nm}^{  12}(z_1 e^{u-v-h (N+ c)}/z_2)^{-1} (D_{[m]}^2)^{-1}\right)
\cdotlr
  R_{nm}^{ 12}(z_1 e^{u-v-hc}/z_2),
$$
the expression in \eqref{ssp9} simplifies to
\begin{align}
 (z_1 -z_2)^p  \ts \Lc_{[m]}^{23}(y)   \big|_{y_j = z_2 e^{v_j}}  \big.
\Lc_{[n]}^{13}(x)\big|_{x_i = z_1 e^{u_i}}\mod U_0.\label{sopa}
\end{align}

Finally, consider the expression which corresponds to the second summand in \eqref{1localitymod}, i.e. which is obtained by applying 
$(z_1 -z_2)^p \ts Y_W(z_2)(1\ot Y_W(z_1))$
to $T_{[n]}^{+14}(u)\ts T_{[m]}^{+23}(v)(\vac\ot\vac)$. Clearly, all its coefficients with respect to monomials \eqref{monomials}   coincide with the corresponding coefficients in \eqref{sopa}, so the $\wht{\Sc}$-locality   follows.
\end{prf}

%%%%%%%%%%%%%%%%%%%%%%%%%%%%%%%%%%%
\subsection{Establishing the restricted \texorpdfstring{$\Uh$}{Uh(glNhat)}-module structure}\label{sec053}
%%%%%%%%%%%%%%%%%%%%%%%%%%%%%%%%%%%
Let $(W,Y_W)$ be a $\phi$-coordinated $\Vc_{\hspace{-1pt}c}(\gl_N)$-module for some $c\in\CC$, where $\phi(z_2,z_0) = z_2 e^{z_0}$. In this subsection, which consists of two lemmas, we finish the proof of the \hyperref[mainthm1]{Main Theorem} in the $\gl_N$ case.

\begin{lem}\label{lemlem5}
Formula \eqref{moduleformula} defines a unique structure of restricted  $\Uh$-module  of level $c$ on $W$.
\end{lem}

\begin{prf}
The uniqueness is clear as \eqref{moduleformula} determines the action of all generators of   $\Uh$ on $W$. 
We now use the Jacobi-type identity  given in Proposition \ref{Jacobi_prop} to check that \eqref{moduleformula} satisfies defining relation \eqref{qc} for the algebra $\Uh$ at the level $c$. 
Let $n\geqslant 0$ be an arbitrary integer. Choose  $p\geqslant 0$ such that the expressions
\begin{align}
&\iotaopdj (z_1 -z_2)^p    R_{12}(z_2 /z_1)^{-1}\ts T_{23}^+ (0)\ts R_{12}(z_2 e^{-hc}/z_1)\ts T_{14}^+(0) (\vac\ot\vac)\Fand\label{jac1}\\
&\iotaopjd  (z_1 -z_2)^p R_{12}(z_2 /z_1)^{-1}\ts T_{23}^+ (0)\ts R_{12}(z_2 e^{-hc}/z_1)\ts T_{14}^+(0) (\vac\ot\vac), \label{jac2}
\end{align}
whose coefficients belong to
$
(\ndo\CC^N )^{\ot 2} \ot \Vc_{\hspace{-1pt}c}(\gl_N)^{\ot 2} 
$,
coincide modulo $  h^n$. 
Note that the embedding map $\iotaopjd$ in \eqref{jac2} can be omitted as   both $R$-matrices are   Taylor series in $z_2 /z_1$, i.e. they consist of nonnegative powers of $z_2/z_1$. Furthermore,  we can assume that the integer $p$ is chosen so that expression \eqref{jac1} modulo $h^n$ is a polynomial in the variables $z_1^{\pm 1}, z_2^{\pm 1}$. Hence the embedding map $\iotaopdj$ can be also omitted when regarding \eqref{jac1} modulo $h^n$.
Applying   first term \eqref{Jacobi1} of the Jacobi identity on \eqref{jac2} we get
\begin{align*}
&(z_2 z)^{-1} \delta\left(\frac{z_1 -z_2}{z_2 z}\right) Y_W(z_1) (1 \ot Y_W(z_2))\\
&\times
(z_1 -z_2)^p    R_{12}(z_2 /z_1)^{-1}\ts T_{23}^+ (0)\ts R_{12}(z_2 e^{-hc}/z_1)\ts T_{14}^+(0) (\vac\ot\vac).
\end{align*}
Using \eqref{moduleformula} we rewrite this as
\begin{align*}
&(z_2 z)^{-1} \delta\left(\frac{z_1 -z_2}{z_2 z}\right) (z_1 -z_2)^p\ts  R_{12}(z_2 /z_1)^{-1}  \Lc_2(z_1)   R_{12}(z_2 e^{-hc}/z_1) \Lc_1(z_2)  .
\end{align*}
Due to the well-known $\delta$-function identity, 
\beq\label{delta567}
x\delta(x)=\delta(x),
\eeq 
by multiplying by $(z_2 z)^{-p}$ and then taking the residue $\rez_{z_2 z}$ we obtain
\beq\label{jac3}
 R_{12}(z_2 /z_1)^{-1} \ts\Lc_2(z_1)\ts R_{12}(z_2 e^{-hc}/z_1)\ts \Lc_1(z_2). 
\eeq

We now turn to  second term \eqref{Jacobi2} of the Jacobi identity. Choose $r\geqslant 0$ such that
\begin{alignat*}{2}
&A(z_1,z_2)\coloneqq\iotaopdj (z_1 - z_2)^r  R_{12}(z_2e^{hc}/z_1)^{-1}, \quad 
&&B(z_1,z_2)\coloneqq\iotaopdj (z_1 - z_2)^r \psi^2 R_{21}(z_1e^{-hc}/z_2),\\
&C(z_1,z_2)\coloneqq\iotaopdj (z_1 - z_2)^r  R_{12}(z_2/z_1), \quad 
&&D(z_1,z_2)\coloneqq\iotaopdj (z_1 - z_2)^r \psi^{-2} R_{21}(z_1/z_2)^{-1}
\end{alignat*}
belong to $(\ndo\CC^N)^{\ot 2} ((z_1,z_2))$ modulo $h^n$. 
As with \eqref{jac2}, observe that the embedding map $\iotaopdj$ can be omitted in the definitions of
$B(z_1,z_2)$ and $D(z_1,z_2)$ above.
By combining \eqref{rplusg} and unitarity property \eqref{uni}  we find
$$
A(z_1,z_2)\big|_{z_1=z_2 e^{z_0}}^{\modd h^n}\big.
=B(z_1,z_2)\big|_{z_1=z_2 e^{z_0}}^{\modd h^n}\big.
\fand
C(z_1,z_2)\big|_{z_1=z_2 e^{z_0}}^{\modd h^n}\big.
=D(z_1,z_2)\big|_{z_1=z_2 e^{z_0}}^{\modd h^n}\big. .
$$
Therefore, by the implication in \eqref{rem28} we conclude that
\beq\label{jac4}
A(z_1,z_2)=B(z_1,z_2)\mod h^n\Fand 
C(z_1,z_2)=D(z_1,z_2)\mod h^n.
\eeq
Consider   \eqref{jac1} modulo $h^n$.
Applying  second term \eqref{Jacobi2} of the Jacobi identity    we get\footnote{Note that, in contrast with \eqref{Jacobi1} and \eqref{Jacobi3}, the vectors $u$ and $v$ in \eqref{Jacobi2} are swapped, so that in \eqref{delta1234} we have $T_{24}^+ (0)$ and $T_{13}^+(0)$ instead of $T_{23}^+ (0)$ and $T_{14}^+(0)$.}
\begin{align}
&-(z_2 z)^{-1} \delta\left(\frac{z_2-z_1}{-z_2 z}\right) Y_W(z_2)(1\ot Y_W(z_1)) \iotaopdj\big( \wht{\Sc}(z_2 /z_1)\big. \non\\
&\times 
(z_1 -z_2)^p  R_{12}(z_2 /z_1)^{-1}\ts T_{24}^+ (0)\ts R_{12}(z_2 e^{-hc}/z_1)\ts T_{13}^+(0) (\vac\ot\vac)\big.\big)\mod h^n .\label{delta1234}
\end{align}
By using explicit formula \eqref{qva2hat} for the  map $\wht{\Sc}(x)$ we rewrite \eqref{delta1234} as
\begin{align*}
&-(z_2 z)^{-1} \delta\left(\frac{z_2-z_1}{-z_2 z}\right) Y_W(z_2)(1\ot Y_W(z_1))  \\
&\times 
(z_1 -z_2)^p\ts   T_{13}^+(0)\iotaopdj\left( R_{12}(z_2e^{hc} /z_1)^{-1}\right)  T_{24}^+ (0)\iotaopdj\left( R_{12}(z_2 /z_1) \right)   (\vac\ot\vac)\mod h^n .
\end{align*}
Next, the application of \eqref{moduleformula} gives us
\begin{align}
&-(z_2 z)^{-1} \delta\left(\frac{z_2-z_1}{-z_2 z}\right) 
(z_1 -z_2)^p  \Lc_1 (z_2)\non \\
&\times \iotaopdj\left( R_{12}(z_2e^{hc} /z_1)^{-1}\right) \Lc_2 (z_1)\iotaopdj\left( R_{12}(z_2 /z_1) \right) \mod h^n .\label{delta6778}
\end{align}
Note that \eqref{delta567} implies
$$
\delta\left(\frac{z_2-z_1}{-z_2 z}\right)=\frac{(z_2-z_1)^{r}}{(-z_2 z)^{r}}\delta\left(\frac{z_2-z_1}{-z_2 z}\right),
$$
so that
we can use both equalities in \eqref{jac4} to rewrite \eqref{delta6778} as
\begin{align*}
&-(z_2 z)^{-1} \delta\left(\frac{z_2-z_1}{-z_2 z}\right) 
(z_1 -z_2)^p  \Lc_1 (z_2)\ts R_{21}(z_1e^{-hc} /z_2)\ts \Lc_2 (z_1)\ts R_{21}(z_1 /z_2)^{-1}   \mod h^n,
\end{align*}
where the embedding maps $\iotaopdj$ are   omitted as   both $R$-matrices    consist of nonnegative powers of $z_1/z_2$.
Finally, multiplying by $(z_2 z)^{-p}$ and   taking the residue $\rez_{z_2 z}$ we get
\beq\label{jac5}
-  \Lc_1 (z_2)\ts R_{21}(z_1e^{-hc} /z_2)\ts \Lc_2 (z_1)\ts R_{21}(z_1 /z_2)^{-1}   \mod h^n. 
\eeq

By applying    third term \eqref{Jacobi3} of the Jacobi identity to  \eqref{jac1}  we get
\begin{align}
&z_1^{-1}\delta\left(\frac{z_2 (1+z)}{z_1}\right) Y_W (z_2)\left(Y(\log(1+z))\ot 1\right)\non\\
&\times\iotaopdj  (z_1 -z_2)^p  R_{12}(z_2 /z_1)^{-1}\ts T_{23}^+ (0)\ts R_{12}(z_2 e^{-hc}/z_1)\ts T_{14}^+(0) (\vac\ot\vac)\mod h^n\label{jac7}.
\end{align}
As before, by  \eqref{rplusg} and  unitarity property \eqref{uni}    there exists $s\geqslant 0$ such that
\beq\label{temp_uni8}
\iotaopdj (z_1-z_2)^s  R_{12}(z_2 e^{-hc}/z_1)
= 
\iotaopdj (z_1-z_2)^s \psi^{-2} R_{21}(z_1 e^{hc}/z_2)^{-1} \mod h^n.
\eeq
Using the $\delta$-function identities 
\beq\label{theiid}
\left(\frac{z_1}{z_2}\right)^l \delta\left(\frac{z_2(1+z)}{z_1}\right)=(1+z)^l\,\delta\left(\frac{z_2(1+z)}{z_1}\right),
\eeq
which follow directly from \eqref{delta567},
one  can easily derive
$$
\delta\left(\frac{z_2(1+z)}{z_1}\right)=\frac{(z_1-z_2)^k}{(z_2 z)^k}\delta\left(\frac{z_2(1+z)}{z_1}\right),
$$
in particular for $k=p,s$.
Therefore, we can employ \eqref{temp_uni8} to rewrite \eqref{jac7} as
\begin{align*}
& \psi^{-2} z_1^{-1}\delta\left(\frac{z_2 (1+z)}{z_1}\right) ( z_2 z)^{p} \ts Y_W (z_2)\left(Y(\log(1+z))\ot 1\right) \\
&\times  \iotaopdj\left( R_{12}(z_2 /z_1)^{-1}\right) T_{23}^+ (0)\iotaopdj\left(  R_{21}(z_1 e^{hc}/z_2)^{-1}\right) T_{14}^+(0) (\vac\ot\vac)\mod h^n.
\end{align*}
Next, using definition \eqref{qva1} of the vertex operator map and \eqref{theiid} we get
\begin{align*}
&\psi^{-2} z_1^{-1}\delta\left(\frac{z_2 (1+z)}{z_1}\right) ( z_2 z)^{p} \ts Y_W (z_2) \iotaopdj\left( R_{12}(z_2 /z_1)^{-1}\right)\\
&\times 
T_{23}^+ (\log(1+z))\ts
T_{23}^* (\log(1+z) +hc/2)^{-1}
\ts R_{21}((1+z) e^{hc})^{-1} \ts T_{13}^+(0) \vac \mod h^n.
\end{align*}
Finally, we use  relation \eqref{rtt3} to swap the operators $T^*_{23}$ and $T^+_{13}$, and then we employ the identity $T^*(x)\vac =\vac$, thus getting
\begin{align}
& \psi^{-2} z_1^{-1}\delta\left(\frac{z_2 (1+z)}{z_1}\right) ( z_2 z)^{p} \ts Y_W (z_2)   \iotaopdj\left(  R_{12}(z_2 /z_1)^{-1}\right) \non\\
&\times T_{23}^+ (\log(1+z))\ts
T_{13}^+(0)
\iotaopdj\left( R_{21}(z_1 /z_2)^{-1}\right)   \vac \mod h^n. \label{thisistheend}
\end{align}
It is clear that the application of the module map $Y_W (z_2)$ in \eqref{thisistheend} will not produce any negative powers of the variable $z$. Therefore, 
multiplying \eqref{thisistheend} by $(z_2 z)^{-p}$ and then taking the residue $\rez_{z_2 z}$ we obtain $0\,\modd h^n$. 
Hence combining the Jacobi-type identity from Proposition \ref{Jacobi_prop}
with  \eqref{jac3} and \eqref{jac5} we obtain the equality
\begin{align*}
&  \Lc_1 (z_2)\ts R_{21}(z_1e^{-hc} /z_2)\ts \Lc_2 (z_1)\ts R_{21}(z_1 /z_2)^{-1}
\\
  & - R_{12}(z_2 /z_1)^{-1} \ts\Lc_2(z_1)\ts R_{12}(z_2 e^{-hc}/z_1)\ts \Lc_1(z_2) =0 \mod h^n    
\end{align*}
for operators on $W$.
As the integer  $n$ was arbitrary, we conclude that
the given equality holds for all $n$. Hence we proved that
  \eqref{moduleformula} satisfies  quantum current commutation relation \eqref{qc}, so that it defines the structure of  $\Uh$-module of level $c$  on $W$, as required. In the end, in order to finish the proof, it remains to observe that $W$ is a topologically free $\CC[[h]]$-module and, furthermore, restricted $\Uh$-module by Definition \ref{phimod}.
\end{prf}

The next lemma completes the proof of the \hyperref[mainthm1]{Main Theorem} for $\g_N =\gl_N$.
\begin{lem}\label{lemlem6}
A topologically free  $\CC[[h]]$-submodule $W_1$ of $W$ is a $\phi$-coordinated $\Vc_{\hspace{-1pt}c} (\gl_N)$-submodule of $W$ if and only if it is an $\Uh$-submodule of $W$.
\end{lem}

\begin{prf}
Suppose that $W_1$ is a $\phi$-coordinated $\Vc_{\hspace{-1pt}c}(\gl_N)$-submodule of $W$. Then  
$$\Lc (z)w=Y_W(T^+(0)\vac,z)w \,\in\, \ndo\CC^N \ot W_1 ((z))[[h]]\quad\text{for any }w\in W_1,$$
so $W_1$ is clearly  an $\Uh$-submodule of $W$. 

Conversely, suppose that $W_1$ is a topologically free $\Uh$-submodule of $W$.
Clearly, $W_1$ is a restricted $\Uh$-module of level $c$, so by Proposition \ref{restricted496}  we have
$$
\Lc_{[n]}(x_1,\ldots ,x_n)w\,\in\,\left(\ndo\CC^N\right)^{\ot n} \ot W_1 ((x_1,\ldots ,x_n))[[h]]\quad\text{for all }n\geqslant 1\text{ and } w\in W_1.
$$
Applying the substitutions $x_i=ze^{u_i}$ with $i=1,\ldots ,n$ we get
$$
 \Lc_{[n]}(x)\big|_{x_i = ze^{u_i}}\big. w =
Y_W(T_{[n]}^+(u)\vac,z) w
\,\in\,\left(\ndo\CC^N\right)^{\ot n} \ot W_1 ((z))[[u_1,\ldots ,u_n,h]].
$$
By \cite[Sect. 3.4]{EK3}, see also \cite[Prop. 2.4]{KM}, the coefficients of all matrix entries of $T_{[n]}^+(u)$, $n\geqslant 1$, and $\vac$ span an $h$-adically dense $\CC[[h]]$-submodule of $\Vc_{\hspace{-1pt}c}(\gl_N)$, so we conclude that $W_1$ is a $\phi$-coordinated $\Vc_{\hspace{-1pt}c}(\gl_N)$-submodule of $W$, as required.
\end{prf}

%%%%%%%%%%%%%%%%%%%%%%%%%%%%%%%%%%%
\subsection{Proof of the \texorpdfstring{\hyperref[mainthm1]{Main Theorem}}{Main Theorem} in the \texorpdfstring{$\sll_N$}{slN} case}\label{sec054}
%%%%%%%%%%%%%%%%%%%%%%%%%%%%%%%%%%%

For any integer $n=1,\ldots ,N$ set
$$
 u_{[n]}=(u,u-h,\ldots ,u-(n-1)h)\fand x_{[n]}=(x,xe^{-h},\ldots ,xe^{-(n-1)h})  .
$$
Let
$
P^{(n)}\colon x_1\ot\ldots \ot x_n \mapsto x_n\ot\ldots \ot x_1
$ 
be the permutation operator
 on $(\CC^N)^{\ot n}$.
Write
\begin{align*}
&\Lc_{[n]}(x_{[n]})=\Lc_{[n]}(x_1,\ldots ,x_n)\big|_{x_1 = x,\ldots, x_n = xe^{-(n-1)h}}\big. , \\
&\cev{\Lc}_{[n]} (x_{[n]}) =P^{(n)}\,  \Lc_{[n]}(x_n,\ldots ,x_1)\big|_{x_1 = x,\ldots, x_n = xe^{-(n-1)h}}\big. \, P^{(n)} . 
\end{align*}
We first list some useful properties of the anti-symmetrizer $A^{(n)}$ defined by \eqref{anti}.

\begin{lem}  \label{profuzion}
 For any $n=1,\ldots ,N$ we have
\begin{gather}
A^{(n)}\ts \Lc_{[n]}(x_{[n]}) = \cev{\Lc}_{[n]}(x_{[n]}) \ts A^{(n)} ,\label{497c}\\
A^{(n)}\ts D_1\ldots D_n =  D_1\ldots D_n \ts A^{(n)} \label{497a},\\
A^{(N)}\ts \wvr{R}_{1N}^{12}(y/x_{[N]})=
\cev{\wvr{R}}_{1N}^{12}(y/x_{[N]})\ts A^{(N)}=
A^{(N)} e^{-(N-1)h/2} \ts\frac{x-e^{(N-1)h}y}{x-y},\label{497d}
\end{gather}
where the arrow in $\cev{\wvr{R}}_{1N}^{12}(y/x_{[N]})$ indicates  the reversed order of factors.
The coefficients in \eqref{497d} belong to $\ndo\CC^N \ot(\ndo\CC^N)^{\ot N}[[h]]$ and the anti-symmetrizer $A^{(N)}$ is applied on the tensor factors $2,\ldots ,N+1$.
\end{lem}

\begin{prf}
Equality \eqref{497c} is verified by using Yang--Baxter equation \eqref{YBE}, generalized quantum current commutation relation \eqref{qcgen}    and the following   case of the  fusion procedure
for the two-parameter $R$-matrix $\R(x,y)=(xe^{-h/2}-ye^{h/2})\R(x/y)$ going back to \cite{C},
\beq\label{fuzion}
\prod_{i=1,\ldots ,n-1}^{\longrightarrow}\prod_{j=i+1,\ldots ,n}^{\longrightarrow}\R_{ij}(xe^{-(i-1)h},xe^{-(j-1)h})
=
n!\, x^{\frac{n(n-1)}{2}}
\prod_{0\leqslant i< j\leqslant n-1} (e^{-ih}-e^{-jh})\ts A^{(n)}.
\eeq
Equality \eqref{497a} follows from the identities
$$
D_1 D_2 =D_2 D_1\fand \R(x,y) D_1 D_2 =D_2 D_1 \R(x,y)
$$
 while \eqref{497d} is established in the proof of \cite[Lemma 4.3]{FJMR}.
\end{prf}

Let $W$ be a  restricted $\Uhsl$-module of level $c\in \CC$. 

\begin{lem}\label{lemlem11}
Formula \eqref{formula}, together with $Y_W(\vac,z)=1_W$, defines a unique
 structure of $\phi$-coordinated $\Vc_{\hspace{-1pt}c}(\sll_N)$-module on $W$, where $\phi(z_2,z_0) = z_2 e^{z_0}$.
\end{lem}

\begin{prf}
In order to prove the lemma, it is sufficient to verify that \eqref{formula}, together with $Y_W(\vac,z)=1_W$, defines a   $\CC[[h]]$-module map $\Vc_{\hspace{-1pt}c}(\sll_N)\ot W \to W((z))[[h]]$. Indeed, all other properties of the aforementioned map are recovered by arguing as in the $\g_N =\gl_N$ case; see Subsection \ref{sec052}.
Therefore, 
we have to show that the map $v\mapsto Y_W (v,z)$ preserves the ideal of relations \eqref{rtt} and \eqref{qdetvoa2}. However, it is sufficient to consider \eqref{qdetvoa2} as relations \eqref{rtt} are already taken care of in the proof of Lemma \ref{lemlem1}. 

Let $n$ and $m$ be nonnegative integers.  Introduce the families of variables $v=(v_1,\ldots ,v_n)$ and $w=(w_1,\ldots ,w_m)$. Consider the image of the expression
\beq\label{231}
T_{[n]}^{+13}(v) \ts\qdet T^+(u)\ts T_{[m]}^{+23}(w)\vac \in (\ndo\CC^N)^{\ot (n+m)} \ot \Vc_{\hspace{-1pt}c}(\sll_N)[[v_1,\ldots, v_n,u,w_1,\ldots w_m]]
\eeq
with respect to $Y_W(z)$. Introduce the tensor product
\beq\label{232}
\smalloverbrace{(\ndo\CC^N)^{\ot n}}^{1} \ot 
\smalloverbrace{(\ndo\CC^N)^{\ot N}}^{2} \ot
\smalloverbrace{(\ndo\CC^N)^{\ot m}}^{3} \ot
\smalloverbrace{\Vc_{\hspace{-1pt}c}(\sll_N)}^{4}  
\eeq
and write
$$
A^{(N)}_2=1^{\ot n}\ot A^{(N)} \ot 1^{\ot m}\fand
D_{[N]}^2=1^{\ot n}\ot D_{[N]} \ot 1^{\ot m}=D_{n+1}\ldots D_{n+N}.
$$
By the definition of quantum determinant given by \eqref{qdetvoa}, using the labels in \eqref{232} to indicate the corresponding tensor factors,  \eqref{231} can be expressed as
\beq\label{233}
\tr_{n+1,\ldots , n+N}\,A^{(N)}_2\ts T_{[n]}^{+14}(v) \ts  T_{[N]}^{+24}(u_{[N]})\ts T_{[m]}^{+34}(w)\ts D_{[N]}^2 \vac .
\eeq
By \eqref{formula}, the image of \eqref{233} with respect to $Y_W(z)$ equals
\beq\label{234}
\tr_{n+1,\ldots , n+N}\,A^{(N)}_2\ts \Lc_{[n+N+m]}(\wvr{x},x_{[N]},\wvr{y})\Big|_{\substack{x_1 = ze^{v_1},\,\ldots, \, x_n = ze^{v_n},\, x=ze^u\\  y_1 = ze^{w_1},\, \ldots,\, y_n = ze^{w_m}}}\Big.
\ts D_{[N]}^2,
\eeq
 where 
$
\wvr{x}=(x_1,\ldots ,x_n)$  and $\wvr{y}=(y_1,\ldots ,y_m).
$
Using    generalized quantum current commutation relation \eqref{qcgen} we transform $A^{(N)}_2\Lc_{[n+N+m]}(\wvr{x},x_{[N]},\wvr{y})$ and bring it to the form
\begin{align}
&A^{(N)}_2\ts \Lc_{[n]}^{14}(\wvr{x})\ts R_{nN}^{21}(x_{[N]}e^{-hc}/\wvr{x})\ts
R_{nm}^{31}(\wvr{y}e^{-hc}/\wvr{x})\ts \Lc_{[N]}^{24}(x_{[N]})\non
 \\
&\times 
R_{Nm}^{32}(\wvr{y}e^{-hc}/x_{[N]})\ts
\Lc_{[m]}^{34}(\wvr{y})\ts
R_{Nm}^{32}(\wvr{y}/x_{[N]})^{-1}\ts R_{nm}^{31}(\wvr{y} /\wvr{x})^{-1}\ts
R_{nN}^{21}(x_{[N]} /\wvr{x})^{-1}\label{235}.
\end{align}
By employing \eqref{fqhqf} and \eqref{497d} one can verify the following identities: 
$$
\begin{aligned}
&A^{(N)}_2\ts R_{nN}^{21}(x_{[N]}e^{-hc}/\wvr{x}) = e^{-n(N-1)h/2} A^{(N)}_2,\quad&
&A^{(N)}_2\ts R_{nN}^{21}(x_{[N]} /\wvr{x})^{-1} = e^{n(N-1)h/2} A^{(N)}_2,\\
&A^{(N)}_2\ts R_{Nm}^{32}(\wvr{y}e^{-hc}/x_{[N]}) = e^{-m(N-1)h/2} A^{(N)}_2,\quad&
&A^{(N)}_2\ts R_{Nm}^{32}(\wvr{y}/x_{[N]})^{-1} = e^{m(N-1)h/2} A^{(N)}_2.
\end{aligned}
 $$
As the anti-symmetrizer $A^{(N)}_2$   commutes with the terms $R_{nm}^{31}(\wvr{y}e^{-hc}/\wvr{x})$, $R_{nm}^{31}(\wvr{y} /\wvr{x})^{-1}$, $\Lc_{[n]}^{14}(\wvr{x})$ and $\Lc_{[m]}^{34}(\wvr{y})$, by combining the above identities and
 \eqref{497c}, we rewrite \eqref{235} as
\begin{align}\label{237}
 \Lc_{[n]}^{14}(\wvr{x})\ts  
R_{nm}^{31}(\wvr{y}e^{-hc}/\wvr{x})\ts A^{(N)}_2\ts \Lc_{[N]}^{24}(x_{[N]})\ts
\Lc_{[m]}^{34}(\wvr{y})\ts
R_{nm}^{31}(\wvr{y} /\wvr{x})^{-1}
 .
\end{align}

Note that the expression in \eqref{234} is obtained from \eqref{237} by applying the substitutions
$$
x_1 = ze^{v_1},\,\ldots,\, x_n = ze^{v_n},\, x=ze^u,\,y_1 = ze^{w_1},\,\ldots,\, y_n = ze^{w_m},
$$
then multiplying by $D_{[N]}^2$ from the right and, finally, taking the trace $\tr_{n+1,\ldots , n+N}$. However, as $D_{[N]}^2$ commutes with the terms    $\Lc_{[m]}^{34}(\wvr{y})$ and $R_{nm}^{31}(\wvr{y} /\wvr{x})^{-1}$, it is clear that applying the aforementioned transformations to  \eqref{237} and using   definition of quantum determinant \eqref{qdet497} results in
\begin{align}
 &\Lc_{[n]}^{13}(\wvr{x})\ts  
R_{nm}^{21}(\wvr{y}e^{-hc}/\wvr{x})\ts\qdet\Lc (x)\ts
\Lc_{[m]}^{23}(\wvr{y})\ts
R_{nm}^{21}(\wvr{y} /\wvr{x})^{-1}\Big|_{\substack{x_1 = ze^{v_1},\,\ldots,\, x_n = ze^{v_n},\,x=ze^u  \\ y_1 = ze^{w_1},\,\ldots,\, y_n = ze^{w_m}}}\Big. ,\label{238}
\end{align}
where, due to   application of the trace, the tensor factors in \eqref{238} are now labeled    in accordance with \eqref{231}. As $\qdet\Lc (x)=1$ in $\Uhsl$ we conclude by quantum current commutation relation \eqref{qcgen} that \eqref{238} is equal to
\begin{align*}
\Lc_{[n+m]} (\wvr{x},\wvr{y}) \Big|_{\substack{x_1 = ze^{v_1},\,\ldots,\, x_n = ze^{v_n}    \\ y_1 = ze^{w_1},\,\ldots,\, y_n = ze^{w_m}}}\Big. 
=Y_W(T_{[n]}^{+13}(v)  \ts T_{[m]}^{+23}(w)\vac,z)
. 
\end{align*}
Therefore, the images of \eqref{231} and $T_{[n]}^{+13}(v)  T_{[m]}^{+23}(w)\vac$  with respect to $Y_W(z)$ coincide, so we conclude that the $\CC[[h]]$-module map $\Vc_{\hspace{-1pt}c}(\sll_N)\ot W \to W((z))[[h]]$ is well-defined by \eqref{formula}, as required.
\end{prf}

Let $(W,Y_W)$ be a $\phi$-coordinated $\Vc_{\hspace{-1pt}c}(\sll_N)$-module for some $c\in\CC$, where $\phi(z_2,z_0) = z_2 e^{z_0}$. 
In order to prove that \eqref{moduleformula} defines a unique structure of restricted  $\Uhsl$-module  of level $c$ on $W$, we need the following identity.

\begin{lem}\label{lemlemC}
For any positive integer $n$ the   identity  
\begin{align}
Y_W\left((T_1^+((n-1)h)T_2^+((n-2)h)\ldots T_n^+(0)\vac, ze^{-(n-1)h}\right)&\non\\
\qquad=\Lc_{[n]}(x_1,x_2,\ldots ,x_n)\big|_{x_1 = z , x_2=ze^{-h},\ldots, x_n = ze^{-(n-1)h}}\big.  &
\label{id460}
\end{align}
holds for operators on $W$,
where the action of $\Lc_{[n]}(x_1,\ldots ,x_n)$ on $W$ is given by formula \eqref{Ln} with $\Lc(x)=Y_W(T^+(0)\vac,x)$.
\end{lem}

\begin{prf}
We derive \eqref{id460} using  the weak associativity property.   Let   $k$ be a positive integer. By  \eqref{associativitymod0} and \eqref{associativitymod} there  exists an integer $p\geqslant 0$ such that
$$
(z_1 -z_2)^p\ts Y_W(T_1^+(0)\vac,z_1)\ts Y_W(T_2^+(0)\vac,z_2)
=(z_1 -z_2)^p\ts \Lc_1(z_1)\ts \Lc_2(z_2)
$$
belongs to $(\ndo\CC^N)^{\ot 2}\ot\om (W,W((z_1,z_2)))$ modulo $h^k$ and such that  
\beq\label{lab1}
\left((z_1 -z_2)^p\ts \Lc_1(z_1)\ts \Lc_2(z_2)\right) \big|_{z_1=z_2 e^{z_0}}^{\modd h^k}  \big.
\fand
z_2^p (e^{z_0}-1)^p\ts Y_W(Y(T_1^+(0)\vac,z_0)T_2^+(0)\vac,z_2)
\eeq
coincide modulo $h^k$. Using relation \eqref{rtt3} and then the first crossing symmetry property in \eqref{csym_equiv}    we express the second term in \eqref{lab1}   as
\beq\label{ned4}
z_2^p (e^{z_0}-1)^p\ts (D_1 R(e^{z_0+h(c+N)})D_1^{-1})\cdotrl\left(Y_W(T_1^+(z_0) T_2^+(0)\vac,z_2)\ts R(e^{z_0})^{-1}\right).
\eeq
The first crossing symmetry property in \eqref{csym_equiv} and unitarity   \eqref{uni} imply the identities
$$
R_{21}(e^{-z_0-hc})\cdotrl (D_1 R(e^{z_0+h(c+N)})D_1^{-1}) =1
\fand
R(e^{z_0})^{-1}R_{21}(e^{-z_0})^{-1} =1,
$$
which
enable us to move the $R$-matrices appearing in \eqref{ned4} from the second term in \eqref{lab1}  to the first term in \eqref{lab1}.
Hence we find that
\begin{gather}
\left(R_{21}(e^{-z_0-hc})\cdotrl
\left((z_1 -z_2)^p\ts \Lc_1(z_1)\ts \Lc_2(z_2)\right)\big|_{z_1=z_2 e^{z_0}}^{\modd h^k}  \big.\right) R_{21}(e^{-z_0})^{-1}\label{lab2}
\\
\text{and}\quad z_2^p (e^{z_0}-1)^p\ts   Y_W(T_1^+(z_0) T_2^+(0)\vac,z_2) \label{lab3} 
\end{gather}
coincide modulo $h^k$.   Without loss of generality we can assume that the integer $p$ is sufficiently large, so that we       conclude by Lemma \ref{poleslemma}  that \eqref{lab2} is equal to
\begin{align}
&\left((z_1 -z_2)^p\ts \Lc_1(z_1)\ts R_{21}(z_2 e^{-hc}/z_1)\ts \Lc_2(z_2) \ts R_{21}(z_2/z_1)^{-1}\right)\big|_{z_1=z_2 e^{z_0}}^{\modd h^k}  \big. .\label{mpim}
\end{align}
By employing \eqref{Ln} for $n=2$    and   
the relation
$$\Lc_{[2]}(z_1,z_2)\in(\ndo \CC^N)^{\ot 2}\ot \om (W,W((z_1,z_2))[[h]]),$$ 
which is
verified by arguing as in the proof of Proposition \ref{restricted496},
we rewrite \eqref{mpim} as
\begin{align} 
\left((z_1 -z_2)^p\ts \Lc_{[2]}(z_1,z_2) \right)\big|_{z_1=z_2 e^{z_0}}^{\modd h^k}  \big.  
=z_2^p (e^{z_0}-1)^p
\left(  \Lc_{[2]}(z_1,z_2) \right)\big|_{z_1=z_2 e^{z_0}}^{\modd h^k}  \big.  .\label{lab4} 
\end{align}
Thus we proved that  \eqref{lab3} and \eqref{lab4} coincide modulo $h^k$.
Hence, multiplying \eqref{lab3} and \eqref{lab4} by $z_2^{-p} (e^{z_0}-1)^{-p}$ we find that
$$
   \Lc_{[2]}(z_1,z_2)  \big|_{z_1=z_2 e^{z_0}}  \fand  Y_W(T_1^+(z_0) T_2^+(0)\vac,z_2)
$$
coincide modulo $h^k$. Moreover, by setting $z_0 =h$ and $z_2=ze^{-h}$ we conclude that
$$
   \Lc_{[2]}(z_1,z_2)  \big|_{z_1=z,\, z_2=ze^{-h}}   \fand  Y_W(T_1^+(h) T_2^+(0)\vac,ze^{-h})
$$
coincide modulo $h^k$. As the integer $k>0$ was arbitrary, this implies   equality \eqref{id460} for $n=2$. The general case is proved by induction on $n$.
\end{prf}

The next two lemmas complete the proof of the \hyperref[mainthm1]{Main Theorem} for $\g_N =\sll_N$. The second lemma follows by the same arguments as for   Lemma \ref{lemlem6},  so  we omit its proof.

\begin{lem}\label{lemlemD}
Formula \eqref{moduleformula} defines a unique structure of restricted  $\Uhsl$-module  of level $c$ on $W$.
\end{lem}

\begin{prf}
Due to the proof of Lemma \ref{lemlem5}, it is sufficient to verify the equality $\qdet\Lc(z)=1$   on $W$, where the action  of $\Lc(z)$ on $W$ is given by \eqref{moduleformula}. By \eqref{qdet497} and \eqref{id460}, the action of  quantum determinant of $\Lc(z)$ on $W$ is given by
\begin{align*}
&\,\tr_{1,\ldots, N}\, A^{(N)}\ts
Y_W\left((T_1^+((N-1)h)T_2^+((N-2)h)\ldots T_N^+(0)\vac, ze^{-(N-1)h}\right)\ts D_1\ldots D_N\\
=&\, Y_W\left(\tr_{1,\ldots ,N}\, A^{(N)}(T_1^+((N-1)h)T_2^+((N-2)h)\ldots T_N^+(0)\vac \ts D_1\ldots D_N, ze^{-(N-1)h}\right).
\end{align*}
By applying \eqref{qdetvoa} with $u=(N-1)h$ and \eqref{qdetvoa2} the given expression takes the form
$$\ts
Y_W\left(\qdet T^+ ((N-1)h)\vac , ze^{-(N-1)h}\right)=Y_W\left(\vac , ze^{-(N-1)h}\right).
$$
Finally,  Definition \ref{phimod} implies $Y_W\left(\vac , ze^{-(N-1)h}\right)=1$, 
which completes the proof.
\end{prf}

\begin{lem}\label{lemlemE}
A topologically free  $\CC[[h]]$-submodule $W_1$ of $W$ is a $\phi$-coordinated $\Vc_{\hspace{-1pt}c} (\sll_N)$-submodule of $W$ if and only if $W_1$ is an $\Uhsl$-submodule of $W$.
\end{lem}

%%%%%%%%%%%%%%%%%%%%%%%%%%%%%%%%%%%
%%%%%%%%%%%%%%%%%%%%%%%%%%%%%%%%%%%
\section{Image of the center of the quantum affine vertex algebra}\label{newsec02}
%%%%%%%%%%%%%%%%%%%%%%%%%%%%%%%%%%%
%%%%%%%%%%%%%%%%%%%%%%%%%%%%%%%%%%%

In this  section, we briefly discuss a connection between families of central elements for  the quantum affine vertex algebra
and the quantum affine algebra 
established by the $\phi$-coordinated module map from  the \hyperref[mainthm1]{Main Theorem}.

%%%%%%%%%%%%%%%%%%%%%%%%%%%%%%%%%%%
\subsection{Noncritical level}\label{newsec022}
%%%%%%%%%%%%%%%%%%%%%%%%%%%%%%%%%%%

Following \cite{JKMY}, we define the {\em center} of the  quantum  vertex algebra  $\Vccgll $ at the level $c\in \CC$
 as the $\CC[[h]]$-submodule  
$$
\z(\Vccgll) =\left\{
v\in \Vccgll\,:\, Y(w,z)v\in \Vccgll[[z]]\text{ for all }w\in\Vccgll
\right\}.
$$
For more details on the notion of center of quantum vertex algebra see \cite[Thm. 1.4]{DGK} and \cite[Sect. 3.2]{JKMY}. 
Observe that \eqref{formula} implies the identity
\beq\label{noncrit}
Y_W(\qdet T^+(0)\vac, z)=\qdet \Lc (z)
\eeq
on any restricted $\Uh$-module $W$ of level $c\in\CC$. 
By \cite[Prop. 3.10]{KM} the coefficients of the quantum determinant $\qdet T^+(u)$, as given by \eqref{detkoef}, belong to the center of the quantum   vertex algebra $\Vccgll$ for any $c\in\CC$.
The next proposition, which is well-known, provides a quantum affine algebra counterpart of this  fact; cf. \cite{FJMR}. 
We formulate the proposition and  outline its proof in terms of Ding's quantum current realization for completeness.
\begin{pro}\label{qdetpro}
 For any $c\in\CC$
all coefficients $d_r$ of the quantum determinant $\qdet \Lc (z)$, as given by \eqref{qdetc},  belong to the center of the  quantum affine   algebra   $\Uh_c$.
\end{pro}

\begin{prf}
It is sufficient to prove the equality
\beq\label{781}
\Lc (y)\, \qdet\Lc(x)=\qdet\Lc(x)\, \Lc (y)
\eeq
in $\ndo\CC^N\ot \Uh_c$. 
By \eqref{qdet497} the left hand side in \eqref{781} equals
\beq\label{782}
 \tr_{1,\ldots ,N} \, \Lc_0 (y)\ts A^{(N)}\ts\Lc_{[N]}(x_{[N]}) \ts D_{[N]},\quad\text{where }D_{[N]}= D_1\ldots D_N
\eeq
and
 the coefficients of the expression under the trace belong  to the tensor product $\ndo\CC^N \ot(\ndo\CC^N)^{\ot N}\ot \Uh_c$. The   copies of $\ndo\CC^N$ in \eqref{782} are labeled by $0,\ldots, N$. The matrix $\Lc(y)$ is applied on the tensor factor $0$ while the remaining terms,   $A^{(N)}$, $\Lc_{[N]}(x_{[N]})$ and $D_{[N]}$ are applied on the tensor factors $1,\ldots ,N$. 
By   $\Lc_0 (y)  A^{(N)}=  A^{(N)}\Lc_0 (y)$ and  generalized   quantum current commutation relation \eqref{qcgen} we rewrite \eqref{782} as
\beq\label{498}
 \tr_{1,\ldots ,N} \,  A^{(N)}\left( A\cdotrl\left(\left(B\ts\Lc_{[N]}(x_{[N]})\ts C\ts  \Lc_0 (y)\right)E\right)\right)
  D_{[N]},  
	\eeq
	where
	$$
\begin{aligned}
&A=D_{[N]}^{-1}\ts R_{1N}^{21}(x_{[N]}e^{-(N+c)h}/y)^{-1}D_{[N]},\quad&
&B=R_{1N}^{12}(y/x_{[N]})^{-1},\\
&C=R_{1N}^{12}(ye^{-hc}/x_{[N]}),&
&E=R_{1N}^{21}(x_{[N]}/y).
\end{aligned}
$$
Note that the  element $A$ is found via the second crossing symmetry property in \eqref{csym}; see also Remark \ref{csrem}.
Next, by using \eqref{fqhqf}  and \eqref{497d} one can verify the following equalities:
\beq\label{497b}
A^{(N)}\ts Z=\lambda_Z\ts A^{(N)} \text{ for }  Z=A,B,C,E\text{ and } \lambda_A=\lambda_B=\lambda_C^{-1}=\lambda_E^{-1}= e^{(N-1)h/2}.
\eeq
Using \eqref{497c} and \eqref{497b}    we move the anti-symmetrizer   in \eqref{498} to the right, thus getting
$$
 \tr_{1,\ldots ,N} \,      \cev{\Lc}_{[N]}(x_{[N]}) \ts \Lc_0 (y)\ts    A^{(N)} 
  D_{[N]} = 
	\tr_{1,\ldots ,N} \,      \cev{ \Lc}_{[N]}(x_{[N]}) \ts    A^{(N)} 
  D_{[N]}   \ts \Lc_0 (y) .
$$
Finally, we use  \eqref{497c} to move the anti-symmetrizer  $A^{(N)}$ to the left, thus getting the right hand side in \eqref{781}, as required.
\end{prf}

Following \cite[Sect. 3.3]{loop}, we define the {\em submodule of invariants} of the   vacuum module  $\Vccgl$ as the $\CC[[h]]$-submodule  
$$
\z (\Vccgl)=\left\{v\in \Vccgl\,:\, \Lc(z)v\in\Vccgl[[z]]    \right\}.
$$
Recall Corollary \ref{maincor}.
By setting  $W=\Vccgl$ in \eqref{noncrit} and then applying the resulting equality on $1\in\Vccgl$ one recovers the  invariants of the vacuum module; cf. \cite{FJMR}.
\begin{kor}
For any $c\in\CC$ all coefficients of the series
$$
\wvr{\ell}_N (z)\coloneqq Y_{\Vccgl} (\qdet T^+(0), z) 1 = \qdet \Lc (z)1 \in  \Vccgl  [[z]].
$$
belong to the submodule of invariants $\z (\Vccgl )$.
\end{kor}

\begin{prf}
	The Corollary follows by applying identity \eqref{781} on $1\in\Vccgl$.
\end{prf}

%%%%%%%%%%%%%%%%%%%%%%%%%%%%%%%%%%%
\subsection{Critical level}\label{newsec021}
%%%%%%%%%%%%%%%%%%%%%%%%%%%%%%%%%%%
Consider the quantum affine vertex algebra at the critical level $\Vcccrib=\VccNb$. 
   The following family of central elements for the  quantum   vertex algebra   $\Vcccrib$ was given by Molev and the author \cite[Prop. 3.5]{KM}.
\begin{pro}
All coefficients  of the series
$$
\phi_n(u)\coloneqq \tr_{1,\ldots ,n} \,A^{(n)}\ts T_{[n]}^+(u,u-h,\ldots ,u-(n-1)h) \ts D_1\ldots D_n \vac \in\Vcccrib[[u]]
$$
with $n=1,\ldots ,N$
belong to the center of the quantum   vertex algebra   $\Vcccrib $.
\end{pro}

Now consider the  quantum affine   algebra at the critical level 
  $\Uhcri =\Uh_{-N}$.
The next theorem goes back  to Frappat, Jing, Molev and Ragoucy \cite[Thm. 3.2]{FJMR}. Although it is originally given in terms of the $RLL$ realization of the quantum affine algebra,  we  formulate the theorem
 using  Ding's quantum current realization.
The direct proof in terms of Ding's realization   is carried out by arguing as in the proof of \cite[Thm. 2.14]{c11} and using Lemma \ref{profuzion}.

\begin{thm}
All coefficients  of the series
$$
\ell_n(z)\coloneqq  \tr_{1,\ldots ,n} \,A^{(n)}\ts\Lc_{[n]}(z_1,\ldots ,z_n)\big|_{z_1 = z,\ldots, z_n = ze^{-(n-1)h}}\big. \ts D_1\ldots D_n \in \Uhcri[[z^{\pm 1}]]
$$
with  $n=1,\ldots ,N$
belong to the center of the   algebra  $\Uhcri $. 
\end{thm}

Finally, let $W$ be any  restricted $\Uh$-module of level $-N$. Then the identities
\beq\label{cen1}
Y_W(\phi_n(0), z)=\ell_n (z) \quad\text{for } n=1,\ldots ,N
\eeq
hold for operators on $W$, where the map $Y_W(z)$ is given by \eqref{formula}. 
Recall Corollary \ref{maincor}.
By setting  $W=\Vcccri$ in \eqref{cen1} and then applying the resulting  equality on $1\in\Vcccri$ one recovers the  invariants of the vacuum module; see \cite[Corollary 3.3]{FJMR}.
\begin{kor}
All coefficients of the series
$$
\wvr{\ell}_n (z)\coloneqq Y_{\Vcccri} (\phi_n(0), z) 1 = \ell_n (z)1 \in \Vcccri  [[z]]
$$
with $n=1,\ldots ,N$
belong to the submodule of invariants $\z (\Vcccri )$.
\end{kor}

%34 pages. Main Theorem extended to sl_N. Subsect.3.4 and Sect.4 added. Other minor changes. Comments are welcome.

%\newpage
\section*{Acknowledgement}
The author would like to thank Naihuan Jing and Mirko Primc for stimulating discussions. 
The research reported in this paper was finalized during the author's visit to 
Max Planck Institute for Mathematics in Bonn. The author is grateful to the Institute for  its hospitality and financial support.
This work has been supported in part by Croatian Science Foundation under the project 8488.


\begin{thebibliography}{9}
\bibitem{BK}
B. Bakalov, V. G. Kac,
{\em Field algebras}, 
Int. Math. Res. Not. (2003), no. 3, 123--159;
\href{http://arxiv.org/abs/math/0204282}{arXiv:math/0204282 [math.QA]}.

\bibitem{B}
R. Borcherds  
{\em Vertex algebras, Kac--Moody algebras, and the Monster},
Proc. Natl. Acad. Sci. USA \textbf{83}  (1986) 3068--3071.

\bibitem{BJK}
M. Butorac, N. Jing, S. Ko\v{z}i\'{c},
{\em $h$-Adic quantum vertex algebras associated with rational $R$-matrix in types $B$, $C$ and $D$},  Lett. Math. Phys. \textbf{109} (2019), 2439--2471;
\href{https://arxiv.org/abs/1904.03771}{arXiv:1904.03771 [math.QA]}.

\bibitem{C}
I. V. Cherednik, 
{\em A  new  interpretation  of  Gelfand--Tzetlin bases},
Duke Math. J. \textbf{54} (1987), 563--577.

\bibitem{DGK}
A. De Sole, M. Gardini, V. G. Kac,
{\em On the structure of quantum vertex algebras},
J. Math. Phys. \textbf{61} (2020), 011701 (29pp);
\href{https://arxiv.org/abs/1906.05051}{arXiv:1906.05051 [math.QA]}.

\bibitem{D}
J. Ding, 
{\em Spinor Representations of $U_q(\hat{\gl} (n))$ and Quantum Boson-Fermion Correspondence},
Comm. Math. Phys. {\bf 200} (1999), 399--420;
\href{https://arxiv.org/abs/q-alg/9510014}{arXiv:q-alg/9510014}.

\bibitem{DF}
J. Ding, I. B. Frenkel,
{\em Isomorphism of two realizations of quantum affine algebra $U_q(\wht{\gl} (n))$},
Comm. Math. Phys. {\bf 156} (1993), 277--300.

\bibitem{DI}
J. Ding, K. Iohara,
{\em Generalization of Drinfeld Quantum Affine Algebras},
Lett. Math. Phys. {\bf 41} (1997), 181--193;
\href{https://arxiv.org/abs/q-alg/9608002}{arXiv:q-alg/9608002}.

\bibitem{EK3}
P. Etingof, D. Kazhdan,
{\em Quantization of Lie bialgebras, III}, Selecta Math. (N.S.) \textbf{4} (1998), 233--269;
\href{https://arxiv.org/abs/q-alg/9610030}{arXiv:q-alg/9610030}.

\bibitem{EK4}
P. Etingof and D. Kazhdan,
{\em Quantization of Lie bialgebras, IV},
Selecta Math. (N.S.) {\bf 6} (2000), 79--104;
\href{https://arxiv.org/abs/math/9801043}{arXiv:math/9801043 [math.QA]}.

\bibitem{EK}
P. Etingof, D. Kazhdan,
{\em Quantization of Lie bialgebras, V}, Selecta Math. (N.S.) \textbf{6} (2000), 105--130;
\href{http://arxiv.org/abs/math/9808121}{arXiv:math/9808121 [math.QA]}.

\bibitem{FRT}
N. Yu. Reshetikhin, L. A. Takhtadzhyan and L. D. Faddeev, 
{\em Quantization of Lie groups and Lie algebras}, 
Algebra i Analiz \textbf{1} (1989), no. 1, 178--206 (Russian); 
English transl., 
Leningrad Math. J. \textbf{1} (1990), no. 1, 193--225.

\bibitem{FJMR}
L. Frappat, N. Jing, A. Molev and E. Ragoucy,
{\em Higher Sugawara operators for the quantum affine
algebras of type $A$},
Comm. Math. Phys. {\bf 345} (2016), 631--657;
\href{https://arxiv.org/abs/1505.03667}{arXiv:1505.03667 [math.QA]}.

\bibitem{loop}
E. Frenkel,
{\em Langlands correspondence for loop groups}, 
Cambridge Studies in Advanced Mathematics, 103. Cambridge University Press, Cambridge, 2007.

\bibitem{FBZ}
E. Frenkel, D. Ben-Zvi,
{\em Vertex Algebras, Algebraic Curves}, 
Mathematical Surveys and Monographs, vol. 88, Second ed., American Mathematical Society, Providence, RI, 2004.

\bibitem{FR2}
E. Frenkel, N. Reshetikhin, 
{\em Towards deformed chiral algebras}, 
preprint
\href{http://arxiv.org/abs/q-alg/9706023}{arXiv:q-alg/9706023}.

\bibitem{FJ}
I. B. Frenkel, N. Jing, 
{\em Vertex representations of quantum affine algebras}, 
Proc. Natl. Acad. Sci. USA, \textbf{85} (1988), 9373--9377.

\bibitem{FLM}
I. Frenkel, J. Lepowsky, A. Meurman,
{\em Vertex operator algebras and the Monster},
Pure and Applied Mathematics, 134. Academic Press, Inc., Boston, MA, 1988.

\bibitem{FR}
I. B. Frenkel and N. Yu. Reshetikhin,
{\em Quantum affine algebras and holonomic difference equations},
Comm. Math. Phys. {\bf 146} (1992), 1--60.

\bibitem{FZ}
I. B. Frenkel and Y.-C. Zhu,
{\em Vertex operator algebras associated to representations of affine and Virasoro algebras},
Duke Math. J. \textbf{66} (1992), 123--168.

\bibitem{J}
M. Jimbo,
{\em A $q$-difference analogue of U(G) and the Yang--Baxter equation},
 Lett. Math. Phys. \textbf{10} (1985) 63--69.

\bibitem{JKMY}
N. Jing, S. Ko\v{z}i\'{c}, A. Molev, F. Yang,
{\em Center of the quantum affine vertex algebra in type $A$},
 J. Algebra \textbf{496} (2018), 138--186;
\href{https://arxiv.org/abs/1603.00237}{arXiv:1603.00237 [math.QA]}.

\bibitem{Kac} 
V. G. Kac, 
{\em Infinite-dimensional Lie algebras}, 
3rd ed., Cambridge University Press, Cambridge, 1990.

\bibitem{Kac2}
V. Kac,
{\em Vertex algebras for beginners}, 
University Lecture Series, 10. American Mathematical Society, Providence, RI, 1997.

\bibitem{Kas}
C. Kassel,
{\em Quantum Groups}, 
Graduate texts in mathematics; vol. \textbf{155}, Springer-Verlag, 1995.

\bibitem{KM}
S. Ko\v{z}i\'{c}, A. Molev, 
{\em Center of the quantum affine vertex algebra associated with trigonometric $R$-matrix}, 
J. Phys. A: Math. Theor. \textbf{50} (2017) 325201 (21pp); 
\href{https://arxiv.org/abs/1611.06700}{arXiv:1611.06700 [math.QA]}.

\bibitem{c11}
S. Ko\v{z}i\'{c}, 
{\em Quantum current algebras associated with rational $R$-matrix}, 
 Adv. Math. {\bf 351} (2019), 1072--1104;
\href{https://arxiv.org/abs/1801.03543}{arXiv:1801.03543 [math.QA]}.

\bibitem{LLi}
J. Lepowsky, H.-S. Li,
{\em Introduction to Vertex Operator Algebras and Their Representations},
Progress in Math., Vol. 227, Birkhauser, Boston, 2004.

\bibitem{LiG1}
H.-S. Li, 
{\em Axiomatic $G_1$-vertex algebras}, 
Commun. Contemp. Math. \textbf{5} (2003), 281--327;
\href{https://arxiv.org/abs/math/0204308}{arXiv:math/0204308 [math.QA]}.

\bibitem{Li}
H.-S. Li,
{\em $\hbar$-adic quantum vertex algebras and their modules},
Comm. Math. Phys. {\bf 296} (2010), 475--523;
\href{http://arxiv.org/abs/0812.3156}{arXiv:0812.3156 [math.QA]}.

\bibitem{Li1}
H.-S. Li,
{\em $\phi$-Coordinated Quasi-Modules for Quantum Vertex Algebras},
Comm. Math. Phys. {\bf 308} (2011), 703--741;
\href{https://arxiv.org/abs/0906.2710}{arXiv:0906.2710 [math.QA]}.

\bibitem{LTW}
H.-S. Li, S. Tan, Q. Wang,
{\em Ding--Iohara algebras and quantum vertex algebras},
J. Algebra {\bf 511} (2018), 182--214;
\href{https://arxiv.org/abs/1706.03636}{arXiv:1706.03636 [math.QA]}.

\bibitem{Lian}
B.-H. Lian,
{\em On the classification of simple vertex operator algebras},
Comm. Math. Phys. {\bf 163} (1994), 307--357.


\bibitem{PS}
J. H. H. Perk, C. L. Schultz,
{\em New families of commuting transfer matrices in $q$-state vertex models},
Phys. Lett. A \textbf{84} (1981), 407--410.

\bibitem{RS}
N. Yu. Reshetikhin, M. A. Semenov-Tian-Shansky,
{\em Central extensions of quantum current groups}, 
Lett. Math. Phys. {\bf 19} (1990), 133--142.

\end{thebibliography}
\end{document}